\documentclass[11pt,reqno]{amsart}
\usepackage{float}
\usepackage{fullpage}
\usepackage{amsmath,amsthm,verbatim,amssymb,amsfonts,amscd, graphicx,bbm,bm}
\usepackage{graphics}
\usepackage{mathtools}
\usepackage{mathrsfs}
\usepackage{esint}
\usepackage{color}
\usepackage{tikz,tikz-cd}
\usepackage[all,cmtip]{xy}
\usetikzlibrary{decorations.pathreplacing,arrows.meta,calc,patterns,angles,quotes,decorations.markings}

\usepackage[pdfstartview=FitH, bookmarksnumbered=true,bookmarksopen=true, colorlinks=true, pdfborder={0 0 1}, citecolor=blue, linkcolor=blue,urlcolor=blue]{hyperref}

\newtheorem{theorem}{Theorem}
\newtheorem{proposition}[theorem]{Proposition}
\newtheorem{lemma}[theorem]{Lemma}
\newtheorem{corollary}[theorem]{Corollary}

\theoremstyle{definition}
\newtheorem{remark}{Remark}

\newtheorem{definition}[theorem]{Definition}

\newcommand{\cref}[1]{Corollary~\ref{c.#1}}

\usepackage{xcolor}

\numberwithin{equation}{section}
\numberwithin{theorem}{section}

\newcommand{\R}{\mathbb{R}}

\newcommand{\eps}{\varepsilon}

\newcommand{\Z}{\mathbb{Z}}

\newcommand{\uinc}{u^{\rm in}}

\title{Homogenization of the Scattered Wave and Scattering Resonances for Periodic High-Contrast Subwavelength Resonators}

\author{Yuxin Du}
\address[Y. Du]{Qiuzhen College, Tsinghua University, Beijing 100084}
\email{duyx23@mails.tsinghua.edu.cn}

\author{Xin Fu}
\address[X. Fu]{}
\email{1550862645cf@gmail.com}

\author{Wenjia Jing}
\address[W. Jing]{Yau Mathematical Sciences Center, Tsinghua University, Beijing 100084 and Beijing Institute of Mathematical Sciences and Applications, Beijing 101408, P.R. China}
\email{wjjing@tsinghua.edu.cn}

\date{\today}

\begin{document}

\begin{abstract}
    We study time-harmonic scattering by a periodic array of penetrable, high-contrast obstacles with small period, confined to a bounded Lipschitz domain. The strong contrast between the obstacles and the background induces subwavelength resonances. We derive a frequency-dependent effective model in the vanishing-period limit and prove quantitative convergence of the heterogeneous scattered wave to the effective scattered wave. We also identify the limiting set of scattering resonances and establish convergence rates. Finally, we establish convergence rates for the far-field pattern of the heterogeneous problem to that of the effective model.

    \smallskip
    

    \smallskip
    
    \noindent{\bf Mathematics subject classification (MSC 2020)}: 35B27, 35B34, 35J70, 35P25, 74J20
\end{abstract}

\maketitle

\tableofcontents

\section{Introduction}

Over the past decades, metamaterials have shown tremendous potential across many areas of science and technology. The great interest in metamaterials arises from their ability to manipulate waves in ways that are unattainable with natural materials. The core idea is to replace the molecules with man-made structures, viewed as “artificial atoms” on a scale much less than the relevant wavelength. In this way, the metamaterial can be described using a small number of effective parameters. In 1968, Veselago \cite{Viktor} theoretically considered materials with both negative permittivity and permeability, and predicted several counterintuitive properties, such as negative refraction. However, at that time, no actual examples of such materials existed. In 1999, J. Pendry \textit{et al.} \cite{798002} proposed the first practical route to metamaterials by designing structures at a subwavelength scale, and built the first electromagnetic metamaterial by using the blocks of split-ring resonators to achieve negative magnetic permeability. Following this breakthrough, Smith \emph{et al.} realized electromagnetic materials exhibiting both negative permittivity and permeability \cite{PhysRevLett.84.4184}. Since then, the field has expanded rapidly from electromagnetics to acoustics, optics, and mechanics, with a great amount of applications, including negative-index media \cite{smith2004metamaterials}, cloaking \cite{cai2007optical}, subwavelength focusing \cite{guenneau2007acoustic}, band-gap engineering \cite{sugino2016mechanism}, tunability \cite{lapine2009structural} and topological insulators \cite{khanikaev2013photonic}. 

Metamaterials can exhibit unusual wave phenomena when their elementary building blocks resonate at wavelengths significantly larger than their size. These elements, known as subwavelength resonators, are engineered by introducing a singular feature through either extreme physical parameter or a carefully designed singular geometry \cite{feppon2023homogenization, schweizer2017resonance}. Typical realizations include split-ring resonators \cite{farhat2009negative, lipton2018effective}, Helmholtz resonators \cite{hu2008homogenization, tachet2025homogenized}, and high-contrast resonators \cite{davies2024landscape, dang2025homogenization}. These resonators are often arranged periodically to simplify analysis and fabrication.

In this work we consider time-harmonic scattering by a metamaterial whose building blocks are high-contrast subwavelength resonators. The composite occupies a bounded Lipschitz set $\Omega\subset\mathbb{R}^d$, inside which lies a periodic array of inclusions with period $\varepsilon$. The ratio of the relevant material parameters in the inclusions to those in the background scales like $\varepsilon^2$, producing local resonances in the subwavelength regime.

A simplified mathematical model is the scalar Helmholtz equation in divergence form. For a wavenumber $k>0$ and incident field $\uinc$, the total field $u_\varepsilon$ solves
\begin{equation}
    (\nabla\cdot A_{\varepsilon} \nabla +k^2)u_{\varepsilon} =0 \qquad \mathrm{in}\ \mathbb{R}^d,
\end{equation}
together with the Sommerfeld radiation condition at infinity. Here $A_\varepsilon\equiv I$ in the background, while $A_\varepsilon$ is of order $\varepsilon^2$ in the inclusions, encoding the high-contrast perturbation responsible for local resonance.

From a mathematical viewpoint, this is a homogenization problem for non-uniformly elliptic operators of divergence form on unbounded domains. Prior work in electromagnetics and elasticity, notably by Bouchitt\'{e}, Felbacq, and \'{A}vila \textit{et al.} \cite{bouchitte2004homogenization,felbacq2005theory,felbacq2005negative,avila2008multiscale}, established two-scale limits to frequency-dependent effective models that may exhibit negative effective permeability and associated band-gap/dispersive phenomena. However, to the best of our knowledge, most existing results are qualitative, or quantitative only under restrictive assumptions (such as dilute arrays or full-space periodicity; see the discussion in Section \ref{previouswork}). In particular, convergence rates for scattering problems with bounded metamaterials remain largely open, limiting practical applications and numerical computation.

This paper aims to address these gaps. We derive a frequency-dependent effective model as $\varepsilon \rightarrow 0$ and provide quantitative error bounds. The analysis combines two-scale expansions with corrector estimates that account for boundary layers near the boundary of the composite. A key point in the proof is characterizing the limiting distribution of scattering resonances for the operator $\nabla \cdot A_{\varepsilon } \nabla$ as $\varepsilon \rightarrow 0$. We also establish a convergence rate of the scattering resonances towards the limiting distribution, which is of independent interest in scattering theory. As a corollary of the above results, we obtain convergence rates for the far-field pattern of the highly oscillatory scattering problem toward the far-field pattern of the effective model. 

These results quantify the accuracy of the effective model and identify the frequency intervals (band gaps) in which the effective parameters are negative. The techniques developed here also suggest a path toward rigorous error control for numerical simulations of resonant metamaterials, where accuracy near resonance is critical.

The remainder of the paper is organized as follows. In Section \ref{problemset}, we introduce the problem set-up and formally derive the effective model. We then present the main results. The relation to previous works is also discussed. In Section \ref{secconrateprob}, assuming suitable boundary layer estimates and uniform resolvent estimates, we establish the \emph{a priori} $L^2$ convergence rate and the weighted $H^1$ convergence rate. In Section \ref{secbound}, we establish boundary layer estimates using the method of layer potentials. In Section \ref{secresout}, we construct the meromorphic continuation of the resolvent and study the limiting distribution of the scattering resonances. As a result, we obtain the uniform resolvent estimates. In Section \ref{secoptmail}, we improve the $L^2$ convergence rate to optimal under certain smooth assumptions on the boundary of the composite, and we obtain the convergence rate for the far-field pattern. In Section \ref{secrema}, we conclude this paper and outline potential future work. Appendix \ref{secunfold} is devoted to the periodic unfolding method used in Section \ref{subsecrateout}. Appendix \ref{appendixA} collects results on the scattering resonances of the effective model, used in Section \ref{secresout}.

\section{Problem set-up and main results}\label{problemset}

\subsection{Geometry set-up and the scattering problem}

\begin{figure}[h]
    \centering
    \includegraphics[width=0.5\linewidth]{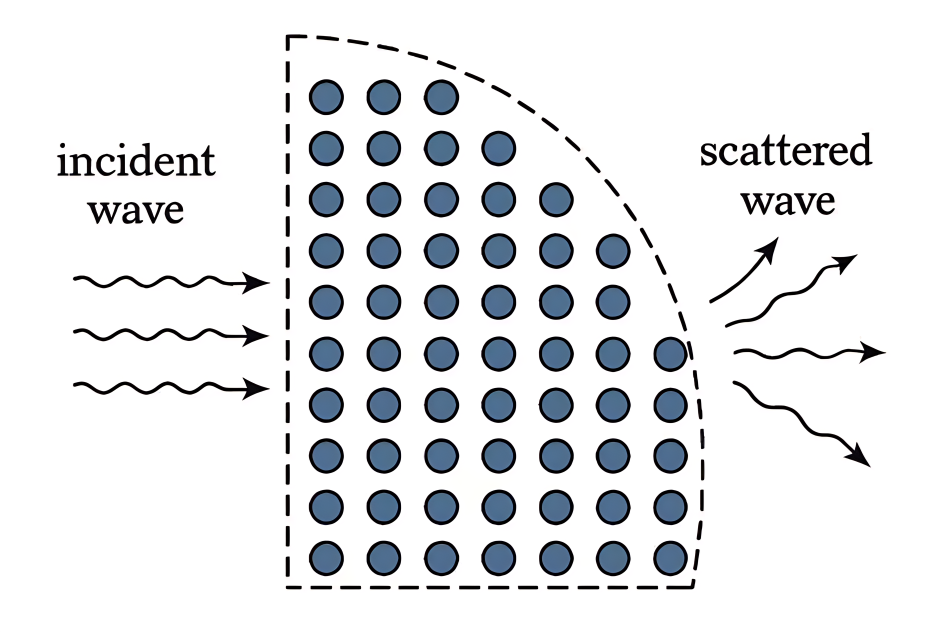}
    \caption{An illustration for the scattering by a periodic array of small inclusions.}
    \label{figscatter}
\end{figure}

Let $Y=(0, 1)^d$ denote the unit cell of $\R^d$, $d\ge 2$. Let $D$ be a $C^{1,\alpha}$ domain with connected boundary compactly contained in $Y$ for some $\alpha>0$. Let $\Omega$ be a bounded connected Lipschitz domain in $\mathbb{R}^d$. For any $\varepsilon \in (0,1)$, we define
\begin{equation}\label{jepsi}
    J_{\varepsilon}:=\{ \mathbf{m} \in \mathbb{Z}^d : \varepsilon (Y+\mathbf{m}) \subset \Omega \},
\end{equation}
and
\begin{equation}\label{def0D}
    D_{\varepsilon}:= \bigcup_{\mathbf{m}\in J_{\varepsilon}} \varepsilon (D+\mathbf{m}) .
\end{equation}
Roughly speaking, $D_{\varepsilon}$ consists of periodically distributed obstacles of size $\varepsilon$ contained in $\Omega$ and away from its boundary; see Figure \ref{figscatter}.

Fix a wave number $k>0$. Let $\uinc$ be a free incident wave solving
\begin{equation}
    (\Delta  + k^2) \uinc =0 \qquad \mathrm{in} \ \mathbb{R}^d.
\end{equation}
A typical choice of $\uinc$ is the plane wave $\mathrm{e}^{\mathrm{i}k\omega \cdot x}$ propagating along the direction $\omega \in S^{d-1}$. 

Let 
\begin{equation}
    L_{\varepsilon} := \nabla\cdot A_{\varepsilon}\nabla
\end{equation}
be the second-order elliptic operator of divergence form with the coefficient
\begin{equation}\label{defAr}
    A_{\varepsilon}: = \varepsilon^2 \mathbbm{1}_{D_{\varepsilon}} + \mathbbm{1}_{\mathbb{R}^d\setminus \overline{D_{\varepsilon}}} ,
\end{equation}
where $\mathbbm{1}_E$ denotes the indicator function of a set $E\subset \mathbb{R}^d$. 

A function $u$ is said to satisfy the Sommerfeld radiation condition at wavenumber $k$ (written $u\in\mathrm{SRC}(k)$) if
\begin{equation}\label{src}
    (\partial_{|x|}  - \mathrm{i} k) u(x) =\mathcal{O}(|x|^{-(d+1)/2}) \qquad \mathrm{as}\ x\rightarrow \infty,
\end{equation}
where $\partial_{|x|} := \frac{x}{|x|}\cdot\nabla$ denotes the radial derivative.

We consider the scattering problem
\begin{equation}\label{maineq}
\left\{\begin{aligned}
    & (L_\varepsilon  + k^2) u_\varepsilon = 0 \quad \mathrm{in}\ \mathbb{R}^d,\\
    & u_{\varepsilon} -\uinc\in \mathrm{SRC}(k).
\end{aligned}\right.
\end{equation}
The main goal of this paper is to study the limiting behavior of $u_{\varepsilon}$ as $\varepsilon \rightarrow 0$.

\subsection{Two-scale expansion and homogenized scattering problem}

In this section, we formally derive the asymptotic expansion of $u_{\varepsilon}$ by using the method of two-scale expansion \cite{bensoussan_asymptotic_1978}, and, on this basis, introduce the homogenized scattering problem.
 
Denote the flux $ v_\varepsilon := A_\varepsilon \nabla u_\varepsilon $. Then, the equation in \eqref{maineq} can be rewritten as
\begin{equation}\label{ansatzeq2}
    \nabla \cdot v_\varepsilon + k^2 u_\varepsilon=0  \qquad \mathrm{in} \ \mathbb{R}^d,
\end{equation}
with the transmission conditions
\begin{equation}\label{ansatz_trans_cond}
   \begin{aligned}
       & u_\varepsilon|_+ = u_\varepsilon|_-  \\
       &\mathbf{n}\cdot v_\varepsilon|_+ = \mathbf{n}\cdot v_\varepsilon|_- 
   \end{aligned}  \qquad \mathrm{on}\ \partial D_{\varepsilon},
\end{equation}
where $\mathbf{n}$ denotes the unit normal pointing outward from $D_\varepsilon$, and $|_{\pm}$ denotes the traces taken from the exterior and interior of $\partial D_\varepsilon$.

We adopt the following two-scale ansatz:
\begin{equation}\label{ansatz_inside}
\begin{aligned}
    & u_\varepsilon(x) = u_0(x,x/\varepsilon) + \varepsilon u_1(x,x/\varepsilon)  + \cdots \\
    & v_\varepsilon (x)= \eps^{-1} v_{-1}(x,x/\eps) + v_0(x,x/\varepsilon) + \varepsilon v_1(x,x/\varepsilon) + \cdots
\end{aligned} \qquad \mathrm{in} \ \Omega,
\end{equation}
where, for each $j\geq -1$, the functions $u_j$ and $v_j$ are $Y$-periodic in the fast variable $y:=x/\varepsilon$ and depend on the slow variable $x$. Substituting \eqref{ansatz_inside} into \eqref{ansatzeq2}, using the chain rule
\begin{equation}
    \nabla = \nabla_x +\varepsilon^{-1 }\nabla_y ,
\end{equation}
and equating coefficients of like powers of $\varepsilon$, yield
\begin{align}
    & \nabla_y \cdot v_{-1}=0 & \mathrm{in} \ \Omega \times Y, \label{Yv-1}\\
    & \nabla_y \cdot v_0 + \nabla_x \cdot v_{-1} = 0, & \mathrm{in} \ \Omega \times Y,  \label{Yv0}\\
    & \nabla_y \cdot v_1 + \nabla_x \cdot v_0 + k^2 u_0 = 0 & \mathrm{in} \ \Omega \times Y.\label{Yv01}
\end{align}
Substituting \eqref{ansatz_inside} into the constitutive relation $v_{\varepsilon} = A_{\varepsilon} \nabla u_{\varepsilon}$, using the definition \eqref{defAr} of $A_{\varepsilon}$ and equating coefficients of like powers of $\varepsilon$, we get 
\begin{align}
& \nabla_y u_0 = v_{-1} && \mathrm{in} \ \Omega \times (Y\setminus \overline{D}) ,\label{u0outside}\\
    & \nabla_y u_1 + \nabla_x u_0  = v_0 && \mathrm{in} \ \Omega \times (Y\setminus \overline{D}), \label{u01outside} \\
     & \nabla_y u_2 + \nabla_x u_1  = v_1 && \mathrm{in} \ \Omega \times (Y\setminus \overline{D}), \label{u12outside} \\
    & 0= v_{-1} &&\mathrm{in} \ \Omega \times D,\label{v-1inside}\\
    & 0= v_0 &&\mathrm{in} \ \Omega \times D,\label{v0inside}\\
     &  \nabla_y u_0 = v_1 &&\mathrm{in} \ \Omega \times D.\label{u0inside}
\end{align}
From \eqref{Yv-1}\eqref{u0outside} and \eqref{v-1inside} we get $v_{-1} = 0$ in $Y$ and hence can delete $v_{-1}$ from the ansatz \eqref{ansatz_inside}. Substituting \eqref{ansatz_inside} into the transmission conditions \eqref{ansatz_trans_cond}, equating coefficients of like powers of $\varepsilon$, and using \eqref{u01outside} and \eqref{v0inside}, we obtain
\begin{align}
    &  u_0|_+ = u_0|_- & \mathrm{on} \ \Omega \times \partial D, \label{transcond1}  \\
    &  \mathbf{n}(y) \cdot \nabla_x u_0|_+ + \mathbf{n} (y)\cdot \nabla_y u_1 |_+ = 0 & \mathrm{on} \ \Omega \times \partial D . \label{transcond3}
\end{align}

Next we determine $u_0$ (and subsequently $u_1$) from \eqref{Yv-1}–\eqref{transcond3}. From \eqref{u0outside}, there exists a function $\widehat{u}_0$ on $\Omega$ such that
\begin{equation}
    u_0(x,y) = \widehat{u}_0(x) \qquad \mathrm{for} \ (x,y)\in\Omega \times (Y\setminus \overline{D}).
\end{equation}
Substituting \eqref{v0inside}-\eqref{u0inside} into \eqref{Yv01} yields
\begin{equation}
    \Delta_y u_0 +k^2u_0 =0 \qquad \mathrm{in}\  \Omega \times D.
\end{equation}
Using the transmission condition \eqref{transcond1}, we conclude that $u_0$ factorizes as
\begin{equation}
    u_0(x,y) = \Lambda(y) \widehat{u}_0(x) \qquad \mathrm{for}\ (x,y)\in \Omega \times Y,
\end{equation}
where $\Lambda$ satisfies the equation 
\begin{equation}\label{eq_lambda}
    \left\{
    \begin{aligned}
        &(\Delta  + k^2) \Lambda = 0 && \mathrm{in} \ D , \\
        &\Lambda \equiv 1 && \mathrm{in} \ Y \setminus D.
    \end{aligned}
    \right. 
\end{equation}
To ensure that $\Lambda$ is well-defined, we need to require that $k^2$ not coincide with any Dirichlet eigenvalues of $-\Delta$ on $D$ with an eigenfunction of nonzero mean.

Applying $\mathrm{div}_y$ to \eqref{u01outside}, using \eqref{Yv0} and noting that, in $\Omega \times (Y\setminus \overline{D})$, $u_0$ depends only on $x$, we get
\begin{equation}
    \Delta_y u_1=0 \qquad \mathrm{in}\ \Omega \times (Y\setminus \overline{D} ).
\end{equation}
By the transmission condition \eqref{transcond3}, $u_1$ admits the representation
\begin{equation}
    u_1(x,y)=\nabla \widehat{u}_0(x) \cdot \chi(y)\  \qquad \mathrm{for}\ (x,y)\in \Omega \times (Y\setminus \overline{D}),
\end{equation}
where $\chi=(\chi_1,\dots,\chi_d)^\top$ is $Y$–periodic with zero mean and, for each $j=1,\dots,d$,
\begin{equation}\label{def0chi}
    \left\{
    \begin{aligned}
        &\Delta\chi_j=0 && \mathrm{in}\ Y\setminus \partial D ,\\
        &\mathbf{n}\cdot\nabla(\chi_j+y_j) |_+=0 && \mathrm{on} \ \partial D,
    \end{aligned}
    \right. 
\end{equation}
where $y_j$ is the $j$-th component of $y$.

Integrating \eqref{Yv01} over $Y$ yields
\begin{equation}
    \nabla_x \cdot \left(\int_Yv_0(x,y)\,dy \right) +  k^2 \left(\int_Y \Lambda(y)\,dy\right) \widehat{u}_0(x) = 0 \qquad  \mathrm{in} \ \Omega .
\end{equation}
Using \eqref{u01outside}–\eqref{v0inside} we compute
\begin{equation}
    \int_Yv_0(x,y)\,dy   = \left(\int_{Y\setminus \overline{D}} (I + \nabla \chi)(y) \,dy \right) \nabla \widehat{u}_0(x), 
\end{equation}
where $I$ denotes the $d\times d$ identity matrix. Therefore $\widehat{u}_0$ solves
\begin{equation}
    \nabla\cdot \left( \int_{Y\setminus \overline{D}} (I + \nabla \chi)(y) \,dy\right)\nabla \widehat{u}_0 + k^2 \left(  \int_Y \Lambda(y) \,dy  \right)\widehat{u}_0=0  \qquad \mathrm{in} \ \Omega.
\end{equation}

We now define the homogenized operator
\begin{equation}
    L_0 : = \nabla \cdot A_0(x) \nabla,
\end{equation}
where the (piecewise constant) coefficient field $A_0$ is defined by
\begin{equation}\label{homogenizedcoe}
    A_0(x) := \left(\int_{Y\setminus \overline{D}}  (I +\nabla\chi)(y) \,dy \right) \mathbbm{1}_{\Omega} (x)+  \mathbbm{1}_{\mathbb{R}^d\setminus \overline{\Omega}} (x)
\end{equation}
The homogenized scattering problem reads
\begin{equation}\label{equ0}
    \left\{ 
        \begin{aligned}
            &(L_0  + k^2 \mu^k_0 ) \widehat{u}_0 = 0 && \mathrm{in} \ \R^d, \\
            & \widehat{u}_0 - \uinc \in \mathrm{SRC}(k).
        \end{aligned}
        \right.
\end{equation} 
Here, the (piecewise constant) coefficient function $\mu^k_0$ is defined by
\begin{equation}
   \mu^k_0 (x):= \left( \int_Y \Lambda(y)\,dy \right) \mathbbm{1}_{\Omega} (x)+  \mathbbm{1}_{\mathbb{R}^d\setminus \overline{\Omega}}(x).
\end{equation}
For brevity, when no confusion can arise, we also write
\begin{equation}\label{A0mu0k}
    A_0 = \int_{Y\setminus D}  (I +\nabla\chi)(y) \,dy \qquad \mathrm{and} \qquad \mu_0^k=\int_Y \Lambda(y)\,dy .
\end{equation}

\subsection{Main results} We introduce some notations. Denote
\begin{equation}
    \Sigma_D:= \{ z\in \mathbb{C}: z^2 \textrm{ is a Dirichlet eigenvalue of } -\Delta \textrm{ on } D\}.
\end{equation}
For $z \in \Sigma_D$, we write $z \in \Sigma_{D,1}$ if $z^2$ is a Dirichlet
eigenvalue of $-\Delta$ on $D$ with an eigenfunction of nonzero mean. We then proceed to decompose $\Sigma_D$ into 
\begin{equation}
    \Sigma_D = \Sigma_{D,0} \cup \Sigma_{D,1},
\end{equation}
where $\Sigma_{D,0}$ consists of $z \in \Sigma_D$ such that all Dirichlet eigenfunctions associated with the eigenvalue $z^2$ are mean-zero.

Throughout the paper, we use the shorthand $f^{\varepsilon}(x):= f(x/\varepsilon)$ for any function $f$. We  denote by $B_r$ the open ball of radius $r$ centered at the origin in $\mathbb{R}^d$. Define
\begin{equation}
    \Lambda_{\varepsilon}(x) : = \Lambda(x/\varepsilon) \mathbbm{1}_{\Omega}(x) + \mathbbm{1}_{\mathbb{R}^d \setminus \Omega} (x), \qquad \forall x\in \mathbb{R}^d.
\end{equation}

We now state the main results of the paper.

\begin{theorem}\label{mainresult1}
    Fix $k\in (0,\infty) \setminus \Sigma_D$. Let $ u_\varepsilon $ be the solution of the scattering problem \eqref{maineq}, and let $\widehat{u}_0$ be the solution of the homogenized scattering problem \eqref{equ0}. Then, for any $r>0$ such that $\overline{\Omega}\subset B_r$, we have
\begin{equation}\label{L2formulaLam}
        \| u_\varepsilon - \Lambda_\varepsilon \widehat{u}_0  \|_{L^2(B_r)}  \leq C\varepsilon^{1/2} \|\uinc\|_{H^1(B_{r+1})}  ,
\end{equation}
and
\begin{equation}\label{L2formulaLam11}
        \| u_\varepsilon - \widehat{u}_0 -\varepsilon \chi^{\varepsilon}   \cdot\eta_{\varepsilon} S_\varepsilon ( \nabla \widehat{u}_0 )  \|_{H^1(B_r\setminus \overline{D_{\varepsilon}})}  \leq C\varepsilon^{1/2} \|\uinc\|_{H^1(B_{r+1})}  ,
\end{equation}
where $\eta_{\varepsilon}$ is the cut-off function defined in \eqref{cut-off_def}, $S_{\varepsilon}$ is the smoothing operator defined in \eqref{epsi_smoother}, and $C>0$ is a constant that depends only on $d,k,r,D$ and $\Omega$.
\end{theorem}

Several remarks are in order.

\begin{remark}\label{rem1rk}
    Since $\Lambda_{\varepsilon} \equiv 1$ on $\mathbb{R}^d\setminus \overline{D_{\varepsilon}}$, and $\mathrm{supp}\,\eta_{\varepsilon} \subset \Omega$, the estimates \eqref{L2formulaLam}-\eqref{L2formulaLam11} imply 
    \begin{equation}
        \| u_\varepsilon - \widehat{u}_0  \|_{L^2(B_r\setminus \overline{D_{\varepsilon}} )} + \| u_\varepsilon - \widehat{u}_0  \|_{H^1(B_r\setminus \overline{\Omega} )}  \leq C\varepsilon^{1/2} \|\uinc\|_{H^1(B_{r+1})},
    \end{equation}
    i.e., $\widehat{u}_0$ is an effective $L^2$ approximation to $u_\varepsilon$ away from the obstacles $D_\varepsilon$, and an effective $H^1$ approximation to $u_\varepsilon$ away from $\Omega$, both with error $\mathcal{O}(\varepsilon^{1/2})$. Moreover, since 
    \begin{equation}
        (\Delta+k^2)(u_{\varepsilon} - \widehat{u}_0) =0 \qquad \mathrm{in} \ \mathbb{R}^d\setminus \overline{\Omega},
    \end{equation}
    interior elliptic regularity implies that $\widehat{u}_0$ is an effective $C^{\infty}$ approximation to $u_\varepsilon$ in any bounded domain $U$ that is compactly contained in $\mathbb{R}^d\setminus \overline{\Omega}$. More precisely, for any $\alpha>0$ and $U'$ compactly contained in $U$, there exists a constant $C>0$ that depends only on $d,k,\alpha,U',U,D$ and $\Omega$ such that 
    \begin{equation}
        \| u_\varepsilon - \widehat{u}_0  \|_{C^\alpha(U')} \leq C\varepsilon^{1/2} \| \uinc \|_{H^1(U)}.
    \end{equation}
\end{remark}

\begin{remark}
    The convergence rate $\mathcal{O}(\varepsilon^{1/2})$ comes from boundary layer estimates near $\partial \Omega$ for the homogenized scattering problem \eqref{equ0}; see Theorem \ref{blethm} below. For general Lipschitz domains $\Omega$, this rate is not expected to improve. Nevertheless, if $\Omega$ is smoother (e.g., $C^{1,1}$), a duality argument yields the optimal rate $\mathcal{O}(\varepsilon)$; see Theorem \ref{mainresult2} below. Alternatively, one can improve the convergence rate by constructing higher-order boundary correctors; see \cite{cakoni2015homogenization,cakoni2016homogenization} for details.
\end{remark}

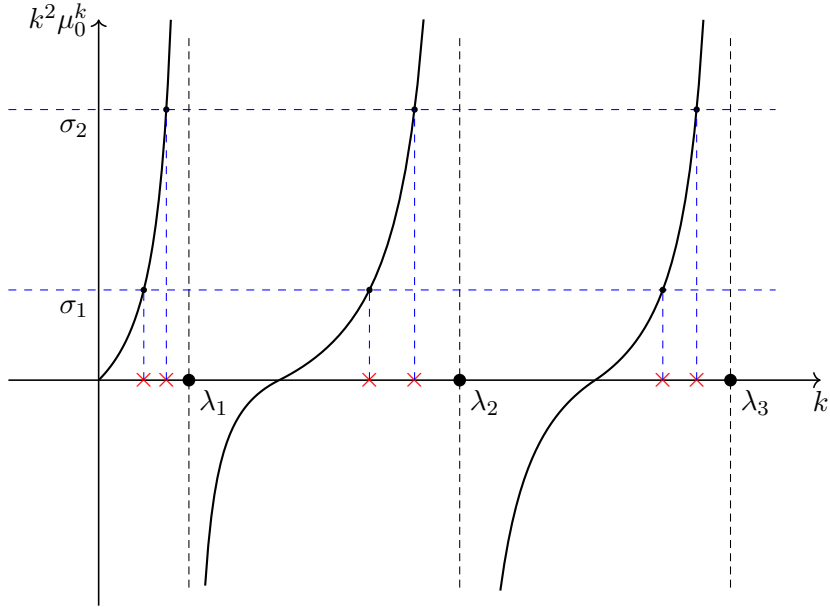
\begin{figure}[h]
\begin{center}
    \begin{tikzpicture}[scale=1.2]
    \draw[->, line width = 0.6pt] (-1,0) -- ( 8, 0) node[below]{$k$};
    \draw[->, line width = 0.6pt] (0, -2.5 ) -- ( 0, 4) node[left]{$k^2 \mu_0^k$};
    \foreach \x in{1,4,7}
          {\fill(\x,0)circle(2pt);}
    \node[below right] at ( 1,0 ){$\lambda_1$};
    \node[below right] at ( 4,0 ){$\lambda_2$};
    \node[below right] at ( 7,0 ){$\lambda_3$};
    \draw[densely dashed] (1, -2.3 ) -- ( 1, 3.8);
    \draw[densely dashed] (4, -2.3 ) -- ( 4, 3.8);
    \draw[densely dashed] (7, -2.3 ) -- ( 7, 3.8);
    \draw[blue, dashed, domain=-1:7.5]plot(\x,1);
    \node[below left] at (0,1){$\sigma_1$};
    \draw[blue, dashed, domain=-1:7.5]plot(\x,3); 
    \node[below left]at (0,3) {$\sigma_2$};
    \draw[line width = 0.8pt, domain=0:0.8]plot(\x,{1/(1-\x)-1});
    \draw[line width = 0.8pt,domain=2:3.6]plot(\x,{2/(4 -\x)-1});
    \draw[line width = 0.8pt,domain=1.18:2]plot(\x,{1/(2-2*\x)+1/2});
    \draw[line width = 0.8pt,domain=5.5:6.7]plot(\x,{3/2/(7-\x)-1});
    \draw[line width = 0.8pt,domain=4.45:5.5]plot(\x,{3/2/(4-\x)+1});
    \fill(1/2,1)circle(1pt);\draw[red] (1/2,0) node {$\times$};
    \draw[blue, dashed] (1/2, 0) -- ( 1/2, 1);
    \fill(3/4,3)circle(1pt);\draw[red] (3/4,0) node {$\times$};
    \draw[blue, dashed] (3/4, 0) -- ( 3/4, 3);
    \fill(3,1)circle(1pt);\draw[red] (3,0) node {$\times$};
    \draw[blue, dashed] (3, 0) -- ( 3, 1);
    \fill(7/2,3)circle(1pt);\draw[red] (7/2,0) node {$\times$};
    \draw[blue, dashed] (7/2, 0) -- ( 7/2, 3);
    \fill(25/4,1)circle(1pt);\draw[red] (25/4,0) node {$\times$};
    \draw[blue, dashed] (25/4, 0) -- ( 25/4, 1);
    \fill(53/8,3)circle(1pt);\draw[red] (53/8,0) node {$\times$};
    \draw[blue, dashed] (53/8, 0) -- ( 53/8, 3);
    \end{tikzpicture}
\end{center}
\caption{The curves represent the function $k\mapsto k^2\mu_0^k$. On a bounded domain $\Omega$, the limiting spectrum of $-L_{\varepsilon}$ with Dirichlet boundary conditions consists of the spectrum of $-\Delta_D$ together with the red crosses, where each $\lambda_j$ is an element of $\Sigma_{D,1}$ and each $\sigma_j$ is a Dirichlet eigenvalue of $-L_0$ on $\Omega$. It is clear that $\Sigma_{D,1}$ is the essential spectrum, since the red crosses in each interval of $(\lambda_{j-1},\lambda_j)$ converge to each $\lambda_j$.}
\label{figmu0k}
\end{figure}

\begin{figure}[h]
    \begin{center}
        \includegraphics[width=0.7\linewidth]{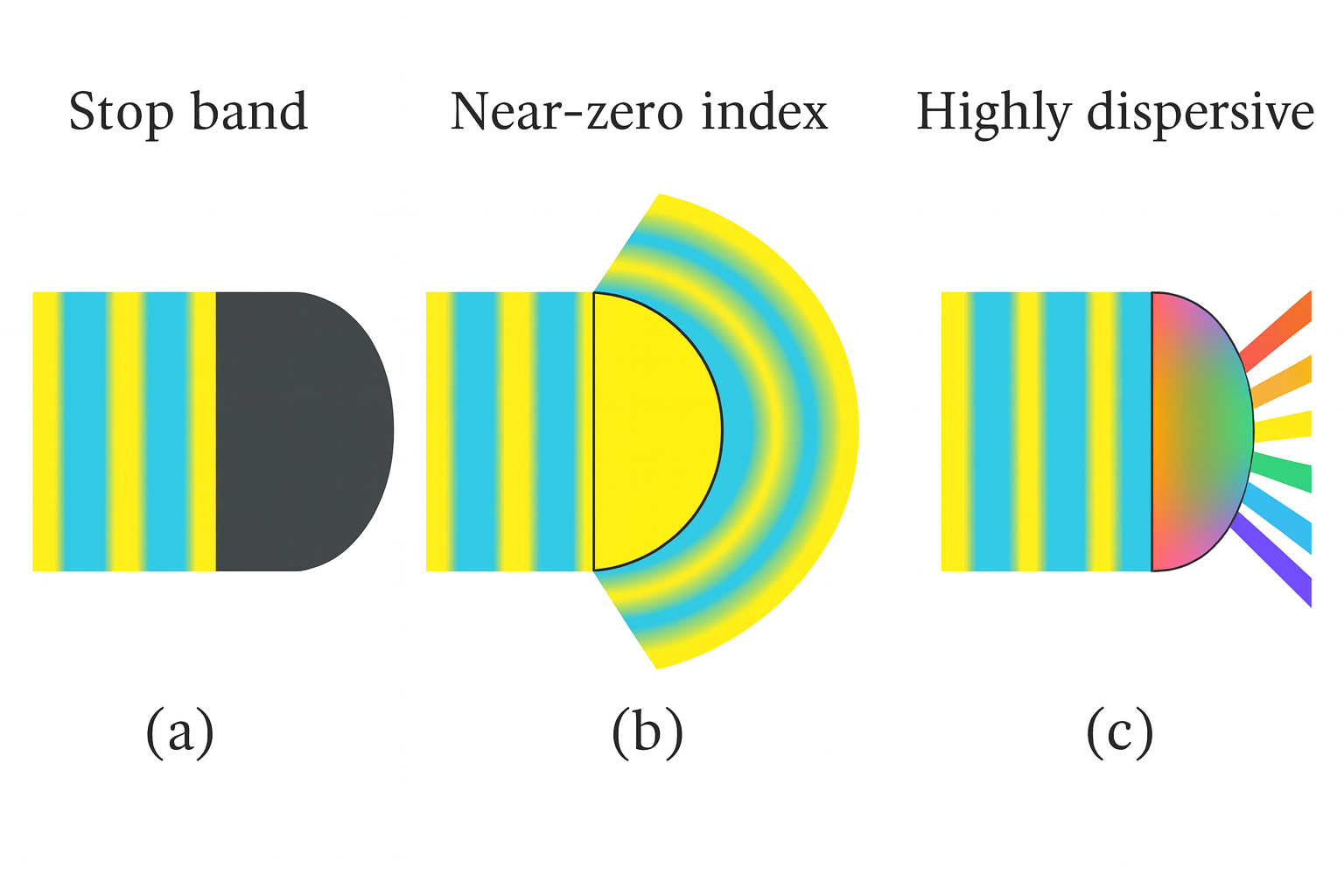}
    \end{center}
    \caption{An illustration of light propagation: (a) $\mu_0^k<0$, light cannot transmit through the medium; (b) $\mu_0^k \approx 0$, light propagates with a uniform phase inside the medium; (c) $\mu_0^k \rightarrow \infty$, the medium looks iridescent as light passes through, due to strong dispersive effects.}
    \label{figband}
\end{figure}

\begin{remark}
By the definitions of $A_0$ and $\mu_0^k$ in \eqref{A0mu0k}, the homogenized tensor $A_0$ is
symmetric and positive-definite, while $\mu_0^k$ changes sign as a function of the frequency $k$.
In particular,
\begin{equation}\label{eq:mu0k_def}
  \mu_0^k
  = 1
  + k^2 \sum_{\lambda_j \in \Sigma_{D,1}}
    \frac{\big(\int_D \varphi_j(y)\,dy\big)^{2}}{\lambda_j^{2}-k^{2}},
\end{equation}
where $\{\varphi_j\}$ are the orthonormal Dirichlet eigenfunctions of eigenvalues $\{\lambda^2_j\}$. A direct computation yields $\frac{d\mu_0^k}{dk}>0$, so $\mu_0^k$ is strictly increasing on every interval $(\lambda_j,\lambda_{j+1})$. Moreover,
\begin{equation}\label{eq:limits}
  \lim_{k \uparrow \lambda_j} \mu_0^k = +\infty,
  \qquad
  \lim_{k \downarrow \lambda_j} \mu_0^k = -\infty,
\end{equation}
hence $\mu_0^k$ admits a unique simple zero in $(\lambda_j,\lambda_{j+1})$; see Figure \ref{figmu0k}. These spectral features have many interesting physical applications. For instance, in optics:
\begin{enumerate}
  \item (Band gaps)
        if $\mu_0^k<0$, the incident light decays exponentially and can not penetrate the media, making the media looks dark \cite{liu2000locally}.
  \item (Near-zero index)
        if $\mu_0^k\approx 0$, then the effective refractive index is near zero, so the light keeps nearly uniform phase in the media \cite{liberal2017near}.
  \item (Strong dispersion)
        as $k \uparrow \lambda_j$ from below, the slope $\frac{d\mu_0^k}{dk}$ diverges, and the medium becomes highly dispersive \cite{pedersen2008limits}.
\end{enumerate}
See Figure \ref{figband} for an illustrative explanation.
\end{remark}

\begin{remark}
By the expression \eqref{eq:mu0k_def}, the singularities of $\mu_0^k$ consist of $\Sigma_{D,1}$. It is natural to ask how large is $\Sigma_{D,1}$ relative to $\Sigma_D$. A recent work \cite{steinerberger2025dirichleteigenfunctionsnonzeromean} established the sharp lower bound, up to a logarithmic factor: for any smooth domain $D\subset \mathbb{R}^d$, let $\{\varphi_k\}_{j=1}^{\infty}$ denote a sequence of Laplacian eigenfunctions with Dirichlet boundary conditions.
There exists a constant $C>0$ that depends only $D$ such that, for all $n$ sufficiently large,
\begin{equation}
    \#\left\{1\leq k \leq n: \int_D \varphi_k(y)\,dy\neq 0\right\} \geq C\frac{n^{1/d}}{\sqrt{\log n}}.
\end{equation}
The extreme case is the ball in $\mathbb{R}^d$, where among the first $n$ eigenfunctions only $\sim n^{1/d}$ have a non-zero mean.
\end{remark}

\begin{remark}\label{kindisp}
    The assumption $k \notin \Sigma_D$ appears to be indispensable. Indeed, \cite{zhikov_extension_2000, fu_homogenization_2023} showed that, for the operator $L_{\varepsilon}$ on a bounded domain with Dirichlet boundary conditions, as $\varepsilon \to 0$ the limiting spectrum consists of two parts: the Dirichlet eigenvalues of $-\Delta$ on $D$, which form the essential spectrum of the `limiting operator', and a collection of discrete real eigenvalues characterized by $L_0+k^2\mu_0^k$, see Figure \ref{figmu0k} for an illustration. 
    Heuristically, upon passing from the bounded-domain problem to the scattering setting, the essential spectrum remains persistent, whereas the discrete spectrum shifts as resonances to the lower half-plane. Since the wave number $k$ is positive, it is natural to require $k\notin \Sigma_D$.
\end{remark}

One key ingredient in the proof of Theorem \ref{mainresult1} is the boundary layer estimates of the homogenized scattering problem \eqref{equ0}. Denote
\begin{equation}
    \mathcal{E}(u) := \varepsilon^{1/2} \| u \|_{H^1(\Omega)} + \| \nabla u \|_{L^2(\mathcal{O}_{(d+2)\varepsilon})} + \varepsilon \| \nabla^2 u \|_{L^2(\Omega \setminus \overline{\mathcal{O}_{(d-1)\varepsilon}})},
\end{equation}
where $\mathcal{O}_{r}$, $r>0$, denotes the boundary layer $\{x\in \Omega: \mathrm{dist}\,(x,\partial \Omega) < r\}$. Recall that $\Omega$ is assumed to be a Lipschitz domain.

\begin{theorem}\label{blethm}
    Fix $k\in (0,\infty)\setminus \Sigma_{D,1}$. Let $\widehat{u}_0$ be the solution of the homogenized scattering problem \eqref{equ0}. Then, for any $r>0$ such that $\overline{\Omega} \subset B_r$, we have
    \begin{equation}\label{Eu0main}
        \mathcal{E}(\widehat{u}_0) \leq C\varepsilon^{1/2} \|\uinc\|_{H^1(B_{r+1})},
    \end{equation}
where $C>0$ is a constant that depends only on $d,k,r,D$ and $\Omega$.
\end{theorem}

\begin{remark}
    In the proof of Theorem \ref{blethm}, we used the condition that $I-A_0$ is either positive or negative semi-definite. This condition ensures the existence of a layer potential representation of $\widehat{u}_0$ and the regularity of the associated density; see Theorem \ref{thm44}. In our setting, the coefficient $A_{\varepsilon}=I$ in the background medium $\Omega\setminus \overline{D_{\varepsilon}}$, which yields that $I-A_0$ is positive definite. 
    
    Our setting can be relaxed as follows. Suppose that $A_{\varepsilon} =A(x/\varepsilon)$ in $\Omega\setminus \overline{D_{\varepsilon}}$ for some symmetric, positive-definite periodic field $A(y)$. If we add the requirement that $I-A_0$ is either positive or negative semi-definite (that is to say, $A_0$ is not very anisotropic), then all results in this work still hold. Alternatively, if we assume that $\Omega$ is smooth, then \eqref{Eu0main} follows directly from elliptic regularity.     Therefore, all our conclusions continue to hold without the necessity of imposing that $I-A_0$ is either positive or negative semi-definite.
\end{remark}

Another key step in the proof of Theorem \ref{mainresult1} is the uniform boundedness of the outgoing resolvent $R_{\varepsilon}(z)$ of the operator $L_{\varepsilon}$ (see Theorem \ref{meromorphic_continuation} for definition and existence), which could be of independent interest. 

As discussed in Remark \ref{kindisp}, the map $z\mapsto \mu_0^z$ sends points of $\Sigma_{D,1}$ to infinity. Hence, the local spectral property of the homogenized operator $L_0+z^2\mu_0^z$ near $\Sigma_{D,1}$ corresponds to the spectral property of $L_0$ at infinity, which is beyond the scope of this paper.
Therefore, in what follows, we restrict our attention to resonances of $L_{\varepsilon}$ lying in 
\begin{equation}
    \mathbb{C}^{\circ} := \mathbb{C} \setminus \Sigma_{D,1}.
\end{equation}
We denote by $\Sigma_{\mathrm{hom}}$ the set of scattering resonances associated with $L_0+z^2\mu_0^z$, which is, by Definition \ref{DefA2z}, contained in $ \mathbb{C}^{\circ}$. 

The following theorems are stated for $d\geq 2$ to be odd. For even $d$, see Remark \ref{oddeven} for the necessary modifications.

\begin{theorem}\label{thm5121main}
    Fix $z\in \mathbb{C}\setminus (\Sigma_D\cup \Sigma_{\mathrm{hom}}) $. For any $r>0$ such that $\overline{\Omega} \subset B_r$, we have
    \begin{equation}
        \|R_{\varepsilon}(z)\|_{L^2(B_r)\rightarrow L^2(B_r)} \leq C,
    \end{equation}
    where $C>0$ is a constant that depends only on $d,z,r,D$ and $\Omega$.
\end{theorem}

By Theorem \ref{thm5121main}, the accumulation points of the scattering resonances of $L_{\varepsilon}$ are contained in the set $\Sigma_{\mathrm{hom}} \cup \Sigma_D$. In fact, we have the following partial converse result.

\begin{definition}
    The limiting set $\Sigma_{\lim} \subset \mathbb{C}^\circ$ for the resonances of $L_{\varepsilon}$ is defined as follows:
    \begin{enumerate}
        \item if there exists $\varepsilon_i\rightarrow 0$ and resonance $z_i $ of $R_{\varepsilon_i}(z)$ such that $z_i \rightarrow z_0 \in \mathbb{C}^\circ$, then $z_0\in \Sigma_{\lim}$;
        \item if $z_0\in \Sigma_{\lim}$, then, in any neighborhood of $z_0$ there exists at least one resonance of $R_{\varepsilon}(z)$ for sufficiently small $\varepsilon$.
    \end{enumerate}
\end{definition}

\begin{theorem}\label{maintheoremresonance}
    We have $\Sigma_{\lim} = \Sigma_{\mathrm{hom}} \cup \Sigma_{D,0}$. Define $\nu:\Sigma_{\lim}\rightarrow (0,1/2]$ by
\begin{equation}
    \nu(z_0) :=\left\{\begin{aligned}
        & \frac{1}{2}, && z_0 \in \Sigma_{D,0}, \\
        & \frac{1}{2(2m(z_0)+1)}, && z_0 \in \Sigma_{\mathrm{hom}},
    \end{aligned} \right. 
\end{equation}
where $m(z_0)>0$ denotes the order of $z_0$ of $(I+K_0(z))^{-1}$ (see \eqref{DefK0} for the definition of $K_0(z)$).
Then, for any $z_0 \in \Sigma_{\lim}$, there exists resonance $z_\varepsilon $ of $R_{\varepsilon}(z)$ such that 
\begin{equation}
    |z_\varepsilon-z_0|<C\varepsilon^{\nu(z_0)-}.
\end{equation}
More precisely, for any $h >0$, there exists a constant $C>0$ that depends only on $d,z_0,h,D$ and $\Omega$ such that $|z_\varepsilon-z_0|<C\varepsilon^{\nu(z_0)-h}$.

\end{theorem}

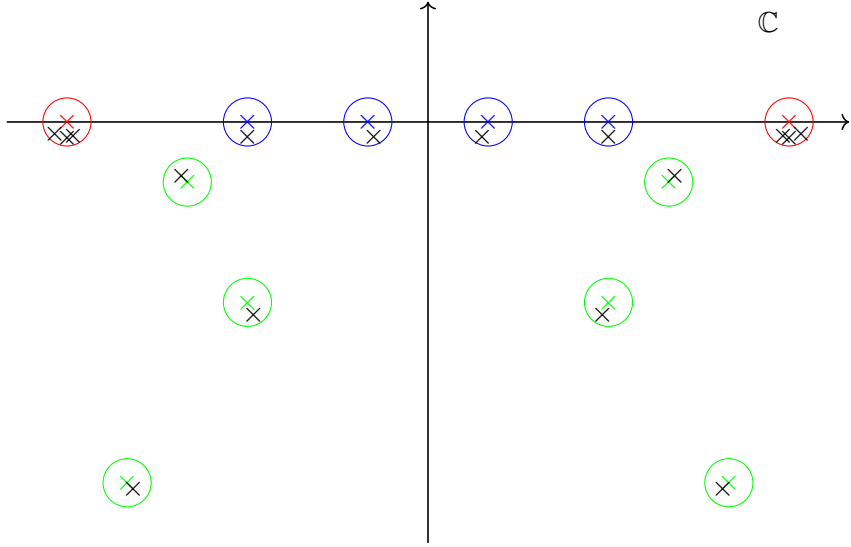
\begin{figure}[h]
\begin{center}
    \begin{tikzpicture}[scale=0.8]
    \node[below left] at ( 9,2 ){$\mathbb{C}$};
    \draw[->, line width = 0.6pt] (-4,0) -- ( 10, 0) ;
    \draw[->, line width = 0.6pt] (3, -7 ) -- ( 3, 2) ;
    \draw[blue] (0,0) circle (4mm); 
    \draw[blue] (0,0) node {$\times$};
    \draw (0,-.25) node {$\times$};
    \draw[blue] (6,0) circle (4mm); 
    \draw[blue] (6,0) node {$\times$};
     \draw (6,-.25) node {$\times$};
    \draw[blue] (2,0) circle (4mm); 
    \draw[blue] (2,0) node {$\times$};
    \draw (2.1,-.25) node {$\times$};
    \draw[blue] (4,0) circle (4mm); 
    \draw[blue] (4,0) node {$\times$};
    \draw (3.9,-.25) node {$\times$};
    \draw[red] (-3,0) circle (4mm); 
    \draw[red] (-3,0) node {$\times$};
    \draw (-3,-.27) node {$\times$};
    \draw (-3.2,-.2) node {$\times$};
    \draw (-2.9,-.23) node {$\times$};
    \draw[red] (9,0) circle (4mm); 
    \draw[red] (9,0) node {$\times$};
    \draw (9,-.27) node {$\times$};
    \draw (9.2,-.2) node {$\times$};
    \draw (8.9,-.23) node {$\times$};
    \draw[green] (0,-3) circle (4mm); 
    \draw[green] (0,-3) node {$\times$};
    \draw (0.1,-3.2) node {$\times$};
    \draw[green] (6,-3) circle (4mm); 
    \draw[green] (6,-3) node {$\times$};
    \draw (5.9,-3.2) node {$\times$};
    \draw[green] (-1,-1) circle (4mm); 
    \draw[green] (-1,-1) node {$\times$};
    \draw (-1.1,-0.9) node {$\times$};
    \draw[green] (7,-1) circle (4mm); 
    \draw[green] (7,-1) node {$\times$};
    \draw (7.1,-0.9) node {$\times$};
    \draw[green] (-2,-6) circle (4mm); 
    \draw[green] (-2,-6) node {$\times$};
    \draw (-1.9,-6.1) node {$\times$};
    \draw[green] (8,-6) circle (4mm); 
    \draw[green] (8,-6) node {$\times$};
    \draw (7.9,-6.1) node {$\times$};
    \end{tikzpicture}
\end{center}
\caption{An illustration of the resonances distribution: the blue crosses denote the eigenvalues in $\Sigma_{D,0}$, the red crosses denote the eigenvalues in $\Sigma_{D,1}$, the green crosses denote the resonances in $\Sigma_{\mathrm{hom}}$, and the black crosses denote the resonances of $L_{\varepsilon}$. }
\label{resodist}
\end{figure}

\begin{remark}
    We establish in Theorem \ref{maintheoremresonance} a convergence rate for the scattering
resonances of $L_{\varepsilon}$ that remain a positive distance from $\Sigma_{D,1}$.
As noted above, points of $\Sigma_{D,1}$ are pathological and correspond to the spectral
behavior of $L_0$ at infinity. Motivated by the bounded-domain case, we formulate the following conjecture:
    \begin{quote}
        \textbf{Conjecture.} Fix $z_0 \in \Sigma_{D,1}$. Denote by $\mathrm{Res}\,L_{\varepsilon}$ the set of scattering resonances of $L_{\varepsilon}$. Then for any $r>0$, 
        \begin{equation}
            \lim_{\varepsilon \rightarrow 0} \# \{z:|z-z_0|<r,\,z\in \mathrm{Res}\,L_{\varepsilon}\} = \infty.
        \end{equation}
    \end{quote}
    See Figure \ref{resodist} for an illustration.
\end{remark}

By combining Theorem \ref{mainresult1} and the dual method, and assuming stronger regularity of $\partial \Omega$, we obtain the optimal $L^2$ convergence rate.

\begin{theorem}\label{mainresult2}
    Under the same conditions and notations as those in Theorem \ref{mainresult1}, we further assume that $\Omega$ is a $C^{1,1}$ domain. Then, we have
    \begin{equation}\label{251optimal}
        \| u_\varepsilon - \Lambda_\varepsilon \widehat{u}_0  \|_{L^2(B_r)}  \leq C\varepsilon  \|\uinc\|_{H^1(B_{r+1})}  ,
    \end{equation}
    where $C>0$ is a constant that depends only on $d,k,r,D$ and $\Omega$.
\end{theorem}

As an application of Theorem \ref{mainresult1} and Theorem \ref{mainresult2}, we can obtain the point-wise convergence rates of the far-field pattern for the heterogeneous scattering problem \eqref{maineq}, which has potential application in inverse scattering problems. More precisely, let $u^s_{\varepsilon} := u_{\varepsilon} - \uinc$ be the scattered wave. Then, it is well known that $u^s_{\varepsilon}$ admits the following asymptotics at infinity:
\begin{equation}\label{633}
    u^s_\varepsilon(x) =  \frac{e^{\mathrm{i}k|x| } }{ |x|^{(d-1)/2} } \Big\{ u^{\infty}_\varepsilon( \theta )  + \mathcal{O}(|x|^{-1}) \Big\} ,\quad x\rightarrow \infty,
\end{equation}
where $\theta= \frac{ x}{ |x|}$; see, e.g., \cite[Chapter 9]{taylor2010partial}. When $\uinc=\mathrm{e}^{\mathrm{i}k\omega \cdot x}$ for some $\omega\in S^{d-1}$, the far-field pattern $u^\infty_\eps$ is the so-called scattering amplitude. The precise formula for $u^\infty_\eps$ is given in \eqref{6341}. We have the following quantitative result for the convergence of $u^\infty_\eps$ to the far-field pattern $u^\infty_0$ of the effective model \eqref{equ0}. 

\begin{theorem}\label{thmamplitd}
    Fix $k\in (0,\infty) \setminus \Sigma_D$. Let $u^{\infty}_\varepsilon :S^{d-1} \rightarrow \mathbb{C}$ and $u^{\infty}_0 :S^{d-1} \rightarrow \mathbb{C}$ be the far-field patterns of the scattering problem \eqref{maineq} and the homogenized scattering problem \eqref{equ0}, respectively. Then, for any integer $n\geq 0$, we have
        \begin{equation}
            \| u^{\infty}_\varepsilon - u^{\infty}_0 \|_{C^n(S^{d-1})} \leq C \varepsilon^{1/2}\| \uinc \|_{H^1(B_{r+1})};
        \end{equation}
        if further assume that $\Omega$ is a $C^{1,1}$ domain, then we have
         \begin{equation}
            \| u^{\infty}_\varepsilon - u^{\infty}_0 \|_{C^n(S^{d-1})} \leq C \varepsilon\| \uinc \|_{H^1(B_{r+1})}.
        \end{equation}
        Here $r>0$ is chosen so that $\overline{\Omega} \subset B_r$, and the constant $C>0$ is a constant that depends only on $d,k,n,r,D$ and $\Omega$.
\end{theorem}

\begin{remark}\label{oddeven}
    Strictly speaking, all of the above results are stated in odd dimensions. In even
dimensions, the fundamental solution $G^{z}$ of $\Delta+z^{2}$ has a logarithmic
singularity at $z=0$ (see \eqref{320Gk}), and the free outgoing resolvent admits a
meromorphic continuation only to the logarithmic cover (often represented by a branch cut
along the negative imaginary axis). For our purposes, this causes no essential difficulty,
since we study the limiting behavior of the scattering problem and the local behavior of
the outgoing resolvent $R_{\varepsilon}(z)$ near its poles. The following modifications
suffice when $d$ is even:
\begin{enumerate}
  \item Theorems \ref{mainresult1}, \ref{blethm}, \ref{mainresult2}, and \ref{thmamplitd} remain valid,
        since $0\notin \Sigma_{D}$.
  \item Theorems \ref{thm5121main} and \ref{maintheoremresonance} hold upon replacing
        $\mathbb{C}$ by the logarithmic cover.
\end{enumerate}
\end{remark}

\subsection{Relation to previous works}\label{previouswork}

This work lies at the intersection of high-contrast homogenization, wave scattering, and the study of subwavelength resonators. 

It was observed early on that media with contrasts of order $\varepsilon^2$ can exhibit macroscopic properties not observed in natural materials over certain frequency ranges \cite{auriault1985dynamique}. The double-porosity model provides one of the first mathematical frameworks for such high-contrast regimes. Arbogast, Douglas, and Hornung \cite{arbogast1990derivation} studied single-phase flow in fractured porous media, where the permeability of the rock matrix is significantly smaller than that of the fracture network. They proved weak $L^2$ convergence of $u_\varepsilon$ to a coupled effective PDE system with fast and slow variables. Allaire \cite{allaire1992homogenization} refined this through two-scale convergence, and Zhikov \cite{zhikov_extension_2000} established two-scale resolvent convergence along with spectral convergence in the Hausdorff (a.k.a. Kuratowski) sense. The Hausdorff spectral convergence was later improved to quantitative convergence by Fu \cite{fu_homogenization_2023}. For Maxwell's equation, Bouchitté and Felbacq \cite{bouchitte2004homogenization} showed that, under $\varepsilon^2$-contrast, frequency-dependent effective parameters can become negative within certain frequency bands.

As one of the major contributions of this work, Theorem \ref{mainresult1} and \ref{mainresult2} extend the previous results in \cite{arbogast1990derivation, allaire_homogenization_1992, zhikov_extension_2000, bouchitte2004homogenization} by improving the two-scale convergence to the convergence rates in $L^2$ and $H^1$ norms. Our proofs, developed for the scattering setting, also apply (with simplifications) to problems in bounded domains. We believe that the techniques developed in this paper can be used to study more complex models in fluid mechanics, electromagnetics, elasticity, and so on.

It is worth mentioning some related works on high-contrast homogenization. \cite{russell_quantitative_2018, shen_large-scale_2021, fu2025uniform} studied high-contrast homogenization for periodic elliptic systems on bounded domains, as well as large-scale regularity. For fully periodic structures, K. Cherednichenko \textit{et al.} \cite{cherednichenko2016resolvent, cherednichenko2018operator, cherednichenko2019time, cherednichenko2020effective, cherednichenko2022sharp, cherednichenko2020order, cherednichenko2022operator} obtained the optimal $L^2$ convergence rates in various settings, including elasticity, electromagnetics, and time-dispersive media. For high-contrast homogenization in random settings, we refer to \cite{cherdantsev2021high, clozeau2024quantitative, bella2025quantitative, bonhomme2025homogenization, armstrong2024renormalization, pinaud2024scaling} for recent progress. Extending this work to the random setting would be of great interest.

Understanding the scattering behavior of highly oscillatory heterogeneous structures is a fundamental problem in many areas. Cakoni \emph{et al.} \cite{cakoni2016homogenization} studied homogenization of scattering by highly oscillatory anisotropic media (without high contrast), and Garnier \emph{et al.} \cite{garnier2023scattered} extended this to random settings. See \cite{MR4026931} also for a special case where random inhomogeneities occur only in the refraction index and geometrically occupy a slab. For scattering resonances, there are some works on Schrödinger equations with highly oscillatory potentials \cite{duchene2011scattering, duchene2014scattering, drouot2018scattering, drouot2018resonances}. To the best of our knowledge, this paper is the first to study homogenization of the scattered wave and the scattering resonances for operators of divergence form with high-contrast coefficients, together with quantitative error estimates.

We should note that there are several different types of high-contrast subwavelength resonators across various high-contrast regimes, and this work only treats one of them. Another important example is the air bubble in water. The strong density and compressibility contrasts between air and water induce the so-called Minnaert resonance \cite{minnaert1933xvi}. Based on layer potentials, the spectral theory of Neumann-Poincar\'{e} operator and capacity matrices that capture bubble interactions, Ammari \emph{et al.} developed a comprehensive theory with applications to subwavelength physics. We list some, but not all, for the interested readers: \cite{ammari2017subwavelength, ammari2017sub, ammari2018minnaert, ammari2019bloch, ammari2020honeycomb, ammari2020high, feppon2023homogenization, ammari2024wave, ammari2024effective, ammari2025analysis}. The key difference from our setting is that the Minnaert resonance corresponds to surface mode and generates stronger interactions between bubbles, whereas the resonators considered in this work correspond to body mode with weaker pairwise coupling. As a consequence, effective medium theories for bubbly media typically require restrictive assumptions such as dilute configurations or precise tuning of the incident frequency \cite{caflisch1985effective, caflisch1985wave, ammari2017effective, alsenafi2023foldy, mukherjee2024dispersive}. It would be of great interest to derive an effective medium theory for bubbly media in densely distributed cases.

\section{Convergence rates for the scattering problem}\label{secconrateprob}

In this section, we establish the $L^2$ and weighted $H^1$ convergence rates for the homogenization of scattering problems \eqref{maineq}, assuming that some suitable estimates hold. These estimates will be proved in the subsequent sections.

\subsection{Preliminaries}

We collect some basic notations and lemmas that are frequently used in this paper. Throughout the paper  it is assumed that $\varepsilon \in (0,1)$.

The boundary cut-off function $\eta_\eps \in C_0^{\infty} (\Omega; [0,1])$ is chosen such that
\begin{equation}\label{cut-off_def}
	\left\{
	\begin{aligned}
		& \mathrm{supp}\,\eta_{\varepsilon}\subset \Omega \setminus \mathcal{O}_{d\varepsilon}, \\
		&\eta_{\varepsilon} \equiv 1 \quad \mathrm{in}\  \Omega \setminus \mathcal{O}_{(d+1)\varepsilon}, \\
		& \|\nabla \eta_{\varepsilon}\|_{L^{\infty}(\Omega)} \leq C\varepsilon^{-1} ,
	\end{aligned}\right.
\end{equation}
where $\mathcal{O}_{r}$, $r>0$, denotes the boundary layer $\{x\in \Omega: \mathrm{dist}\,(x,\partial \Omega) < r\}$.

Fix a function $\xi \in C_0^{\infty}(B_{1/2})$ such that $0\leq \xi (x)\leq 1$ and $\int_{\mathbb{R}^d}\xi (x)\,dx=1$. We define the smoothing operator $S_{\varepsilon}$ by
\begin{equation}\label{epsi_smoother}
	S_{\varepsilon}(v)(x):=\varepsilon^{-d} \int_{\mathbb{R}^d} v(x-z)\xi(z/\varepsilon)\,dz, \qquad \forall v\in L^1_{\rm loc}(\mathbb{R}^d).
\end{equation}

Recall that $f^{\varepsilon}(x):= f(x/\varepsilon)$ for a function $f$.

\begin{lemma}\label{Svarepsilonproperty1}
    Let $U\subset \mathbb{R}^d$ be an open set. For any $ f \in L^2(\mathbb{T}^d)$ and $v \in L^2(\mathbb{R}^d)$, we have
    \begin{equation}
        \| f^{\varepsilon} S_\varepsilon(v)\|_{L^2( U)}\leq C \|f \|_{ L^2(Y ) } \|v \|_{L^2(U(\varepsilon))},
    \end{equation}
    where $U(\varepsilon):=\{x\in \mathbb{R}^d: \mathrm{dist}\,(x,U)<\varepsilon/2\}$ and $C>0$ is a universal constant.
\end{lemma}
\begin{proof}
See \cite[Proposition 3.1.5]{shen2018periodic}. 
\end{proof}

\begin{lemma}
\label{propertySepsilondifference}
    Let $U\subset \mathbb{R}^d$ be an open set. For any $v\in H^1(\mathbb{R}^d)$, we have
    \begin{equation}
        \left\|S_{\varepsilon}( v )-v  \right\|_{L^2( U)}\leq C \varepsilon  \| \nabla v \|_{L^2(U(\varepsilon)) } ,
    \end{equation}
    where $U(\varepsilon):=\{x\in \mathbb{R}^d: \mathrm{dist}\,(x,U)<\varepsilon/2\}$ and the constant $C>0$ depends only on $d$.
\end{lemma}
\begin{proof}
By the density argument, we may assume that $ v \in C_0^1(\mathbb{R}^d)$. By the fundamental theorem of calculus, we have
\begin{equation}
    \begin{aligned}
        S_{\varepsilon}( v )(x)- v(x)  & = \int_{ |y|\le 1/2 } \big( v(x-\varepsilon y) -v(x) \big) \xi(y)  \,dy\\
    & = \int_{ |y|\le 1/2 } \left(\int_0^1  \varepsilon y \cdot \nabla v(x-t\varepsilon y) \,dt \right) \xi(y)\,dy \\
    & = \int_{ |y|\le 1/2 }  \int_0^1   \varepsilon y \cdot\nabla v(x-t\varepsilon y) \xi(y)\,dt  dy.
    \end{aligned}
\end{equation}
It follows that 
\begin{equation}
    \begin{aligned}
        \left\|S_{\varepsilon}( v )-v  \right\|_{L^2( U)}   \leq  \int_{ |y|\le 1/2 }  \int_0^1   \| \varepsilon y \cdot\nabla v(\cdot-t\varepsilon y) \xi(y) \|_{L^2(U)}\,dt  dy  \leq C \varepsilon  \| \nabla v \|_{L^2(U(\varepsilon)) } .
    \end{aligned}
\end{equation}
The proof is complete.
\end{proof}

Recall that $\Omega\subset \mathbb{R}^d$ is a bounded Lipschitz domain, and $\mathcal{O}_r=\{x\in \Omega: \mathrm{dist}\,(x,\partial \Omega) < r\}$.

\begin{lemma}\label{blelem}
For any  $u\in H^1(\Omega)$, we have the boundary layer estimate:
\begin{equation}
    \| u \|_{L^2(\mathcal{O}_\varepsilon)} \leq C \varepsilon^{1/2} \|u \|_{H^1(\Omega)},
\end{equation}
where the constant $C>0$ depends only on $\Omega$.
\end{lemma}
\begin{proof}
    See \cite[Proposition 3.1.7]{shen2018periodic}.
\end{proof}

The flux corrector plays an important role in homogenization theory.

\begin{lemma}
	\label{thmflux}
	Suppose that $b=(b_{ij} )_{i,j=1}^d $ is a tensor field on $Y$ satisfying
	\begin{equation}\label{requireflux}
		\frac{\partial }{\partial y_i} b_{ij} (y) =0, \qquad \mathrm{for\ all\ } y\in Y,\ 1\leq j \leq d,  
	\end{equation}
	and
	\begin{equation}
		\int_Y b_{ij}(y)\,dy =0, \qquad \mathrm{for\ all\ } 1\leq i,j \leq d.
	\end{equation}
	Then there exists a 
	tensor field $\Psi = (\Psi_{kij} )_{k,i,j=1}^d$ on $\mathbb{T}^d$ such that
	\begin{equation}
		\label{eq:FtoPsi}
		\frac{\partial }{\partial y_k} \Psi_{kij} (y)= b_{ij} (y), \qquad \mathrm{for\ all\ } y\in Y, \ 1\leq i,j\leq d,
	\end{equation}
	and
	\begin{equation}	\label{eq:FtoPsi2}
		\Psi_{kij} (y)= - \Psi_{ikj} (y), \qquad \mathrm{for\ all\ } y\in Y,\  1\leq k,i,j\leq d.
	\end{equation}
Moreover, there exists a constant $C>0$, which depends only on $d$, such that 
	\begin{equation}\label{Psib}
		\| \Psi \|_{H^1(\mathbb{T}^d)} \leq C  \|b\|_{L^2(Y)} .
	\end{equation}
	We refer to
	$\Psi$ as the flux corrector associated with $b$. 
\end{lemma}
\begin{proof}
See \cite[Proposition 3.1.1]{shen2018periodic}.
\end{proof}

The following extension lemma is a fundamental tool for studying problems on high-contrast or perforated domains. Recall the definition \eqref{def0D} for $D_{\varepsilon}$.

\begin{lemma}
\label{extensionoperator}
    Let $U\subset \mathbb{R}^d$ be an open set such that $\Omega \subset U$. There exists an extension operator 
    \begin{equation}
        \mathcal{P}_\varepsilon : H^1(U\setminus \overline{D_\varepsilon} ) \rightarrow H^1(U) 
    \end{equation}
such that, for any $u \in H^1(U\setminus \overline{D_\varepsilon} )$,
    \begin{equation}\label{contorl_Pepsilon}
         \|\mathcal{P }_\varepsilon u \|_{H^1(U)} \leq C \| u \|_{H^1( U\setminus \overline{D_\varepsilon} ) }  \quad \mathrm{and} \quad \|\nabla \mathcal{P }_\varepsilon u \|_{L^2(U)} \leq C \| \nabla u \|_{L^2(U\setminus \overline{D_\varepsilon} ) } ,
    \end{equation}
where the constant $C>0$ depends only on $D$.
\end{lemma}
\begin{proof}
    See Step 1 in the proof of \cite[Lemma 3.3]{fu2025uniform}.
\end{proof}

We also need some lemmas to handle the Sommerfeld radiation condition.

\begin{lemma}
    For any $k>0$ and $r>0$. If $u \in H^2_{\mathrm{loc}}(\mathbb{R}^d \setminus \overline{B_r})$ satisfies 
    \begin{equation}
        \left\{
        \begin{aligned}
            &(\Delta +k^2)u =0 && \mathrm{in}\ \mathbb{R}^d\setminus \overline{B_r}, \\
            &u\in \mathrm{SRC}(k).
        \end{aligned}\right.
    \end{equation}
Then,
\begin{equation}\label{Renega}
    \mathrm{Re}\int_{\partial B_r} \frac{\partial u}{\partial \mathbf{n}} \overline{ u }\,dS(x) \leq 0.
\end{equation}
Moreover, if $v$ satisfies the same conditions as $u$, then
\begin{equation}\label{smiden}
    \int_{\partial B_r } \left( u \frac{ \partial v }{ \partial \mathbf{n}} - v \frac{ \partial u }{ \partial \mathbf{n}} \right) dS(x)  =0.
\end{equation}
\end{lemma}
\begin{proof}
    The inequality \eqref{Renega} is well-known; see \cite[(1.52)]{cakoni2022inverse} for the case $d=3$, where the arguments readily adapt to any $d\geq 2$. Therefore, we only prove \eqref{smiden}.

By the far-field expansion of $u$ (\cite[Page 294]{mclean2000strongly}), we have
    \begin{equation}\label{farfield}
        |u(x)| = O( |x|^{-(d-1)/2}) \qquad \mathrm{as} \ x\rightarrow \infty.
    \end{equation}
This can be seen by substituting the standard asymptotic expansion for Hankel function into the Green's representation
   \begin{equation}\label{Greenrep}
        u(x)= \int_{\partial B_r} \frac{\partial u}{ \partial \mathbf{n}_y }(y)  G^k(x,y) -\frac{\partial G^k}{ \partial \mathbf{n}_y }(x,y) u(y)  \,dS(y), \qquad x\in \mathbb{R}^d\setminus \overline{B_r},
    \end{equation}
where $G^k $ is the fundamental solution to the Helmholtz operator $\Delta +k^2$:
\begin{equation}\label{320Gk}
    G^k(x,y) := -\frac{\mathrm{i}}{4} \left( \frac{k}{2\pi |x-y|} \right)^{\frac{d-2}{2}} H^1_{\frac{d-2}{2}} (k|x-y|), \qquad d\geq 2,
\end{equation}
(\cite{agmon1990representation, nachman1988reconstructions}) and $H^1_{\mu}(z)$ denotes the Hankel function of the first kind of order $\mu$.

By Green's identity and the equation $(\Delta+k^2)u = (\Delta+k^2)v=0$ in $\mathbb{R}^d\setminus \overline{B_r}$, we have
    \begin{equation}
    \begin{aligned}
        \int_{\partial B_r } \left( u \frac{ \partial v }{ \partial \mathbf{n}} - v \frac{ \partial u }{ \partial \mathbf{n}} \right) dS(x) & = \lim_{R\rightarrow \infty}\int_{\partial B_R } \left( u \frac{ \partial v }{ \partial \mathbf{n}} - v \frac{ \partial u }{ \partial \mathbf{n}} \right)  dS(x) \\
        &=\lim_{R\rightarrow \infty}\int_{\partial B_R } u(\partial_{|x|} -\mathrm{i}k) v -  v(\partial_{|x|} -\mathrm{i}k) u \, dS(x) \\
        & = \lim_{R\rightarrow \infty} \int_{\partial B_R } O(R^{-d})\,dS(x) \\
        & = 0,
    \end{aligned}
\end{equation}
where the third equality follows from \eqref{farfield} and the Sommerfeld radiation conditions for both $u $ and $v$.
\end{proof}

\begin{lemma}\label{rellichunique}
    Let $U\subset \mathbb{R}^d$ be a bounded Lipschitz domain. For any $k>0$, if $u \in H^2_{\mathrm{loc}}(\mathbb{R}^d \setminus \overline{U})$ satisfies 
    \begin{equation}
        \left\{
        \begin{aligned}
            &(\Delta +k^2)u =0 && \mathrm{in}\ \mathbb{R}^d\setminus \overline{U}, \\
            &u\in \mathrm{SRC}(k),
        \end{aligned}\right.
    \end{equation}
and
\begin{equation}\label{Imgre}
    \mathrm{Im}\int_{\partial U} \frac{\partial u}{\partial \mathbf{n}} \overline{ u }\,dS(x) \leq 0,
\end{equation}
then $u\equiv 0$ in $\mathbb{R}^d \setminus \overline{U}$.
\end{lemma}
\begin{proof}
    This is the Rellich uniqueness theorem; see \cite[Lemma 9.9]{mclean2000strongly}.
\end{proof}

\subsection{\texorpdfstring{$L^2$ convergence rate}{L2 convergence rate}}

In this section, we establish the $L^2$ convergence rate. The following lemma is the starting point of the analysis, which extends \cite[Theorem 4.2]{fu2025uniform} to the scattering setting. Consider the following function defined in \cite{zhikov_extension_2000}:
\begin{equation}\label{defJikovfunc}
    \beta(z): = \int_Y (\Delta_D +z^2)^{-1}[1] (y)\,dy.
\end{equation}
We call it Zhikov's function in this paper. It is clear that $\mu_0^k = 1-k^2 \beta(k)$.

\begin{lemma}\label{basiclemmasca}
    Let $k\in (0,\infty) \setminus \Sigma_D$, and $U\subset \mathbb{R}^d$ be a Lipschitz domain such that $\overline{\Omega}\subset U$. Let $w_\varepsilon,\, \widehat{w}_0  \in H^1(U) $ satisfy
\begin{equation}\label{generaleq}
    (L_{\varepsilon}+k^2) w_\varepsilon =f,\quad \mathrm{and}\quad  (L_0 +k^2 \mu^k_0 ) \widehat{w}_0 = \mu_0^k f  \qquad \mathrm{in}\ U 
\end{equation}
for some $f\in L^2(U)$. Let 
    \begin{equation}
        z_{\varepsilon}: = w_\varepsilon - \Lambda_{\varepsilon}\widehat{w}_0 - \varepsilon \chi^{\varepsilon}   \cdot \eta_\varepsilon S_\varepsilon (  \nabla \widehat{w}_0  ) .
    \end{equation}
Then, for any $\psi ,\,\varphi \in H^1(U)$, we have
\begin{equation}
    \begin{aligned}
        &\big\langle ( L_{\varepsilon} +k^2) z_{\varepsilon}, \psi \big\rangle_{H^{-1}(U),H^1(U)} =\\
        &\varepsilon\int_{\Omega } \Psi^{\varepsilon} \nabla \big(\eta_{\varepsilon}  S_{\varepsilon} (\nabla\widehat{w}_0) \big)\cdot \overline{\nabla \varphi}\,dx + k^2 \int_{\Omega}  (\mu^k_0 -\Lambda_{\varepsilon} )\widehat{w}_0  \overline{\varphi} \,dx + \varepsilon \int_{D_{\varepsilon}}  (\nabla\Lambda)^{\varepsilon} \widehat{w}_0 \cdot \overline{\nabla \varphi}  \,dx \\
        &+ \varepsilon\int_{\Omega} A_{\varepsilon} \chi^{\varepsilon} \cdot \nabla\big( \eta_{\varepsilon}  S_{\varepsilon} (\nabla\widehat{w}_0) \big) \overline{\nabla \varphi} \,dx + \varepsilon^2 \int_{D_{\varepsilon}} (   \Lambda I +  \nabla \chi )^{\varepsilon} \eta_{\varepsilon}  S_{\varepsilon} (\nabla\widehat{w}_0) \cdot \overline{\nabla \varphi}\,dx \\
        & + \int_{\Omega} (\Lambda_{\varepsilon} A_{\varepsilon} - A_0 ) \big( \nabla\widehat{w}_0 - \eta_{\varepsilon}  S_{\varepsilon} (\nabla\widehat{w}_0)  \big)\cdot \overline{\nabla \varphi} \,dx  + \int_U (f - k^2 \Lambda_{\varepsilon} \widehat{w}_0)( \overline{\psi} - \overline{\varphi} )\,dx  \\
        &+ \varepsilon \int_{D_{\varepsilon}}  (\nabla  \Lambda)^{\varepsilon}\widehat{w}_0   \cdot \nabla (\overline{\psi} -\overline{\varphi})  \,dx+ \int_U A_{\varepsilon}  \big( \Lambda_{\varepsilon} \nabla\widehat{w}_0 + \nabla(\varepsilon \chi^{\varepsilon} \eta_{\varepsilon}  \cdot S_{\varepsilon} (\nabla\widehat{w}_0) )\big)  \cdot \nabla (\overline{\psi} -\overline{\varphi}) \,dx \\
        & + k^2 \beta(k)\int_{\Omega}  f \overline{\varphi} \,dx - k^2 \int_{\Omega} \varepsilon \chi^\varepsilon\cdot \eta_\varepsilon S_\varepsilon (  \nabla \widehat{w}_0  ) \overline{\psi} \,dx - \int_{\partial U} \frac{\partial \widehat{w}_0}{\partial \mathbf{n}} (\overline{\psi} - \overline{\varphi})\,dS(x) ,
    \end{aligned}
\end{equation}
where $\Psi$ is the flux corrector associated with $b(y) = A_0 -\big( I + \nabla \chi(y)  \big)\mathbbm{1}_{Y\setminus D}(y)$.
\end{lemma}

\begin{proof}
    Using the equations in \eqref{generaleq}, via direct computation we get
    \begin{equation}\label{bas1}
        \begin{aligned}
             &\big\langle (\nabla \cdot A_{\varepsilon} \nabla +k^2) (w_{\varepsilon}- \Lambda_{\varepsilon} \widehat{w}_0), \psi \big\rangle_{H^{-1}(U),H^1(U)} \\
             & = \big\langle \nabla \cdot A_0 \nabla  \widehat{w}_0  , \varphi \big\rangle_{H^{-1}(U),H^1(U)} - \big\langle \nabla \cdot A_{\varepsilon} \nabla ( \Lambda_{\varepsilon} \widehat{w}_0), \varphi \big\rangle_{H^{-1}(U),H^1(U)} \\
             &\quad\ + k^2 \int_{\Omega}  (\mu_0^k -\Lambda_{\varepsilon} )\widehat{w}_0  \overline{\varphi} \,dx + \int_U (f-k^2 \Lambda_{\varepsilon} \widehat{w}_0)(\overline{\psi}-\overline{\varphi})\,dx + k^2 \beta(k)\int_{\Omega}  f \overline{\varphi} \,dx \\
             & \quad\ +\int_U A_{\varepsilon} \nabla(\Lambda_{\varepsilon}\widehat{w}_0 )\cdot \nabla (\overline{\psi}- \overline{\varphi}) \,dx - \int_{\partial U} \frac{\partial \widehat{w}_0}{\partial \mathbf{n}} (\overline{\psi} - \overline{\varphi})\,dS(x),
        \end{aligned}
    \end{equation}
and
\begin{equation}\label{bas2}
    \begin{aligned}
        &-\big\langle \nabla \cdot A_{\varepsilon} \nabla  (\varepsilon \chi^{\varepsilon} \cdot \eta_\varepsilon S_\varepsilon (  \nabla \widehat{w}_0  )), \psi \big\rangle_{H^{-1}(U),H^1(U)} \\
        &=\int_{\Omega} A_{\varepsilon}(\nabla \chi)^{\varepsilon} \cdot \eta_{\varepsilon}  S_{\varepsilon} (\nabla\widehat{w}_0)  \cdot \overline{\nabla \varphi} \,dx + \varepsilon\int_{\Omega} A_{\varepsilon} \chi^{\varepsilon} \cdot \nabla\big( \eta_{\varepsilon}  S_{\varepsilon} (\nabla\widehat{w}_0) \big) \overline{\nabla \varphi} \,dx \\
        &\quad\ + \int_U A_{\varepsilon}  \nabla\big( \varepsilon \chi^{\varepsilon} \cdot \eta_{\varepsilon}  S_{\varepsilon} (\nabla\widehat{w}_0) \big) \cdot \nabla (\overline{\psi} -\overline{\varphi})   \,dx .
    \end{aligned}
\end{equation}
Summing \eqref{bas1}-\eqref{bas2} yields
\begin{equation}
    \big\langle (\nabla \cdot A_{\varepsilon} \nabla +k^2) z_{\varepsilon}, \psi \big\rangle_{H^{-1}(U),H^1(U)} = \mathrm{I} + \mathrm{II},
\end{equation}
where
\begin{equation}
    \begin{aligned}
        \mathrm{I}& = \big\langle \nabla \cdot A_0 \nabla  \widehat{w}_0  , \varphi \big\rangle_{H^{-1}(U),H^1(U)} - \big\langle \nabla \cdot A_{\varepsilon} \nabla ( \Lambda_{\varepsilon} \widehat{w}_0), \varphi \big\rangle_{H^{-1}(U),H^1(U)} \\
        &\quad\ + \int_{\Omega} A_{\varepsilon}(\nabla \chi)^{\varepsilon} \cdot \eta_{\varepsilon}  S_{\varepsilon} (\nabla\widehat{w}_0)  \cdot \overline{\nabla \varphi} \,dx,
    \end{aligned}
\end{equation}
and
\begin{equation}\label{termii}
    \begin{aligned}
        &\mathrm{II} =\\
        &k^2 \int_U  (\mu_0^k -\Lambda_{\varepsilon} )\widehat{w}_0  \overline{\varphi} \,dx - k^2 \int_{\Omega} \varepsilon \chi^\varepsilon\cdot \eta_\varepsilon S_\varepsilon (  \nabla \widehat{w}_0  ) \overline{\psi} \,dx + \varepsilon\int_{\Omega} A_{\varepsilon} \chi^{\varepsilon} \cdot \nabla\big( \eta_{\varepsilon}  S_{\varepsilon} (\nabla\widehat{w}_0) \big) \overline{\nabla \varphi} \,dx  \\
        & + \int_U (f-k^2 \Lambda_{\varepsilon} \widehat{w}_0)(\overline{\psi} - \overline{\varphi})\,dx+ k^2 \beta(k)\int_{\Omega}  f \overline{\varphi} \,dx + \varepsilon \int_{D_{\varepsilon}}  (\nabla  \Lambda)^{\varepsilon}\widehat{w}_0   \cdot \nabla (\overline{\psi} -\overline{\varphi})  \,dx \\
        & + \int_U A_{\varepsilon}  \big( \Lambda_{\varepsilon} \nabla\widehat{w}_0 + \nabla(\varepsilon \chi^{\varepsilon} \cdot \eta_{\varepsilon}  S_{\varepsilon} (\nabla\widehat{w}_0) )\big)  \cdot \nabla (\overline{\psi} -\overline{\varphi}) \,dx  - \int_{\partial U} \frac{\partial \widehat{w}_0}{\partial \mathbf{n}} (\overline{\psi} - \overline{\varphi})\,dS(x) .
    \end{aligned}
\end{equation}

It remains to deal with term I. By the definition, we have
 \begin{equation}\label{I1}
    \begin{aligned}
        \mathrm{I}& = \int_{\Omega } \big[ ( I + (\nabla \chi)^{\varepsilon}  )\mathbbm{1}_{\Omega \setminus D_{\varepsilon}}- A_0 \big] \eta_{\varepsilon}  S_{\varepsilon} (\nabla\widehat{w}_0) \cdot \overline{\nabla \varphi} \,dx \\
        &\quad \ + \varepsilon \int_{D_{\varepsilon}}  (\nabla\Lambda)^{\varepsilon} \widehat{w}_0 \cdot \overline{\nabla \varphi} \,dx + \varepsilon^2 \int_{D_{\varepsilon}} (   \Lambda I +  \nabla \chi )^{\varepsilon} \eta_{\varepsilon}  S_{\varepsilon} (\nabla\widehat{w}_0) \cdot \overline{\nabla \varphi} \,dx\\
        &\quad \ + \int_{\Omega} (\Lambda_{\varepsilon} A_{\varepsilon} - A_0 ) \big( \nabla\widehat{w}_0 - \eta_{\varepsilon}  S_{\varepsilon} (\nabla\widehat{w}_0)  \big)\cdot \overline{\nabla \varphi} \,dx    .
    \end{aligned}
\end{equation}
Thanks to the boundary cut-off $\eta_{\varepsilon}$, one has $\big[ ( I + (\nabla \chi)^{\varepsilon}  )\mathbbm{1}_{\Omega \setminus D_{\varepsilon}}- A_0 \big] \eta_{\varepsilon} = -b^{\varepsilon} \eta_{\varepsilon}$. Hence, by Lemma \ref{thmflux}, the first term on the right-hand side of \eqref{I1} becomes
\begin{equation}
    \varepsilon\int_{\Omega } \Psi^{\varepsilon}  \nabla \big(\eta_{\varepsilon}  S_{\varepsilon} (\nabla\widehat{w}_0) \big)\cdot \overline{\nabla \varphi} \,dx .
\end{equation}
Summing \eqref{termii}-\eqref{I1} up together, we get the desired conclusion.
\end{proof}

To present the main results of this section, we first give some notation and definitions. 

For $k>0$, we define the resolvent operator $R_{\varepsilon}(k)$ by 
\begin{equation}\label{defRp}
    R_{\varepsilon}(k):g \mapsto v_{\varepsilon},
\end{equation}
where $v_{\varepsilon}$ is the solution of
    \begin{equation}\label{defRp2}
        \left\{ 
        \begin{aligned}
           & (L_\varepsilon  +k^2) v_\varepsilon = g && \mathrm{in} \ \mathbb{R}^d, \\
           & v_{\varepsilon} \in \mathrm{SRC}(k). 
        \end{aligned}
        \right.
    \end{equation}

For any bounded Lipschitz domain $U\subset \mathbb{R}^d$ such that $\overline{\Omega}\subset U$, we define the Hilbert space $\mathcal{H}_{\varepsilon}(U)$ as the Sobolev space $H^1(U)$ equipped with the following norm:
\begin{equation}\label{weightednorm}
    \| u \|_{\mathcal{H}_{\varepsilon}(U)} := \| u \|_{L^2(U)} + \| \nabla u\|_{L^2(U\setminus D_{\varepsilon})}  + \varepsilon \|\nabla u \|_{L^2(D_{\varepsilon})} .
\end{equation}

\begin{lemma}\label{lemRLH}
    Let $k>0$. Let $r>0$ such that $\overline{\Omega} \subset B_r$. We have
    \begin{equation}\label{revoHL}
        \|R_{\varepsilon}(k)\|_{L^2(B_r)\rightarrow \mathcal{H}_{\varepsilon}(B_r)} \leq C(1+\|R_{\varepsilon}(k)\|_{L^2(B_r)\rightarrow L^2(B_r)} ),
    \end{equation}
where $C>0$ is a constant that depends only on $d$ and $k$.
\end{lemma}
\begin{proof}
    For any $g \in L^2(B_r) $, let $ v_\varepsilon =R_\varepsilon(k)g$. Using integration by parts, we get
\begin{equation}
\begin{aligned}
    & \int_{B_r\setminus \overline{D_{\varepsilon}}} |\nabla v_{\varepsilon}|^2\,dx + \varepsilon^2 \int_{D_{\varepsilon}} |\nabla v_{\varepsilon}|^2\,dx \\
    &  = k^2 \int_{B_r}|v_{\varepsilon}|^2\,dx + \mathrm{Re} \int_{\partial B_r} \overline{v_\varepsilon} \frac{\partial v_\varepsilon }{\partial \mathbf{n}} \, dS(x) - \mathrm{Re} \int_{B_r} g \overline{ v_\varepsilon}\,dx.
\end{aligned}
\end{equation}
Using \eqref{Renega}, we get
\begin{equation}
    \int_{B_r\setminus \overline{D_{\varepsilon}}} |\nabla v_{\varepsilon}|^2\,dx + \varepsilon^2 \int_{D_{\varepsilon}} |\nabla v_{\varepsilon}|^2\,dx \leq (1+k^2)  \int_{B_r}|v_{\varepsilon}|^2\,dx + \int_{B_r}|g|^2\,dx .
\end{equation}
We conclude that
\begin{equation}
    \| v_{\varepsilon} \|_{\mathcal{H}_{\varepsilon}(B_r)} \leq C(1+k^2)\| v_{\varepsilon} \|_{L^2(B_r)}  + C\|g \|_{L^2(B_r)}.
\end{equation}
This implies \eqref{revoHL}.
\end{proof}

\begin{theorem}
\label{L2velocity}
Let $k\in (0,\infty) \setminus \Sigma_D$, and $r>0$ such that $\overline{\Omega} \subset B_r$. Let $ u_\varepsilon $ be the solution of the scattering problem \eqref{maineq}, and $\widehat{u}_0$ be the solution of the homogenized scattering problem \eqref{equ0}. Then, we have
\begin{equation}
        \| u_\varepsilon - \Lambda_\varepsilon \widehat{u}_0 -\varepsilon \chi^{\varepsilon}   \cdot\eta_{\varepsilon} S_\varepsilon ( \nabla \widehat{u}_0 )  \|_{L^2(B_r)}  \leq C\mathcal{E}(\widehat{u}_0)  (1+\| R_{\varepsilon}(k) \|_{L^2(B_r) \rightarrow L^2(B_r)} ),
\end{equation}
where $C>0$ is a constant that depends only on $d,k,D$ and $\Omega$.
\end{theorem}
\begin{proof}
For any $g \in L^2(B_r) $, let $v_\varepsilon =R_{\varepsilon}(k)g$. Denote 
\begin{equation}
     z_\varepsilon := u_\varepsilon - \Lambda_\varepsilon \widehat{u}_0 -\varepsilon \chi^{\varepsilon}  \cdot\eta_{\varepsilon} S_\varepsilon ( \nabla \widehat{u}_0 ) .
\end{equation}
Using integration by parts and \eqref{smiden}, we have
\begin{equation}
    \int_{B_r} z_\varepsilon g\, dx= \big\langle (L_{\varepsilon}+k^2)z_{\varepsilon}, \overline{v_{\varepsilon}} \big\rangle_{H^{-1}(B_r),H^1(B_r)}   .
\end{equation}

We now apply Lemma \ref{basiclemmasca} for $U=B_r$, $\psi = \overline{v_{\varepsilon}} |_{B_r}$ and $\varphi = \mathcal{P}_{\varepsilon} (\overline{v_{\varepsilon}}|_{B_r\setminus \overline{D_{\varepsilon}}} )$. For simplicity, we denote $\mathcal{P}_{\varepsilon} (v_{\varepsilon}|_{B_r\setminus \overline{D_{\varepsilon}}} )$ by $\mathcal{P}_{\varepsilon} v_{\varepsilon}$ below. By Cauchy-Schwarz inequality and Lemma \ref{extensionoperator}, we get
\begin{equation}\label{Mep0}
\begin{aligned}
    \left|\int_{B_r} z_\varepsilon g\, dx \right| &\leq C  \mathcal{M}(\widehat{u}_0) \| v_{\varepsilon} \|_{\mathcal{H}_{\varepsilon}(\Omega)} + C \left| \int_{\Omega}  (\mu^k_0 -\Lambda_{\varepsilon} )\widehat{u}_0  \mathcal{P}_{\varepsilon}v_{\varepsilon} \,dx\right| \\
    &\quad\,+\left|\varepsilon \int_{D_{\varepsilon}}  (\nabla  \Lambda)^{\varepsilon}\widehat{u}_0   \cdot \nabla (v_{\varepsilon} - \mathcal{P}_{\varepsilon} v_{\varepsilon})  \,dx -\int_{D_{\varepsilon}} k^2\Lambda^{\varepsilon} \widehat{u}_0 (v_{\varepsilon} - \mathcal{P}_{\varepsilon} v_{\varepsilon}) \,dx \right| ,
\end{aligned}
\end{equation}
where
\begin{equation}\label{Mep1}
\begin{aligned}
    \mathcal{M}(\widehat{u}_0) & := \varepsilon \|\Psi^{\varepsilon} \nabla \big(\eta_{\varepsilon}  S_{\varepsilon} (\nabla\widehat{u}_0) \big) \|_{L^2(\Omega)} + \varepsilon \| \chi^{\varepsilon}  \nabla\big( \eta_{\varepsilon}  S_{\varepsilon} (\nabla\widehat{u}_0) \big) \|_{L^2(\Omega)}+ \varepsilon \| \chi^\varepsilon \cdot\eta_\varepsilon S_{\varepsilon} (\nabla\widehat{u}_0) \|_{L^2(\Omega)} \\
    &\quad\ \, +\varepsilon \| \nabla \Lambda \|_{L^{\infty}(D)} \| \widehat{u}_0 \|_{L^2(\Omega)} +(1+\|\nabla \chi\|_{L^2(Y)}) \|\nabla\widehat{u}_0 - \eta_{\varepsilon}  S_{\varepsilon} (\nabla\widehat{u}_0)  \|_{L^2(\Omega)}  \\
    & \quad\ \,  + \varepsilon^2 \|( \Lambda I +  \nabla \chi )^{\varepsilon} \eta_{\varepsilon} S_{\varepsilon} (\nabla\widehat{u}_0) \|_{L^2(\Omega)} +\varepsilon \| \Lambda_{\varepsilon} \nabla\widehat{u}_0 +\nabla\big( \varepsilon \chi^{\varepsilon} \cdot\eta_{\varepsilon}  S_{\varepsilon} (\nabla\widehat{u}_0) \big)\|_{L^2(D_{\varepsilon})} .
\end{aligned}
\end{equation}
The rest proof are divided into three steps.

\emph{Step 1}. By the choice of $\eta_{\varepsilon}$ and Lemma \ref{Svarepsilonproperty1}, we have
\begin{equation}\label{Mep2}
\begin{aligned}
    &\varepsilon \|\Psi^{\varepsilon} \nabla \big(\eta_{\varepsilon}  S_{\varepsilon} (\nabla\widehat{u}_0) \big) \|_{L^2(\Omega)} \\
    & \leq C \varepsilon\|\Psi^{\varepsilon} \big(\nabla \eta_{\varepsilon} \otimes S_{\varepsilon} (\nabla\widehat{u}_0) \big) \|_{L^2(\mathrm{supp}\,\nabla \eta_{\varepsilon})} + C\varepsilon \|\Psi^{\varepsilon} S_{\varepsilon} (\nabla^2 \widehat{u}_0)  \|_{L^2(\mathrm{supp}\, \eta_{\varepsilon})}  \\
    &\leq C\| \Psi\|_{L^2(Y)}\mathcal{E}(\widehat{u}_0).
\end{aligned}
\end{equation}
The same reasoning yields 
\begin{equation}\label{Mep3}
    \varepsilon \|\chi^{\varepsilon} \nabla \big(\eta_{\varepsilon}  S_{\varepsilon} (\nabla\widehat{u}_0) \big) \|_{L^2(\Omega)}  \leq C\| \chi\|_{L^2(Y)} \mathcal{E}(\widehat{u}_0) ,
\end{equation}
and
\begin{equation}
  \varepsilon \| \chi^\varepsilon \cdot\eta_\varepsilon S_{\varepsilon} (\nabla\widehat{u}_0) \|_{L^2(\Omega)} \leq C\|\chi\|_{L^2(Y)} \mathcal{E}(\widehat{u}_0) ,
\end{equation}
\begin{equation}\label{Mep4}
    \varepsilon^2 \|( \Lambda I +  \nabla \chi )^{\varepsilon}  \eta_{\varepsilon} S_{\varepsilon} (\nabla\widehat{u}_0) \|_{L^2(\Omega)}  \leq C \varepsilon   \| \Lambda +  \nabla \chi  \|_{L^2(Y)}  \mathcal{E}(\widehat{u}_0) .
\end{equation}

By Lemma \ref{propertySepsilondifference}, we have
\begin{equation}\label{Mep5}
\begin{aligned}
    \|\nabla\widehat{u}_0 - \eta_{\varepsilon}  S_{\varepsilon} (\nabla\widehat{u}_0)  \|_{L^2(\Omega)} &\leq C \|\nabla\widehat{u}_0 - \eta_{\varepsilon} \nabla\widehat{u}_0 \|_{L^2(\Omega)} + C\| \eta_{\varepsilon} \nabla\widehat{u}_0 - \eta_{\varepsilon}  S_{\varepsilon} (\nabla\widehat{u}_0)  \|_{L^2(\Omega)} \\
    &\leq C \mathcal{E}(\widehat{u}_0).
\end{aligned}
\end{equation}

For the last term on the right-hand side of \eqref{Mep1}, we have
\begin{equation}\label{Mep6}
    \varepsilon \|\Lambda_{\varepsilon} \nabla\widehat{u}_0 +\nabla\big( \varepsilon \chi^{\varepsilon} \cdot\eta_{\varepsilon}  S_{\varepsilon} (\nabla\widehat{u}_0) \big)\|_{L^2(D_{\varepsilon})}  \leq C  (\|\Lambda\|_{L^{\infty}(Y)} + \|\chi\|_{H^1(Y)} ) \mathcal{E}(\widehat{u}_0).
\end{equation}
Combining \eqref{Mep1}-\eqref{Mep6} together yields
\begin{equation}\label{ME}
    \mathcal{M}(\widehat{u}_0)  \leq C ( 1 + \|\Psi \|_{L^2(Y)} +   \|\Lambda\|_{W^{1,\infty}(Y)} + \|\chi\|_{H^1(Y)} ) \mathcal{E}(\widehat{u}_0) .
\end{equation}

\emph{Step 2}. For the second term on the right-hand side of \eqref{Mep0}, recall the definition \eqref{jepsi} of $J_{\varepsilon}$, we have
\begin{equation}\label{Stack}
\begin{aligned}
     \left| \int_{\Omega}  (\mu^k_0 -\Lambda_{\varepsilon} )\widehat{u}_0  \mathcal{P}_{\varepsilon}v_{\varepsilon} \,dx\right|&\leq  \left| \sum_{\mathbf{m}\in J_\varepsilon} \int_{\varepsilon (Y+\mathbf{m})}  (\mu^k_0 -\Lambda_{\varepsilon} )\widehat{u}_0  \mathcal{P}_{\varepsilon}v_{\varepsilon} \,dx\right| \\
     &\quad\,+   C\varepsilon^{1/2}\| \Lambda \|_{L^\infty(Y ) } \|\widehat{u}_0\|_{H^1(\Omega)}  \| v_\varepsilon \|_{H^1( \Omega \setminus \overline{ D_{\varepsilon}}  ) } ,
\end{aligned}
\end{equation}
where we used Lemma \ref{blelem} for $\widehat{u}_0$. Since $\mu^k_0 = \int_Y \Lambda(y)\,dy$, we have
\begin{equation}
\begin{aligned}
& \left| \sum_{\mathbf{m}\in J_\varepsilon} \int_{\varepsilon (Y+\mathbf{m})}  (\mu^k_0 -\Lambda_{\varepsilon} )\widehat{u}_0  \mathcal{P}_{\varepsilon}v_{\varepsilon} \,dx\right| \\
&=\left| \sum_{\mathbf{m}\in J_\varepsilon} \int_{\varepsilon (Y+\mathbf{m})}  (\mu^k_0 -\Lambda_{\varepsilon} ) \left(\widehat{u}_0  \mathcal{P}_{\varepsilon}v_{\varepsilon} - \varepsilon^{-d} \int_{\varepsilon (Y+\mathbf{m})} ( \widehat{u}_0  \mathcal{P}_{\varepsilon}v_{\varepsilon})(y)\,dy \right) \,dx\right|   \\
&\leq C \varepsilon\|\Lambda\|_{L^{\infty}(Y)} \| \nabla ( \widehat{ u}_0 \mathcal{ P } v_\varepsilon  ) \|_{ L^1(\Omega) }\\
& \leq C \varepsilon\| \Lambda \|_{L^\infty (Y ) } \| \widehat{ u }_0 \|_{H^1( \Omega )} \|v_\varepsilon \|_{ H^1 (\Omega\setminus \overline{ D_{\varepsilon}} ) } ,
\end{aligned}
\end{equation}
where we used the (1,1)-Poincar\'{e} inequality for $\widehat{ u}_0 \mathcal{ P } v_\varepsilon $.

For the last term on the right-hand side of \eqref{Mep0}, we have
\begin{equation}\label{Stack2}
    \begin{aligned}
        &\left|\varepsilon \int_{D_{\varepsilon}}  (\nabla  \Lambda)^{\varepsilon}\widehat{u}_0   \cdot \nabla (v_{\varepsilon} - \mathcal{P}_{\varepsilon} v_{\varepsilon})  \,dx -\int_{D_{\varepsilon}} k^2\Lambda^{\varepsilon} \widehat{u}_0 (v_{\varepsilon} - \mathcal{P}_{\varepsilon} v_{\varepsilon}) \,dx \right| \\
        & = \left|\varepsilon \int_{D_{\varepsilon}}  (\nabla  \Lambda)^{\varepsilon}\widehat{u}_0   \cdot \nabla (v_{\varepsilon} - \mathcal{P}_{\varepsilon} v_{\varepsilon})  \,dx - \varepsilon\int_{D_{\varepsilon}} (\nabla \Lambda)^{\varepsilon} \nabla(\widehat{u}_0 (v_{\varepsilon} - \mathcal{P}_{\varepsilon} v_{\varepsilon})) \,dx \right| \\
        & = \left| \varepsilon\int_{D_{\varepsilon}} (\nabla \Lambda)^{\varepsilon} \cdot \nabla\widehat{u}_0 (v_{\varepsilon} - \mathcal{P}_{\varepsilon} v_{\varepsilon}) \,dx \right| \\
        &\leq C\varepsilon \| \nabla \Lambda \|_{L^{\infty}(D)} \| \nabla \widehat{u}_0 \|_{L^2(\Omega)} \| v_{\varepsilon} \|_{\mathcal{H}_{\varepsilon}(\Omega)}  ,
    \end{aligned}
\end{equation}
where we used the equation $(\Delta +k^2)\Lambda =0$ in $D$ for the first equality and the Poincar\'{e} inequality for $v_{\varepsilon} - \mathcal{P}_{\varepsilon} v_{\varepsilon} \in H_0^1(D_{\varepsilon})$ in the last inequality.

\emph{Step 3}. Combining estimates \eqref{Mep0} and \eqref{ME}-\eqref{Stack2}  above, we get
\begin{equation}
\left|\int_{B_r} z_\varepsilon g\, dx \right| \leq  C\big( 1 + \|\Psi \|_{L^2(Y)} +   \|\Lambda\|_{W^{1,\infty}(Y)} + \|\chi\|_{H^1(Y)} \big) \mathcal{E}(\widehat{u}_0) \| v_{\varepsilon} \|_{\mathcal{H}_{\varepsilon}(\Omega)} .
\end{equation}
By Lemma \ref{thmflux}, $\|\Psi \|_{L^2(Y)} \leq C +C\|\chi\|_{H^1(Y)} \leq C$. By the definition \eqref{eq_lambda} of $\Lambda$, $\|\Lambda\|_{W^{1,\infty}(Y)} \leq C$, since $D$ is a $C^{1,\alpha}$ domain. Finally, the conclusion follows by Lemma \ref{lemRLH}.
\end{proof}

\subsection{\texorpdfstring{Weighted $H^1$ convergence rate}{Weighted H1 convergence rate}} This section is devoted to the $\mathcal{H}_{\varepsilon}$ convergence rate.
\begin{theorem}
\label{H1_velocity}
Under the same notations and conditions as Theorem \ref{L2velocity}. Then, there exist $\varepsilon_0,\,C>0$ that depend only on $d,k,D$ and $\Omega$ such that for any $0<\varepsilon<\varepsilon_0 $, we have
\begin{equation}
    \| u_\varepsilon - \Lambda_\varepsilon \widehat{u}_0 -\varepsilon \chi^{\varepsilon}  \cdot\eta_{\varepsilon} S_\varepsilon ( \nabla \widehat{u}_0 )  \|_{ \mathcal{H}_{\varepsilon}(B_r) } \leq C \mathcal{E}(\widehat{u}_0) (1+ \| R_{\varepsilon}(k) \|_{L^2(B_r)\rightarrow L^2(B_r)} ).
\end{equation}
\end{theorem}
\begin{proof}
Let $ z_\varepsilon := u_\varepsilon - \Lambda_\varepsilon \widehat{u}_0 -\varepsilon \chi^{\varepsilon}  \cdot \eta_{\varepsilon} S_\varepsilon ( \nabla \widehat{u}_0 ) $. By applying Lemma \ref{basiclemmasca} for $U= B_r$, $\psi = \varphi = \mathcal{P}_{\varepsilon} z_{\varepsilon} $, and using a similar (and simpler) argument to the proof of Theorem \ref{L2velocity}, we get
\begin{equation}\label{H1esti1}
 \left|\int_{B_r} A_{\varepsilon} \nabla z_{\varepsilon} \cdot \overline{\nabla \mathcal{P}_{\varepsilon} z_{\varepsilon} } \,dx - k^2 \int_{B_r} z_{\varepsilon} \overline{ \mathcal{P}_{\varepsilon} z_{\varepsilon} } \,dx - \int_{\partial B_r} \frac{\partial z_{\varepsilon}}{\partial \mathbf{n}} \overline{ z_{\varepsilon}}\,dS(x) \right| \leq C \mathcal{E}(\widehat{u}_0) \|z_{\varepsilon} \|_{H^1(\Omega\setminus D_{\varepsilon})} .
\end{equation}
Taking the real part of the left-hand side and using \eqref{Renega}, it follows that
\begin{equation}\label{353}
    \| \nabla z_{\varepsilon} \|_{L^2(B_r\setminus D_{\varepsilon})}  \leq C \mathcal{E}(\widehat{u}_0) +  C\| z_{\varepsilon} \|_{L^2(B_r)} +  C\varepsilon^2  \|\nabla z_{\varepsilon} \|_{L^2(D_{\varepsilon})}  .
\end{equation}
Therefore, by Theorem \ref{L2velocity}, we have
\begin{equation}\label{H1es2}
    \| z_{\varepsilon} \|_{\mathcal{H}_{\varepsilon}(B_r)}  \leq C \mathcal{E}(\widehat{u}_0) (1+ \| R_{\varepsilon}(k) \|_{L^2(B_r) } ) + C \varepsilon  \|\nabla z_{\varepsilon} \|_{L^2(D_{\varepsilon})} .
\end{equation}

We estimate the last term $\varepsilon \|\nabla z_{\varepsilon}\|_{L^2(D_{\varepsilon})}$ on the right-hand side of \eqref{H1es2} as follows. Denote $u_{1,\varepsilon} =\chi^{\varepsilon} \cdot \eta_{\varepsilon} S_\varepsilon ( \nabla \widehat{u}_0 )  $, then $z_{\varepsilon} = u_{\varepsilon} - \Lambda_{\varepsilon} \widehat{u}_0 - \varepsilon u_{1,\varepsilon}$. By the proof of Theorem \ref{L2velocity}, it is clear that
\begin{equation}\label{355}
    \| u_{1,\varepsilon} \|_{L^2(\Omega)} +\varepsilon \| \nabla u_{1,\varepsilon} \|_{L^2(\Omega)}   \leq C \varepsilon^{-1/2}\mathcal{E}(\widehat{u}_0) .
\end{equation}
Let $\zeta = z_\varepsilon-\mathcal{P}_{\varepsilon} z_\varepsilon \in H^1_0(D_\varepsilon)$. A direct computation yields
\begin{equation}\label{356}
    \begin{aligned}
        \varepsilon^2\int_{ D_\varepsilon} \nabla z_\varepsilon\cdot \overline{\nabla\zeta } \, dx 
        &=  k^2 \int_{D_\varepsilon}   (z_\varepsilon+\varepsilon u_{1,\varepsilon} ) \overline{\zeta} \,dx + k^2  \int_{D_{\varepsilon}} \Lambda^\varepsilon \widehat{ u}_0 \overline{\zeta} \,dx \\
        & \quad\,- \varepsilon \int_{D_{\varepsilon}} (\nabla \Lambda)^\varepsilon \widehat{u}_0 \cdot \overline{\nabla\zeta } \,dx -\varepsilon^2 \int_{D_{\varepsilon}}  (\Lambda^\varepsilon \nabla \widehat{u}_0+\varepsilon \nabla u_{1,\varepsilon} )  \cdot \overline{\nabla\zeta } \, dx \\
        & = k^2 \int_{D_\varepsilon}   (z_\varepsilon+\varepsilon u_{1,\varepsilon} ) \overline{\zeta} \,dx 
        + \varepsilon  \int_{D_\varepsilon}  (\nabla \Lambda)^{\varepsilon} \cdot  \nabla \widehat{u}_0 \overline{\zeta}\, dx \\
        &\quad\,-\varepsilon^2 \int_{D_{\varepsilon}}  (\Lambda^\varepsilon \nabla \widehat{u}_0+\varepsilon \nabla u_{1,\varepsilon} )  \cdot \overline{\nabla\zeta } \, dx,
    \end{aligned}
\end{equation}
where the last equality follows from the equation $(\Delta +k^2)\Lambda =0$ in $D$. Using \eqref{355} and the Poincar\'{e} inequality $\|\zeta \|_{L^2(D_\varepsilon)} \leq C\varepsilon \left\|\nabla \zeta \right\|_{L^2(D_\varepsilon)} $, we obtain
\begin{equation}\label{357}
\begin{aligned}
    \| \varepsilon \nabla z_{\varepsilon} \|^2_{L^2(D_{\varepsilon})} &\leq C \| \varepsilon \nabla z_{\varepsilon} \|_{L^2(D_{\varepsilon})} \| \varepsilon \nabla z_{\varepsilon} \|_{L^2(B_r\setminus \overline{D_{\varepsilon}})} \\
    &\quad\, + C \mathcal{E}(\widehat{u}_0) (1+ \| R_{\varepsilon}(k) \|_{L^2(B_r) \rightarrow \mathcal{H}_{\varepsilon}(\Omega)} ) ( \| \varepsilon \nabla z_{\varepsilon} \|_{L^2(D_{\varepsilon})} + \| \varepsilon \nabla z_{\varepsilon} \|_{L^2(B_r\setminus \overline{D_{\varepsilon}})}  ).
\end{aligned}
\end{equation}

Denote $X = \| \varepsilon \nabla z_{\varepsilon} \|_{L^2(D_{\varepsilon})}$ and $H = \mathcal{E}(\widehat{u}_0) (1+ \| R_{\varepsilon}(k) \|_{L^2(B_r) \rightarrow \mathcal{H}_{\varepsilon}(\Omega)} ) $. Then, \eqref{353} and \eqref{357} imply that
\begin{equation}
   X^2 \leq C( \varepsilon^2 X^2 +HX+ H^2).
\end{equation}
Therefore, there exists $\varepsilon_0>0$ that depends only on $d,k,\Omega,D$, such that for any $0<\varepsilon <\varepsilon_0$, we have $X\leq C H$.
The proof is complete with Lemma \ref{lemRLH}.
\end{proof}

It follows from Theorems \ref{L2velocity} and \ref{H1_velocity} that, to complete the proof of Theorem \ref{mainresult1}, it remains to establish Theorem \ref{blethm} (Section \ref{secbound}) and to establish Theorem \ref{thm5121main} (Section \ref{secresout}).

\section{Boundary layer estimate for the homogenized scattering problem}\label{secbound}

The goal of this section is to estimate
\begin{equation}
    \mathcal{E}(\widehat{u}_0) = \varepsilon^{1/2} \| \widehat{u}_0 \|_{H^1(\Omega)} + \| \nabla \widehat{u}_0 \|_{L^2(\mathcal{O}_{(d+2)\varepsilon})} + \varepsilon \| \nabla^2 \widehat{u}_0 \|_{L^2(\Omega \setminus \overline{\mathcal{O}_{(d-1)\varepsilon}})},
\end{equation}
where $\widehat{u}_0$ solves the homogenized scattering problem \eqref{equ0}, which should be viewed as the scattering problem by a penetrable anisotropic obstacle (possibly with negative refraction, depending on the sign of $\mu^k_0$) with a Lipschitz boundary.

\subsection{Existence of weak solutions}

Let $k>0$, $s\in \mathbb{R}$, and $A $ be a symmetric, positive-definite, real-valued $d\times d$ matrix. We consider the following problem
    \begin{equation}\label{4231}
        \left\{\begin{aligned}
           & (\Delta+k^2)u=0 && \mathrm{in}\ \R^d \setminus \overline{\Omega} ,\\
        &( \nabla\cdot A\nabla + s )u =0 &&  \mathrm{in}\ \Omega, \\
        &u|_{+} = u|_{-} && \mathrm{on}\ \partial \Omega ,\\
        &\frac{\partial u}{\partial \mathbf{n}} \Big|_{+} = \frac{\partial u}{\partial \nu_A} \Big|_{-}  && \mathrm{on}\ \partial \Omega, \\
        & u - \uinc\in \mathrm{SRC}(k).
        \end{aligned}\right.
    \end{equation}
Note that $\widehat{u}_0$ solves problem \eqref{4231} with $A=A_0$, $s = k^2\mu^k_0$ and the same choice of $\uinc$. 

Fix $r>0$ such that $\overline{\Omega} \subset B_r$. Define the Dirichlet-to-Neumann (DtN) map $\Lambda$ by
\begin{equation}
    \Lambda: v\mapsto \frac{\partial v}{\partial \mathbf{n}} \Big|_{\partial B_r} ,
\end{equation}
where $v$ is the solution of 
    \begin{equation}\label{scaprov}
        \left\{
        \begin{aligned}
            &(\Delta +k^2)v =0 && \mathrm{in}\ \mathbb{R}^d\setminus \overline{B_r}, \\
            &v\in \mathrm{SRC}(k).
        \end{aligned}\right.
    \end{equation}
Let $w:= u - \uinc$, with the help of the DtN map, the problem \eqref{4231} is equivalent to the following variational problem on $B_r$ \cite[Page 26]{cakoni2022inverse}: find $w\in H^1(B_r)$ such that
\begin{equation}\label{varian}
    a_1(w,v) + a_2(w,v) = H(v) \qquad \forall v \in H^1(B_r),
\end{equation}
where
\begin{equation}
\begin{aligned}
    a_1(\phi,v) &:= \int_{\Omega} A\nabla \phi \cdot \overline{\nabla v}\,dx + \int_{B_r\setminus \overline{\Omega}} \nabla \phi\cdot \overline{\nabla v}\,dx \\
    &\quad\, + \int_{B_r} \phi\overline{v}\,dx - \langle \Lambda \phi, v\rangle_{H^{-1/2}(\partial B_r),H^{1/2}(\partial B_r)} ,
\end{aligned}
\end{equation}
\begin{equation}
    a_2 (\phi,v) := - (s+1) \int_{\Omega} \phi\overline{v}\,dx  - (k^2+1)\int_{B_r\setminus \overline{\Omega}}  \phi\overline{v}\,dx ,
\end{equation}
\begin{equation}
    H(v):= \int_{\Omega} (I-A) \nabla \uinc \cdot \overline{\nabla v}\,dx +(s-k^2)\int_{\Omega} \uinc \overline{v}\,dx .
\end{equation}

\begin{theorem}\label{weakso}
Let $k>0$, $s\in \mathbb{R}$, and $A $ be a symmetric, positive-definite, real-valued $d\times d$ matrix. There exists a unique solution $u \in H^1_{\mathrm{loc}}(\mathbb{R}^d)$ to problem \eqref{4231}. Moreover, for any $r>0$ such that $\overline{\Omega} \subset B_r$, we have 
\begin{equation}
    \| u \|_{H^1(B_r)} \leq C\|\uinc\|_{H^1(B_r)},
\end{equation}
where $C>0$ is a constant depending only on $d,k,s,r, A$ and $\Omega$.
\end{theorem}
\begin{proof}
The proof is divided into three steps.

\emph{Step 1}. We consider $a_1$ as a bilinear form on $H^1(B_r) \times H^1(B_r)$. Since $A $ is positive-definite and real-valued, $a_1$ is coercive:
    \begin{equation}
        \mathrm{Re}\,a_1(v,v) \geq C\| v \|_{H^1(B_r)}^2 - \mathrm{Re}\,\langle \Lambda v, v\rangle_{H^{-1/2}(\partial B_r),H^{1/2}(\partial B_r)} \geq C\| v \|_{H^1(B_r)}^2,
    \end{equation}
where we used \eqref{Renega} for the second inequality. On the other hand, $a_1$ is clearly bounded. Hence, by the Lax-Milgram theorem, there exists an invertible bounded operator $T$ on $H^1(B_r)$ such that $a_1(\phi,v) = \langle T\phi, v\rangle_{H^1(B_r)}$. 

For any $\phi \in L^2(B_r)$, $a_2(\phi,\cdot):H^1(B_r)\rightarrow \mathbb{C}$ is a bounded linear functional. By the Riesz representation theorem, there exists $\phi'\in H^1(B_r)$ such that $a_2(\phi,v) = \langle \phi'  , v \rangle_{H^1(B_r)}$ for any $v\in H^1(B_r)$. We define the operator $K:L^2(B_r)\rightarrow L^2(B_r)$ by $K\phi =  \phi'$.
Note that
\begin{equation}\label{Kphiphi}
    \| K\phi \|_{H^1(B_r)} \leq C(1+|s|+k^2)\|\phi \|_{L^2(B_r)},
\end{equation}
which comes from $\| K\phi \|_{H^1(B_r)}^2 =a_2(\phi,K\phi) \leq C(1+|s|+k^2)\|\phi \|_{L^2(B_r)} \|K\phi \|_{L^2(B_r)}$. 

Similarly, $H:v\mapsto H(v)$ is a bounded linear functional on $H^1(B_r)$. By the Riesz representation theorem, there exists $h\in H^1(B_r)$ such that $H(v) = \langle h , v \rangle_{H^1(B_r)}$ for any $v\in H^1(B_r)$. 

Therefore, the problem \eqref{4231} is equivalent to solve 
\begin{equation}\label{in42}
    (I+T^{-1}K) w= T^{-1}h.
\end{equation}

\emph{Step 2}. We now consider $T^{-1}K$ as an operator on $L^2(B_r)$. Assuming that $I+T^{-1}K$ is invertible, then \eqref{in42} has a unique solution $w \in L^2(B_r)$ and $\| w\|_{L^2(B_r)} \leq C\|h \|_{H^1(B_r)} $. Thus,
\begin{equation}
\begin{aligned}
     \|w \|_{H^1(B_r)}&=\|T^{-1}(h - K w) \|_{H^1(B_r)} \\
     &\leq C( \|h \|_{H^1(B_r)} + \|Kw \|_{H^1(B_r)} ) \\
     &\leq C( \|h \|_{H^1(B_r)} + \|w \|_{L^2(B_r)} ) \\
     &\leq C \|h \|_{H^1(B_r)} .
\end{aligned}
\end{equation}
This completes the proof by noting $u = w+\uinc$ and $\|h \|_{H^1(B_r)} \leq C\| \uinc\|_{H^1(\Omega)}$.

\emph{Step 3}. We now prove $I+T^{-1}K$ is invertible. It follows from \eqref{Kphiphi} that
\begin{equation}
    \|T^{-1 } K \phi \|_{H^1(B_r)} \leq C\|  K \phi \|_{H^1(B_r)} \leq C\|\phi \|_{L^2(B_r)}.
\end{equation}
By the Rellich–Kondrachov theorem, $T^{-1}K$ is compact. By the Fredholm theory, $I+T^{-1}K$ is a Fredholm operator of index zero. Therefore, it suffices to show that $I+T^{-1}K$ is injective. Assume that $(I+T^{-1}K)\widetilde{u}=0$. By the equivalence of problems \eqref{4231} and \eqref{varian}, it turns out that 
\begin{equation}
    u=\left\{ 
    \begin{aligned}
        & \widetilde{u} && \mathrm{in}\ B_r, \\
        & \textrm{solution to problem \eqref{scaprov} with } v|_{\partial\Omega} = \widetilde{u}|_{\partial \Omega}  && \mathrm{in}\ \mathbb{R}^d\setminus \overline{B_r}
    \end{aligned}\right.
\end{equation}
solves the problem \eqref{4231} with $\uinc=0$. By integration by parts, we get
\begin{equation}
\begin{aligned}
    \int_{\partial B_r} \frac{\partial u}{\partial \mathbf{n}} (x)\overline{u(x)}\,dS(x) & = \int_{\Omega}A\nabla u(x)\cdot \overline{\nabla u(x)}\,dx + \int_{B_r\setminus \overline{\Omega}}  |\nabla u(x)|^2\,dx \\
    & \quad \, - s \int_{\Omega} |u(x)|^2\,dx -  k^2\int_{B_r\setminus \overline{\Omega}} |u(x)|^2\,dx.
\end{aligned}
\end{equation}
Taking the imaginary part of the above equality, we obtain
\begin{equation}
    \mathrm{Im}\,\int_{\partial B_r} \frac{\partial u}{\partial \mathbf{n}} (x)\overline{u(x)}\,dS(x) =0.
\end{equation}
It follows from Lemma \ref{rellichunique} (the Rellich uniqueness theorem) that $u\equiv 0$ in $\mathbb{R}^d\setminus \overline{B_r}$. By the unique continuity principle, $u\equiv 0$ in $\mathbb{R}^d\setminus \overline{\Omega}$. Since $u|_{\partial \Omega} = \frac{\partial u}{\partial \nu_A}|_{\partial \Omega} =0$, $u$ can be extended by zero over a larger domain such that the extension still satisfies the equation $\nabla\cdot(A\nabla u)+su=0$ on the larger domain. Using the unique continuity principle again, we obtain $u\equiv 0$ in $\mathbb{R}^d$.
\end{proof}

\subsection{Layer potentials}
The idea is to represent $\widehat{u}_0$ via layer potentials, with the associated density determined by a boundary integral equation. The desired estimate for $\mathcal{E}(\widehat{u}_0) $ is given by the quantitative estimates of layer potentials.

We begin by recalling basic facts from layer potential theory; see \cite{ammari2009layer} and the references therein for more details. The fundamental solution for the Laplacian is
\begin{equation}\label{funwhole}
    \Gamma(x): = \left\{
    \begin{aligned}
        & \frac{\mathrm{1}}{2\pi  } \log |x|, && d=2, \\
        &\frac{1}{(2-d)\omega_d}|x|^{2-d}, && d\geq 3,
    \end{aligned}\right.
\end{equation}
where $\omega_d$ denotes the area of the unit sphere in $\mathbb{R}^d$. Let $A $ be a symmetric, positive-definite, real-valued $d\times d$ matrix. Define
\begin{equation}\label{fundHel}
    \Gamma_A(x) = |\det A|^{-d/2} \Gamma(A^{-1/2}x).
\end{equation}
We have $\nabla \cdot A\nabla \Gamma_A = \delta_0$, where $\delta_0$ denotes the Dirac delta distribution.

We introduce the single layer potential operator
\begin{equation}\label{single}
    \mathcal{S}_A(f)(x):= \int_{\partial \Omega} \Gamma_A(x-y) f(y)\,dS(y), \qquad x\in \mathbb{R}^d\setminus \partial \Omega,
\end{equation}
and the boundary version of \eqref{single}, i.e.,
\begin{equation}
    S_Af(x) :=\int_{\partial \Omega} \Gamma_A(x-y) f(y)\,dS(y), \qquad x\in \partial \Omega.
\end{equation}
The Neumann-Poincar\'{e} operator is defined by
\begin{equation}
    K^{*}_A(f)(x) : = \mathrm{p.v.}\int_{\partial \Omega} \mathbf{n}(x)\cdot A\nabla \Gamma_A (x-y) f(y)\,dS(y), \qquad x\in \partial \Omega,
\end{equation}
where $\mathrm{p.v.}$ denotes the principal value. 

Fix $\alpha > 0$ and for each $x \in \partial \Omega$, we introduce the non-tangential approach regions
\begin{equation}
    \Gamma^+(x) := \Gamma^+_{\alpha}(x) := \{ y \in \mathbb{R}^d\setminus \overline{\Omega} : |x-y| < (1 + \alpha)\,\mathrm{dist}\,(y, \partial \Omega) \},
\end{equation}
and
\begin{equation}
    \Gamma^-(x) := \Gamma^-_{\alpha}(x) := \{ y \in \Omega : |x-y| < (1 + \alpha)\,\mathrm{dist}\,(y, \partial \Omega) \}.
\end{equation}
For $u : \mathbb{R}^d\setminus \partial\Omega \to \mathbb{R}$, we define the non-tangential maximal function of $u$ as follows:
\begin{equation}
      (u)^*_{\pm}(x) := (u)^*_{\alpha,\pm}(x) := \sup \{ |u(y)| : y \in \Gamma^{\pm}_{\alpha}(x) \}, \quad x \in \partial \Omega.
\end{equation}
For each $\alpha,\beta>0$, the $L^2(\partial \Omega)$ norms of $(u)^*_{\alpha,\pm}$ and $(u)^*_{\beta,\pm}$
are comparable \cite[Proposition 2.1.2]{hofmann2010singular}, so the choice of $\alpha$ plays a relatively minor role when measuring the $L^2(\partial \Omega)$ norm of the non-tangential maximal function. Hence, we omit the subscript $\alpha$ below.

\begin{lemma}\label{lemnont1}
We have the following properties:
\begin{enumerate}
    \item For any $f \in L^2(\partial \Omega)$, $\nabla \cdot A\nabla \mathcal{S}_Af =0$ in $\mathbb{R}^d\setminus \partial \Omega$.
    \item \emph{(Jump relation)} For any $f \in L^2(\partial \Omega)$, we have 
    \begin{equation}\label{jump}
    \left. \frac{\partial \mathcal{S}_Af}{\partial \nu_A} \right|_{\pm}(x) = \left( \pm\frac{1}{2} I +K^{*}_A\right) f(x) \quad \textrm{a.e. }x\in\partial \Omega, 
\end{equation}
where the conormal derivative of $u$ is defined by
\begin{equation}
    \left. \frac{\partial u}{\partial \nu_A} \right|_{\pm}(x) :=\lim_{y\rightarrow x \atop y\in \Gamma^{\pm}(x)}\mathbf{n}(x)\cdot A\nabla u (y) , \quad \textrm{a.e. }x\in\partial \Omega.
\end{equation}
\item \emph{(Non-tangential maximal function estimates)} For any $f \in L^2(\partial \Omega)$, we have 
\begin{equation}\label{inenontan}
    \| (\nabla \mathcal{S}_Af)^*_{\pm} \|_{L^2(\partial \Omega)} \leq C \| f\|_{L^2(\partial \Omega)},
\end{equation}
where $C>0$ is a constant that depends only on 
\item For $d\geq 3$, $S_A:L^2(\partial \Omega) \rightarrow H^1(\partial \Omega)$ is invertible. 
\item For $d\geq2$, the operator $B:L^2(\partial \Omega) \times \mathbb{C}\rightarrow H^1(\partial \Omega)\times \mathbb{C}$ by
\begin{equation}
    B(f,a):= \left(S_Af+a, \int_{\partial \Omega} f(x)\,dS(x) \right)
\end{equation}
is invertible.
\end{enumerate}
\end{lemma}
\begin{proof}
    For the case of $A=I$, we refer to \cite{verchota1984layer} for items (1)-(3), and to \cite[Theorem 2.26]{ammari2007polarization} for items (4)-(5). The conclusion for the general matrix $A$ then follows from the relation \eqref{fundHel} between $\Gamma^A$ and $\Gamma$.
\end{proof}

For later use, we introduce the layer potentials on the ball $B_r$, $r>0$. The Green’s function of Laplacian on the ball $B_r$ is given by 
\begin{equation} \label{Gr}
G^r(x,y):= \left\{
    \begin{aligned}
        & \frac{1}{2\pi} \left( \log |x-y|
 - \log \left|\frac{r^2}{|x|}x - \frac{|x|}{r}y \right| \right) , && d=2, \\
 & \frac{1}{(2-d)\omega_d} \left( |x-y|^{2-d}
 - \left|\frac{r^2}{|x|}x - \frac{|x|}{r}y \right|^{2-d} \right), && d\geq 3.
    \end{aligned}\right.
\end{equation}
We have
\begin{equation}\label{proGr} \left\{
    \begin{aligned}
        & \Delta G^r(\cdot,y) = \delta_y && \mathrm{in}\ B_r,\\
 & G^r(\cdot,y) = 0 && \mathrm{on}\ \partial B_r.
    \end{aligned}\right.
\end{equation}
Similar to the case of the whole space, we introduce the single layer potential operator
\begin{equation}\label{single1}
    \mathcal{S}^r(f)(x):= \int_{\partial \Omega} G^r(x,y) f(y)\,dS(y), \qquad x\in B_r\setminus \partial \Omega,
\end{equation}
and the boundary version of \eqref{single1}, i.e.,
\begin{equation}
    S^rf(x) :=\int_{\partial \Omega} G^r(x,y) f(y)\,dS(y), \qquad x\in \partial \Omega.
\end{equation}
The Neumann-Poincar\'{e} operator is defined by
\begin{equation}
    K^{r,*}(f)(x) : = \mathrm{p.v.}\int_{\partial \Omega} \mathbf{n}(x)\cdot \nabla G^r (x,y) f(y)\,dS(y),\qquad x\in \partial \Omega.
\end{equation}

\begin{lemma}\label{lemnont2}
We have the following properties:
\begin{enumerate}
    \item For any $f \in L^2(\partial \Omega)$, $\Delta \mathcal{S}^rf =0$ in $B_r\setminus \partial \Omega$ and $(\mathcal{S}^r f)|_{\partial B_r} =0$.
    \item \emph{(Jump relation)} For any $f \in L^2(\partial \Omega)$, we have 
    \begin{equation}\label{jump1}
    \left. \frac{\partial \mathcal{S}^rf}{\partial \mathbf{n}} \right|_{\pm}(x) = \left( \pm\frac{1}{2} I +K^{r,*} \right) f(x) \quad \textrm{a.e. }x\in\partial \Omega. 
\end{equation}
    \item For $d\geq 2$, $S^r:H^{-1/2}(\partial \Omega)\rightarrow H^{1/2}(\partial \Omega)$ is invertible.
\end{enumerate}
\end{lemma}
\begin{proof}
    Item (1) follows from the property \eqref{proGr} of the Green's function. For item (2), it follows from the definition \eqref{Gr} of $G^r$ that, for any $y\in B_r$, $\Gamma(\cdot - y) - G^r(\cdot,y)$ is smooth in $B_r$. This can also be observed from the equation 
    \begin{equation}
        \left\{\begin{aligned}
            & \Delta (\Gamma(\cdot - y) - G^r(\cdot,y)) = 0 && \mathrm{in} \ B_r, \\
            & \Gamma(\cdot - y) - G^r(\cdot,y) =\Gamma(\cdot - y) && \mathrm{on} \ \partial B_r, 
        \end{aligned} \right.
    \end{equation}
    and the elliptic regularity. Therefore, the jump relation \eqref{jump} for $A=I$ implies \eqref{jump1}.

    For item (3), by \cite[Theorem 7.6]{mclean2000strongly}, $S_I:H^{-1/2}(\partial \Omega)\rightarrow H^{1/2}(\partial \Omega)$ is a Fredholm operator of index zero. Since $\Gamma(\cdot - y) - G^r(\cdot,y)$ is smooth in $B_r$, $S^r$ is also a Fredholm operator of index zero. It suffices to show $S^r$ is injective. Assume that $S^r f=0$. By the jump relation \eqref{jump1}, we have
\begin{equation}
\begin{aligned}
    0 &= \int_{\partial \Omega} f (x)\overline{S^r f(x)} \,dS(x) \\
    &=\int_{\partial \Omega} \frac{\partial \mathcal{S}^rf}{\partial \mathbf{n}} \Big|_+(x) \overline{S^r f(x)}  \,dS(x) -\int_{\partial \Omega} \frac{\partial \mathcal{S}^rf}{\partial \mathbf{n}} \Big|_- (x) \overline{S^r f(x)}  \,dS(x) \\
    &= - \int_{B_r} |\nabla \mathcal{S}^rf|^2\,dx.
\end{aligned}
\end{equation}
This, combining with the Dirichlet boundary condition $\mathcal{S}^r f=0$ on $\partial B_r$, implies that $\mathcal{S}^r f\equiv 0$ in $B_r$. Thus, $f = \frac{\partial \mathcal{S}^rf}{\partial \mathbf{n}} \big|_+ - \frac{\partial \mathcal{S}^rf}{\partial \mathbf{n}} \big|_- =0$. The proof is complete.
\end{proof}

\subsection{Proof of Theorem \ref{blethm}: boundary layer estimates}

In this section, we prove Theorem \ref{blethm}. We first consider the following Dirichlet boundary value problem on $B_r$ with transmission conditions on the interface $\partial \Omega$:
\begin{equation}\label{42310}
        \left\{\begin{aligned}
           & \Delta v=0 && \mathrm{in}\ B_r \setminus \overline{\Omega} ,\\
        &\nabla\cdot A\nabla v =0 &&  \mathrm{in}\ \Omega, \\
        &v|_{+} - v|_{-}=\varphi && \mathrm{on}\ \partial \Omega ,\\
        &\frac{\partial v}{\partial \mathbf{n}} \Big|_{+} - \frac{\partial v}{\partial \nu_A} \Big|_{-} =  \psi  && \mathrm{on}\ \partial \Omega, \\
        & v=0 && \mathrm{on}\ \partial B_r.
        \end{aligned}\right.
    \end{equation}

\begin{theorem}\label{thm44}
   Suppose that $I-A$ is either positive or negative semi-definite. For any $(\varphi,\psi)\in H^1(\partial \Omega) \times L^2(\partial \Omega)$, the solution $v$ to the problem \eqref{42310} has the single-layer potential representation
    \begin{equation}\label{sinv}
        v(x) = \left\{
        \begin{aligned}
            & \mathcal{S}^rf(x), && x\in B_r\setminus \overline{\Omega}, \\
            & \mathcal{S}_A g(x) , && x\in \Omega,
        \end{aligned}\right.
    \end{equation}
where $(f,g)\in L^2(\partial \Omega) \times L^2(\partial \Omega)$ satisfies the boundary integral equation
\begin{equation}\label{sinv2}
    \mathcal{A} \begin{pmatrix}
        f \\ g
    \end{pmatrix}: = \begin{pmatrix}
        S^r & -S_A \\
        \frac{1}{2}I+K^{r,*} & \frac{1}{2}I -K^{*}_{A}
    \end{pmatrix} 
    \begin{pmatrix}
        f \\ g
    \end{pmatrix}
    =\begin{pmatrix}
        \varphi \\ \psi
    \end{pmatrix}.
\end{equation}
For $d=2$, we additionally require that $g \in L^2_0(\partial \Omega)$, where $L^2_0(\partial \Omega):= \{f \in L^2(\partial \Omega):\int_{\partial \Omega }f=0\}$. Moreover, we have
\begin{equation}\label{441fg}
    \|f \|_{L^2(\partial \Omega)} + \|g \|_{L^2(\partial \Omega)}  \leq C( \|\varphi \|_{H^1(\partial \Omega)}  + \|\psi \|_{L^2(\partial \Omega)} ),
\end{equation}
where $C>0$ is a constant that depends only on $d,r,A$ and $\Omega$.
\end{theorem}
\begin{proof}
    Using the jump relations \eqref{jump} and \eqref{jump1}, we find that the function $v$ defined by \eqref{sinv}-\eqref{sinv2} solves the problem \eqref{42310}. 
    
    \emph{Step 1}. For $d\geq 3$. It remains to show that $\mathcal{A}$ is an invertible operator from $L^2(\partial \Omega) \times L^2(\partial \Omega)$ to $H^1(\partial \Omega) \times L^2(\partial \Omega)$. Define the auxiliary operator
\begin{equation}
    \mathcal{A}' : = \begin{pmatrix}
        S_I & -S_A \\
        \frac{1}{2}I+K_I^{*} & \frac{1}{2}I -K^{*}_{A}
    \end{pmatrix} .
\end{equation}
It follows from \cite[Theorem 2]{escauriaza1993regularity} that, with the assumption that $I-A$ is either positive or negative semi-definite, $\mathcal{A}'$ is an invertible operator from $L^2(\partial \Omega) \times L^2(\partial \Omega)$ to $H^1(\partial \Omega) \times L^2(\partial \Omega)$. We note that
\begin{equation}
    \mathcal{A} = \mathcal{A}' + \mathcal{R}, \qquad \mathcal{R} :=\begin{pmatrix}
        S^r -S_I & 0 \\
        K^{r,*} -K^{*}_I & 0
    \end{pmatrix} .
\end{equation}
It follows from \eqref{funwhole} and \eqref{Gr} that $\mathcal{R}$ is a compact operator from $L^2(\partial \Omega) \times L^2(\partial \Omega)$ to $H^1(\partial \Omega) \times L^2(\partial \Omega)$. Consequently, $\mathcal{A}$ is a Fredholm operator of index zero. To complete the proof, it suffices to show $\mathcal{A}$ is injective.

Suppose that $\mathcal{A}(f,g)^{\top} =0$ for some $(f,g)^\top \in L^2(\partial \Omega) \times L^2(\partial \Omega)$, then
\begin{equation}\label{deux}
        u(x): = \left\{
        \begin{aligned}
            & \mathcal{S}^rf(x), && x\in B_r\setminus \overline{\Omega}, \\
            & \mathcal{S}_A g(x) , && x\in \Omega
        \end{aligned}\right.
    \end{equation}
solves the Dirichlet problem
\begin{equation}\label{4Diri}
        \left\{\begin{aligned}
           & \Delta u=0 && \mathrm{in}\ B_r \setminus \overline{\Omega} ,\\
        &\nabla\cdot A\nabla u =0 &&  \mathrm{in}\ \Omega, \\
        &u|_{+} = u|_{-} && \mathrm{on}\ \partial \Omega ,\\
        &\frac{\partial u}{\partial \mathbf{n}} \Big|_{+} = \frac{\partial u}{\partial \nu_A} \Big|_{-}   && \mathrm{on}\ \partial \Omega, \\
        & u=0 && \mathrm{on}\ \partial B_r.
        \end{aligned}\right.
    \end{equation}
It is clear that $u\equiv 0$ in $B_r$. In particular, $S^r f=S_A g=0$. By Lemma \ref{lemnont1}.(4) and Lemma \ref{lemnont2}.(3), we get $f=g=0$. Thus $\mathcal{A}$ is injective.

\emph{Step 2}. For $d=2$. We now consider $\mathcal{A}$ and $\mathcal{A}'$ as operators mapping from $L^2(\partial \Omega) \times L_0^2(\partial \Omega)$ to $H^1(\partial \Omega) \times L^2(\partial \Omega)$. By inspecting the proof of \cite[Theorem 2]{escauriaza1993regularity}, one can still show that, assuming $I-A$ is either positive or negative semi-definite, $\mathcal{A}'$ has a closed and dense range. To complete the proof, we repeat the argument of Step 1. It turns out that it remains to check the injectivity of $\mathcal{A}$ and $\mathcal{A}'$. We only check the injectivity of $\mathcal{A}$, since the same reasoning applies to $\mathcal{A}'$.

Suppose that $\mathcal{A}(f,g)^{\top} =0$ for some $(f,g)^\top \in L^2(\partial \Omega) \times L_0^2(\partial \Omega)$, then the function defined by \eqref{deux} is the solution of problem \eqref{4Diri}.
Therefore, $u\equiv 0$ in $B_r$. In particular, $S^r f=S_A g=0$. By Lemma \ref{lemnont2}.(3), we get $f=0$. By Lemma \ref{lemnont1}.(5), $B(g,0)=(0,0)$ has a unique solution $(0,0)$, so $g=0$. Thus $\mathcal{A}$ is injective. The proof is complete.
\end{proof}

\noindent\textit{Proof of Theorem \ref{blethm}.}
We divide the proof into several steps.

\emph{Step 1}. We verify that $I-A_0$ is positive definite. For any $\xi \in \mathbb{R}^d$, we define $w_0(y): = \xi\cdot (y+ \chi)$ and $w_1(y) := \xi \cdot y$. We have
\begin{equation}
\begin{aligned}
    &\int_{Y\setminus \overline{D}} \nabla w_0 \cdot \nabla(w_0-w_1) \,dy \\
    &=\int_{\partial Y} (w_0-w_1) \frac{\partial w_0}{\partial \mathbf{n}} \Big|_-\,dS(y)  -\int_{\partial D} (w_0-w_1) \frac{\partial w_0}{\partial \mathbf{n}} \Big|_+\,dS(y)  \\
    &\quad\,- \int_{Y\setminus \overline{D}} (w_0-w_1)\Delta w_0\,dy \\
    &  =0,
\end{aligned}
\end{equation}
where the last equality follows from the definition \eqref{def0chi} of $\chi$ and $w_0-w_1 = \xi \cdot \chi$.
Therefore,
\begin{equation}
    \begin{aligned}
       |Y\setminus \overline{D}| \,|\xi|^2 &=\int_{Y\setminus \overline{D}} |\nabla w_1|^2\,dy = \int_{Y\setminus \overline{D}} |\nabla w_0|^2\,dy + \int_{Y\setminus \overline{D}} |\nabla (w_0-w_1)|^2\,dy \\
        &\geq \int_{Y\setminus \overline{D}} |\nabla w_0|^2\,dy = \int_{Y\setminus \overline{D}}  \nabla w_0 \cdot \nabla w_1\,dy  \\
        & = A_0\xi\cdot \xi ,
    \end{aligned}
\end{equation}
where $|A|$ denotes the Lebesgue measure of a Lebesgue measurable set $A\subset \mathbb{R}^d$. We conclude that $(I-A_0)\xi\cdot \xi \geq |D|\,|\xi|^2$, so $I-A_0$ is positive definite.

\emph{Step 2}. We define \begin{equation}
    \left\{ \begin{aligned}
        & \Delta F = -k^2 \widehat{u}_0 && \mathrm{in}\ B_r, \\
        & F=\widehat{u}_0 && \mathrm{on}\ \partial B_r,
    \end{aligned} \right.
\end{equation}
and
\begin{equation}
    \left\{ \begin{aligned}
        & \nabla \cdot A_0\nabla G = -k^2 \mu^k_0\, \widehat{u}_0 && \mathrm{in}\ B_r, \\
        & G=0 && \mathrm{on}\ \partial B_r.
    \end{aligned} \right.
\end{equation}
Let
\begin{equation}
    v:=\left\{ \begin{aligned}
        & \widehat{u}_0-F && \mathrm{in}\ B_r\setminus \overline{\Omega}, \\
        & \widehat{u}_0-G && \mathrm{in}\ \Omega.
    \end{aligned} \right.
\end{equation}
Then, $v$ satisfies problem \eqref{42310} with $A=A_0$, $\varphi := G|_--F|_+$ and $\psi:=\frac{\partial G}{\partial \nu_A} \big|_{-} -\frac{\partial F}{\partial \mathbf{n}} \big|_{+}$. By Theorem \ref{thm44}, we get
\begin{equation}
    \widehat{u}_0 = \mathcal{S}_{A_0} g + G \qquad \mathrm{in} \ \Omega,
\end{equation}
where $g$ is given by the boundary integral equation \eqref{sinv2} with $A=A_0$. We have
\begin{equation}
\begin{aligned}
     \mathcal{E}(G) &= \varepsilon^{1/2} \| G\|_{H^1(\Omega)} + \| \nabla G \|_{L^2(\mathcal{O}_{(d+2)\varepsilon})} + \varepsilon \| \nabla^2 G \|_{L^2(\Omega \setminus \overline{\mathcal{O}_{(d-1)\varepsilon}})}\\
     & \leq C\varepsilon^{1/2} \|G\|_{H^2(B_r)} \leq C\varepsilon^{1/2} \| \widehat{u}_0 \|_{L^2(B_r)} \\
     & \leq C \varepsilon^{1/2} \| \uinc \|_{H^1(B_r)},
\end{aligned}
\end{equation}
where the first inequality follows from Lemma \ref{blelem}, the second inequality follows from the elliptic regularity, and the last inequality follows from Theorem \ref{weakso}.

It follows from the non-tangential maximal function estimate \eqref{inenontan} that, for $t>0$,
\begin{equation}\label{nablaSg}
    \int_{\mathcal{O}_t} |\nabla \mathcal{S}_{ A_0} g(x)|^2\,dx\leq Ct  \int_{\partial \Omega} |\nabla (\mathcal{S}_{ A_0} g)^*_-(x)|^2\,dx \leq Ct\| g \|_{L^2(\partial \Omega)}^2.
\end{equation}
We now use the interior estimate
\begin{equation}
    |\nabla^2 \mathcal{S}_{ A_0}g(x)|\leq \frac{C}{\delta(x)} \left( \frac{1}{|B(x,\delta(x)/8)|} \int_{B(x,\delta(x)/8)} |\nabla \mathcal{S}_{ A_0}g(y)|^2\,dy \right)^{1/2},
\end{equation}
where $\delta(x) = \mathrm{dist}\,(x,\partial \Omega)$. By Fubini theorem, we get
\begin{equation}
    \begin{aligned}
        t\left( \int_{\Omega \setminus \mathcal{O}_t} |\nabla^{2} \mathcal{S}_{ A_0} g (x)|^{2} \, dx \right)^{1/2}
& \leq C t \left( \int_{\Omega \setminus \mathcal{O}_t}  \int_{B(x,\delta(x)/8)}
  \frac{|\nabla \mathcal{S}_{ A_0} g(y)|^{2}}{\delta(x)^{d+2}} \, dydx \right)^{1/2}
\\
&\leq C t^{1/2} \| g \|_{L^2(\partial \Omega)}.
    \end{aligned}
\end{equation}
where we have used \eqref{nablaSg} for the last inequality. By \eqref{441fg} we have
\begin{equation}\label{gesti}
\begin{aligned}
    \| g \|_{L^2(\partial \Omega)} &\leq C\left( \| G|_--F|_+ \|_{H^1(\partial \Omega)}  + \left\| \frac{\partial G}{\partial \nu_A} \Big|_{-} -\frac{\partial F}{\partial \mathbf{n}} \Big|_{+} \right\|_{L^2(\partial \Omega)} \right) \\
    &\leq C( \|G\|_{H^2(B_r)} + \| F\|_{H^2(B_r)}  ) \\
    & \leq C \| \widehat{u}_0 \|_{L^2(B_{r+1})} \\
    &\leq C\| \uinc \|_{H^1(B_{r+1})},
\end{aligned}
\end{equation}
where we used the elliptic regularity in the third inequality. Combining \eqref{nablaSg}-\eqref{gesti}, we get
\begin{equation}
    \mathcal{E}(\mathcal{S}_A g) \leq C\varepsilon^{1/2} \|\uinc\|_{H^1(B_{r+1})}.
\end{equation}
Therefore, we conclude that
\begin{equation}
    \mathcal{E}(\widehat{u}_0) \leq \mathcal{E}(\mathcal{S}_A g) + \mathcal{E}(G) \leq C\varepsilon^{1/2} \|\uinc\|_{H^1(B_{r+1})}.
\end{equation}
The proof is complete. \hfill $\square$

\begin{corollary}\label{col46L}
Assume $k\in (0,\infty)\setminus \Sigma_{D,1}$. Let $r>0$ such that $\overline{\Omega} \subset B_r$. For any $ g\in L^2(\mathbb{R}^d)$ with $\mathrm{supp}\, g\subset B_r $, let $\widehat{w}_0$ be either the solution to 
    \begin{equation}\label{coLg}
    \left\{
    \begin{aligned}
        & (L_0+k^2\mu^k_0)\widehat{w}_0=g  && \mathrm{in}\ \R^d, \\
        & \widehat{w}_0 \in \mathrm{SRC}(k),
    \end{aligned}\right.
    \end{equation}
    or the solution to
    \begin{equation}\left\{
    \begin{aligned}
        & L_0\widehat{w}_0=g  && \mathrm{in}\ B_r, \\
        & \widehat{w}_0 =0 && \mathrm{on} \ \partial B_r.
    \end{aligned}\right.
    \end{equation}
We have
    \begin{equation}
        \mathcal{E}(\widehat{w}_0) \le C \varepsilon^{\frac{1}{2}}\|g\|_{L^2(B_r)},
    \end{equation}
    where $C>0$ is a constant that depends only on $d,k,r,D$ and $\Omega$.
\end{corollary}
\begin{proof}
We only prove the conclusion for the case where $\widehat{w}_0$ is the solution to \eqref{coLg}, as the other case can be treated in a similar manner.

By an argument similar to that of Theorem \ref{weakso}, we have
    \begin{equation}
        \|\widehat{w}_0\|_{H^1(B_{r+1})} \le C \|g\|_{L^2(B_r)}.
    \end{equation}
    Next, define
    \begin{equation}
        \left\{
        \begin{aligned}
            &\Delta F = -k^2\widehat{w}_0 +g && \mathrm{in}\ B_r, \\
            &F=\widehat{w}_0 && \mathrm{on}\ \partial B_r,
        \end{aligned}
        \right.
        \quad \mathrm{and} \quad
        \left\{
        \begin{aligned}
            &\nabla\cdot A_0 \nabla G = -k^2\mu_0^k\widehat{w}_0 +g && \mathrm{in}\ B_r, \\
            &G=\widehat{w}_0 && \mathrm{on}\ \partial B_r.
        \end{aligned}
        \right.
    \end{equation}
    The remainder of the proof then follows by repeating the argument of Theorem \ref{blethm}.
\end{proof}

\section{Scattering resonances for the outgoing resolvent}\label{secresout}

In this section, we extend the outgoing resolvent $R_{\varepsilon}(k)$ for $k>0$ to $R_\eps(z)$ for $z\in \mathbb{C}$, and study the limiting behavior of $R_{\varepsilon}(z)$ as $\varepsilon \rightarrow 0$.

For simplicity, we assume that the dimension $d\geq 2$ is odd in this section. For even $d$, corresponding results can be  obtained in a similar manner, see Remark \ref{oddeven}.

\subsection{Meromorphic continuation}\label{secscrses}
We consider $L_{\varepsilon} = \nabla \cdot A_{\varepsilon} \nabla$ as an unbounded operator on $L^2(\mathbb{R}^d)$ with the domain
\begin{equation}\label{domainD}
    \mathcal{D}(L_{\varepsilon}) := \{ u \in  H^1(\mathbb{R}^d): \nabla \cdot A_{\varepsilon} \nabla u \in L^2(\mathbb{R}^d) \} ,
\end{equation}
where $\nabla \cdot A_{\varepsilon} \nabla u$ is understood in the sense of distribution. 

Similarly, for any $r>0$ such that $\overline{\Omega} \subset B_r$, we consider $L^r_{\varepsilon} = \nabla \cdot A_{\varepsilon} \nabla$ as an unbounded operator on $L^2(B_r)$ with the domain
\begin{equation}\label{domainDr}
    \mathcal{D}(L^r_{\varepsilon}):= \{ u \in  H_0^1(B_r): \nabla \cdot A_{\varepsilon} \nabla u \in L^2(B_r)\} .
\end{equation}
By the theory of quadratic forms (see \cite[Section 8.6]{reed1980methods}), $L_{\varepsilon}$ and $L_{\varepsilon}^r$ are self-adjoint.

For $\mathrm{Im}\,z >0$, let $R_{\Delta}(z) = (\Delta +z^2)^{-1}$ be the free outgoing resolvent on $L^2(\mathbb{R}^d)$, \emph{i.e.}, 
\begin{equation}\label{repRde}
    R_\Delta (z) g (x) = \int_{\mathbb{R}^d} G^z(x,y)g(y)\,dy, \qquad \forall g \in L^2(\mathbb{R}^d),
\end{equation}
where $G^z(x,y)$ is the fundamental solution to the Helmholtz operator $\Delta +z^2$, see \eqref{320Gk}.
From the expression \eqref{repRde}, it is clear that $R_{\Delta}(z)$ admits a holomorphic continuation on $\mathbb{C}$ as a family of operators from $L^2_{\mathrm{comp}}(\mathbb{R}^d) $ to $L^2_{\mathrm{loc}}(\mathbb{R}^d) $. For further details, see \cite{dyatlov2019mathematical}.

\begin{lemma}\label{lemmaresolvent}
Let $r>0$ satisfy $\overline{\Omega} \subset B_r$. Let $\chi \in C_0^{\infty} (\mathbb{R}^d; [0,1])$ such that 
\begin{equation}
    \mathrm{supp}\,\chi \subset B_r \qquad \mathrm{and} \qquad \overline{\Omega} \subset \{\chi =1\}^{\circ}.
\end{equation} 
For any bounded operator $T$ on $L^2(\mathbb{R}^d)$, we define
\begin{equation}\label{defW}
    W_{\varepsilon} := (L^r_{\varepsilon})^{-1} (I-T).
\end{equation}
Then, as bounded operators on $L^2(\mathbb{R}^d)$, the following identity holds:
\begin{equation}\label{eqmainidentity}
    \begin{aligned}
          (L_{\varepsilon} + z^2) Q_{\varepsilon}(z) = I + K_{\varepsilon}(z)  , \quad  \mathrm{Im} \, z >0,
    \end{aligned}
\end{equation}
where
\begin{equation}\label{defQ}
    Q_{\varepsilon}(z) := \chi W_{\varepsilon} +  (1-\chi)R_{\Delta}(z),
\end{equation}
and 
\begin{equation}\label{defK}
 K_{\varepsilon}(z) := \chi (z^2  W_{\varepsilon} -T) + [\Delta,\chi] (W_{\varepsilon} - R_{\Delta}(z) )  .
\end{equation}
Here $[A,B] := AB-BA$ denotes the commutator of operators $A$ and $B$.
\end{lemma}
\begin{proof}
Through direct computation, we obtain
    \begin{equation}
        \begin{aligned}
             (L_{\varepsilon}+z^2) Q_\varepsilon (z) &= (L_{\varepsilon}+z^2) \chi W_{\varepsilon} + (L_{\varepsilon}+z^2)  (1-\chi)R_{\Delta}(z) \\
             &= \chi (L_{\varepsilon} + z^2) W_{\varepsilon} + 
             (1-\chi) (L_{\varepsilon}+z^2)  R_{\Delta}(z)  + [L_{\varepsilon},\chi] (W_{\varepsilon} - R_{\Delta}(z) ).
        \end{aligned}
    \end{equation}
Since $\mathrm{supp}\,\chi \subset B_r$ and $\mathrm{supp}\,(1-\chi) \subset \mathbb{R}^d \setminus \overline{\Omega}$, we have $[L_{\varepsilon},\chi]  = [\Delta, \chi]$, $\chi L_{\varepsilon} (L^r_{\varepsilon})^{-1} = \chi$ and
\begin{equation}
    (1-\chi) (L_{\varepsilon}+z^2)  R_{\Delta}(z) = 1-\chi.
\end{equation}
Hence,
\begin{equation}
    (L_{\varepsilon}+z^2) Q_\varepsilon (z)  = I  + \chi (z^2  W_{\varepsilon} -T) + [\Delta,\chi] (W_{\varepsilon} - R_{\Delta}(z) ) = I+K_{\varepsilon}(z).
\end{equation}
The proof is complete.
\end{proof}

Identities in the form of \eqref{eqmainidentity} are fundamental for constructing the meromorphic continuation of the resolvent $(L_{\varepsilon}+z^2)^{-1}$ and for analyzing the scattering resonances of $L_{\varepsilon}$, as they provide a parametrix $Q_{\varepsilon}(z)$ of $L_{\varepsilon} +z^2$ modulo a compact remainder $K_{\varepsilon}(z)$. 

Indeed, for each fixed $\varepsilon$, the operator $L_{\varepsilon}$ falls in the black-box Hamiltonian framework \cite[Chapter 4]{dyatlov2019mathematical}, where an analogous identity \cite[(4.2.10)]{dyatlov2019mathematical} plays the same role as \eqref{eqmainidentity}. Nevertheless, that abstract construction there does not capture the homogenization limit $\varepsilon \rightarrow 0$ and, in particular, does not yield estimates uniform in $\varepsilon$.

By contrast, as shown below, the operators $Q_{\varepsilon}(z)$ and $K_{\varepsilon}(z)$ defined in \eqref{defQ} and \eqref{defK} are well behaved in the limit $\varepsilon \to 0$. Consequently, the identity \eqref{eqmainidentity} is effective for studying the meromorphic continuation of the resolvent $(L_{\varepsilon}+z^2)^{-1}$ and for analyzing the scattering resonances of $L_{\varepsilon}$ as $\varepsilon \rightarrow 0$. 

To proceed, we need to choose a suitable $T$. In the absence of high-contrast, we simply set $T=0$. While in the high-contrast case, we select $T$ as follows.

For any Lipschitz domain $U$, we denote by $\Delta_U$ the Laplacian on $U$ with the homogeneous Dirichlet boundary condition. Given $\eps>0$, recall the set $D_\eps$ of scatters defined in \eqref{def0D} consists of $\eps$-sized copies of the model set $D$ translated by $\eps \mathbf{m}$ for integers $\mathbf{m} \in J_\eps \subset \Z^d$. 

For $z\in \mathbb{C}\setminus \Sigma_D$, define
\begin{equation}\label{defTeps}
    T_{\varepsilon}(z): = z^2(\varepsilon^2 \Delta_{D_{\varepsilon}} +z^2)^{-1}  + z^2 \mathbbm{1}_{\mathbb{R}^d\setminus \overline{\Omega}} R_{\Delta}(z).
\end{equation}
Here, $(\eps^2\Delta_{D_\eps} + z^2)^{-1}$ is the solution operator of $\eps^2\Delta_{D_\eps} + z^2$ with zero Dirichlet boundary condition. In fact, let $\{(\varphi_j, \lambda_j^2)\}_{j=1}^{\infty}$ be the eigenpairs of $-\Delta_D$, with $\{\varphi_j\}_{j=1}^\infty$ forming an orthonormal basis of $L^2(D)$, then we have
\begin{equation}\label{defdeltaDepsp}
    (\varepsilon^2 \Delta_{D_{\varepsilon}} + z^2)^{-1} g(x)= \sum_{\mathbf{m}\in J_\varepsilon} \sum_{j=1}^\infty \frac{1}{z^2-\lambda_j^2} \big\langle g(\varepsilon(\cdot+\mathbf{m}) ),\varphi_j \big\rangle_{L^2(D)} \varphi_j\left( \frac{x}{\varepsilon} - \mathbf{m} \right),
\end{equation}
Hence, for each $\eps$, $(\varepsilon^2 \Delta_{D_{\varepsilon}} + z^2)^{-1} $ forms a meromorphic family of bounded operators from $L^2(D_{\varepsilon}) $ to $ H_0^1(D_{\varepsilon})$, with poles at $\{\pm \lambda_j\}_{j=1}^\infty$. By zero extension outside $D_\eps$, we can view $(\varepsilon^2 \Delta_{D_{\varepsilon}} + z^2)^{-1} $ as a meromorphic family of operators on $L^2(\mathbb{R}^d)$. Similarly, $\chi T_\eps(z)$ becomes a meromorphic family of compact operators on $L^2(\R^d)$.

In what follows, let $T$ in \eqref{defW} be the operator $T_{\varepsilon}(z)$ defined above. We check that $K_\eps(z)$ defined in \eqref{defK}
 is a compact operator on $L^2(\R^d)$.



\begin{lemma}\label{inverIK}
   $I + K_\varepsilon(z)$ is invertible on $L^2(\mathbb{R}^d)$ for $\mathrm{Im}\,z>0$. 
\end{lemma}
\begin{proof}
Since $K_{\varepsilon}(z)$ is a compact operator on $L^2(\mathbb{R}^d)$, by the Fredholm theory, it suffices to show that $I+K_{\varepsilon}(z)$ is injective. Assume that $(I+K_{\varepsilon}(z))g =0$ for some $g \in L^2(\mathbb{R}^d)$.

Since the spectrum of $L_{\varepsilon}$ lies in $(-\infty,0]$, it follows from \eqref{eqmainidentity} that $Q_{\varepsilon}(z)g =0$. Thus,
\begin{equation}
    W_{\varepsilon} (z) g =0 \quad \mathrm{in} \ \{ \chi =1 \} , \qquad \mathrm{and} \qquad R_\Delta(z)g=0 \quad\mathrm{in} \ \{ \chi = 0 \} ,
\end{equation}
which implies that $  W_\varepsilon(z) g$ and $R_\Delta(z)g$ belong to $H^2(B_r)\cap H^1_0(B_r) $.
It follows that
\begin{equation}\label{514g}
    g =  (I-T_{\varepsilon}(z))g  = L^r_{\varepsilon} W_{\varepsilon}(z)g =\Delta W_{\varepsilon} (z)g =0 \qquad \mathrm{in}\ \Omega \setminus \overline{D_{\varepsilon}} ,
\end{equation}
and
\begin{equation}
    (I-T_{\varepsilon}(z)) g =  L^r_{\varepsilon} W_{\varepsilon}(z)g = \varepsilon^2 \Delta W_{\varepsilon}(z) g =0 \qquad \mathrm{in}\  D_{\varepsilon}.
\end{equation}
This implies that
\begin{equation}\label{516g}
    \Delta g = \Delta T_{\varepsilon} (z)g = \varepsilon^{-2}z^2(I-T_{\varepsilon} (z))g = 0 \qquad \mathrm{in} \ D_{\varepsilon}.
\end{equation}
Combining \eqref{514g} and \eqref{516g}, we get $g =0$ in $\Omega$. Let $u=W_{\varepsilon} (z)g - R_{\Delta}(z)g$, we get
\begin{equation}\label{1u}
    (\Delta+z^2) u = -(\Delta+z^2)R_{\Delta}(z) g = -g =0 \qquad \mathrm{in}\ \Omega,
\end{equation}
and
\begin{equation}\label{2u}
    \begin{aligned}
        \Delta u &= \Delta ( L^r_{\varepsilon})^{-1} (I-T_{\varepsilon}(z)) g - \Delta R_{\Delta}(z) g \\
        & =(I- z^2R_{\Delta}(z))g -\Delta R_{\Delta}(z)g  \\
        &=0
    \end{aligned} \qquad \mathrm{in} \ B_r\setminus \overline{\Omega}.
\end{equation}
From these two equations \eqref{1u}-\eqref{2u} and the regularity $u\in H^2(B_r)$, we find that $u \in H_0^1(B_r)$ satisfies 
\begin{equation}
    (\Delta + z^2 \mathbbm{1}_{\Omega}) u=0 \qquad \mathrm{in}\ B_r.
\end{equation}
Since $\mathrm{Im}\,z>0$, by Proposition \ref{uH01R} below, we get $u\equiv 0$ in $B_r$. Thus,
\begin{equation}
    R_{\Delta}(z)g  = Q_{\varepsilon}(z)g -\chi u = 0 \quad \mathrm{in} \ \mathbb{R}^d.
\end{equation}
Therefore, $g\equiv0$ in $\mathbb{R}^d$, and $I+K_{\varepsilon}(z)$ is injective. This completes the proof.
\end{proof}

\begin{proposition}\label{uH01R}
    For $z\in \mathbb{C}$, if there exists a nonzero $u \in H_0^1(B_r)$ such that $(\Delta + z^2 \mathbbm{1}_{\Omega}) u=0$ in $B_r$, then $z\in \mathbb{R}$. 
\end{proposition}
\begin{proof}
    Multiplying by $ \overline{u} $ on both sides of $(\Delta + z^2 \mathbbm{1}_{\Omega}) u=0$ and taking integral over $ B_r $ gives
    \begin{equation}
        -\int_{B_r} |\nabla u|^2\,dx + \big((\mathrm{Re}\, z )^2-(\mathrm{Im}\, z)^2\big)\int_{\Omega} |u|^2\, dx=-\mathrm{i}\,( \mathrm{Re}\, z)( \mathrm{Im}\, z) \int_{\Omega} |u|^2\, dx.
    \end{equation}
    For $ \mathrm{Re}\,z=0 $, it is clear that $ u\equiv 0 $ in $B_r$. For $ \mathrm{Re}\, z\neq 0 $ and $ \mathrm{Im}\, z\neq 0 $, taking the imaginary part of the above identity gives $ u=0 $ in $\Omega$, which further implies that $ u\equiv 0 $ in $B_r$. Therefore, $z\in \mathbb{R}$.
\end{proof}

\begin{theorem}\label{meromorphic_continuation}
    $(L_{\varepsilon} + z^2)^{-1}:L^2(\mathbb{R}^d) \rightarrow \mathcal{D}(L_{\varepsilon})$ is holomorphic for $\mathrm{Im}\,z>0$, and admits a meromorphic continuation to $\mathbb{C}$, as a family of operators from $L^2_{\mathrm{comp}}(\mathbb{R}^d)$ to $\mathcal{D}_{\mathrm{loc}}(L_{\varepsilon})$, where
    \begin{equation}
        \mathcal{D}_{\mathrm{loc}}(L_{\varepsilon}) : = \big\{ u \in  H_{\mathrm{loc}}^1(\mathbb{R}^d): \rho \in C_0^{\infty}(\mathbb{R}^d), \rho|_{B_r} \equiv 1 \Rightarrow\nabla \cdot A_{\varepsilon} \nabla (\rho u) \in L^2(\mathbb{R}^d) \big\} .
    \end{equation}
This meromorphic continuation is called the outgoing resolvent of $L_{\varepsilon}$, and is denoted by $R_{\varepsilon}(z)$. Moreover, we have the identity
\begin{equation}\label{Repdec2}
    R_{\varepsilon}(z) = Q_{\varepsilon}(z) \big\{I + K_{\varepsilon}(z) \rho\big\}^{-1} \big\{ I-K_{\varepsilon}(z) (1-\rho) \big\},
\end{equation}
where $\rho \in C_0^{\infty} (\mathbb{R}^d; [0,1])$ is a cut-off such that $B_r \subset \{\rho=1\}$.
\end{theorem}
\begin{remark}
    When $z>0$, the outgoing resolvent $R_{\varepsilon}(z)$ defined as above coincides with $R_{\varepsilon}(k)$ defined in \eqref{defRp}-\eqref{defRp2}, which can be proved by following the standard procedure (see \cite[Theorem 3.37]{dyatlov2019mathematical} for the case of potential scattering) and we omit the details. Therefore, there is no ambiguity of notations. 
\end{remark}
\begin{proof}
We divide the proof into two steps.

\emph{Step 1}. In this step, we assume that $\mathrm{Im} \, z >0$. By Lemma \ref{lemmaresolvent} and Lemma \ref{inverIK}, one has
\begin{equation}\label{522L}
    (L_{\varepsilon} + z^2)^{-1} = Q_{\varepsilon}(z)  (I + K_{\varepsilon}(z) )^{-1} ,
\end{equation}
Since $I+K_{\varepsilon}(z)$ is holomorphic and Fredholm, it follows from the analytic Fredholm theory (see, for example, \cite[Theorem C.8]{dyatlov2019mathematical}) and Lemma \ref{inverIK} that $(I + K_{\varepsilon}(z) )^{-1}$ is holomorphic. Since $Q_{\varepsilon}(z)$ is holomorphic, it follows from \eqref{522L} that $(L_{\varepsilon} + z^2)^{-1} $ is holomorphic. This proves the first part of this theorem. 

We now give the second representation of $(L_{\varepsilon} + z^2)^{-1}$, which differs from \eqref{522L}, as follows. For $\rho \in C_0^{\infty} (\mathbb{R}^d; [0,1])$ such that $B_r \subset \{\rho=1\}$, we have $(1-\rho)K_{\varepsilon}(z)=0$. From this, we have
\begin{align}
    & I+K_{\varepsilon}(z) =  \big\{ I+ K_{\varepsilon}(z) (1-\rho)\big\} \big\{I + K_{\varepsilon}(z) \rho\big\}, \label{rhoin1}\\
    &\big\{ I+ K_{\varepsilon}(z)(1-\rho)\big\}^{-1} =  I- K_{\varepsilon}(z)(1-\rho). \label{rhoin2}
\end{align}
By \eqref{rhoin1}, $I + K_{\varepsilon}(z) \rho $ is invertible, and the inverse is given by
\begin{equation}\label{527-1}
    \big(I + K_{\varepsilon}(z) \rho \big)^{-1} = \big(I+K_{\varepsilon}(z)\big)^{-1}   \big\{ I+ K_{\varepsilon}(z) (1-\rho)\big\} .
\end{equation}
Substituting \eqref{rhoin1}-\eqref{rhoin2} into \eqref{522L}, we obtain 
\begin{equation}\label{extenform}
    (L_{\varepsilon} + z^2)^{-1} = Q_{\varepsilon}(z) \big\{I + K_{\varepsilon}(z) \rho\big\}^{-1} \big\{ I-K_{\varepsilon}(z) (1-\rho) \big\}.
\end{equation}

\emph{Step 2}. We now assume $z\in \mathbb{C}$. Since $R_\Delta(z):L^2_{\mathrm{comp}}(\mathbb{R}^d) \rightarrow L^2_{\mathrm{loc}}(\mathbb{R}^d) $ is holomorphic on $\mathbb{C}$, and $T_{\varepsilon}(z): L^2_{\mathrm{comp}}(\mathbb{R}^d) \rightarrow L^2_{\mathrm{loc}}(\mathbb{R}^d) $ is meromorphic on $\mathbb{C}$, we get that
\begin{equation}
    K_{\varepsilon}(z) :L^2_{\mathrm{comp}}(\mathbb{R}^d) \rightarrow L^2_{\mathrm{comp}}(\mathbb{R}^d) 
\end{equation}
is meromorphic on $\mathbb{C}$, and $Q_{\varepsilon}(z):L^2_{\mathrm{comp}}(\mathbb{R}^d) \rightarrow \mathcal{D}_{\mathrm{loc}}(L_{\varepsilon}) $ is meromorphic on $\mathbb{C}$. 

Moreover, since $I+K_{\varepsilon}(z) \rho:L^2(\mathbb{R}^d) \rightarrow L^2(\mathbb{R}^d) $ is meromorphic and Fredholm for $z\in \mathbb{C}$, by the analytic Fredholm theory, 
\begin{equation}
    \big\{I + K_{\varepsilon}(z) \rho\big\}^{-1}:L_{\mathrm{comp}}^2(\mathbb{R}^d) \rightarrow L_{\mathrm{comp}}^2(\mathbb{R}^d)
\end{equation}
is meromorphic on $\mathbb{C}$. Therefore, \eqref{extenform} gives us the desired continuation. 
\end{proof}

\begin{proposition}\label{prp55}
    Let $R>0$ be such that $\overline{\Omega}\subset B_R$. The set of poles (in $z$) of the cutoff resolvent $R_{\varepsilon}(z)\mathbbm{1}_{B_R}$, viewed as a meromorphic family of bounded operators from $L^2(\mathbb{R}^d)$ to $L_{\mathrm{loc}}^2(\mathbb{R}^d)$, is independent of the choice of $R$; these poles are called the scattering resonances of $L_{\varepsilon}$.
\end{proposition}
\begin{proof}
Let $\overline{\Omega} \subset  B_{r_1} \subset B_{r_2} $. Since
\begin{equation}
    R_{\varepsilon}(z) \mathbbm{1}_{B_{r_1}} = R_{\varepsilon}(z) \mathbbm{1}_{B_{r_2}} \circ \mathbbm{1}_{B_{r_1}} ,
\end{equation}
the poles of $R_{\varepsilon}(z) \mathbbm{1}_{B_{r_1}} $ is contained in the poles of $R_{\varepsilon}(z) \mathbbm{1}_{B_{r_2}} $. 

Conversely, for any $g \in L^2(B_{r_2})$, we define
\begin{equation}
    W=(1-\chi_1)R_{\Delta}(z) g,
\end{equation}
where $\chi_1 \in C_0^{\infty} (\mathbb{R}^d; [0,1])$ such that 
\begin{equation}
    \mathrm{supp}\,\chi_1 \subset B_{r_1} \qquad \mathrm{and} \qquad \overline{\Omega} \subset \{\chi_1 =1\}^{\circ}.
\end{equation} 
Then, $g- (\Delta+z^2)W \in L^2(B_{r_1})$, since
\begin{equation}
    g- (\Delta+z^2)W =\chi_1 g +[\Delta,\chi_1] R_{\Delta}(z) g.
\end{equation}
Moreover, we have
\begin{equation}
    \begin{aligned}
       & (L_{\varepsilon} +z^2) \Big\{ R_{\varepsilon} (z) \mathbbm{1}_{B_{r_1}} \big(g- (\Delta+z^2)W \big)  + W\Big\} \\
       &= g- (\Delta+z^2)W+(\Delta+z^2)W \\
       &=g.
    \end{aligned}
\end{equation}
Therefore, 
\begin{equation}
\begin{aligned}
    R_{\varepsilon} (z) \mathbbm{1}_{B_{r_2}}  &= R_{\varepsilon} (z) \mathbbm{1}_{B_{r_1}} \big\{I- (\Delta+z^2)(1-\chi_1)R_{\Delta}(z)  \big\} \\
    &\quad \, + (1-\chi_1)R_{\Delta}(z) .
\end{aligned}
\end{equation}
Since $R_{\Delta}(z) $ is holomorphic, we conclude that the poles of $R_{\varepsilon} (z) \mathbbm{1}_{B_{r_2}} $ is contained in the poles of $R_{\varepsilon} (z) \mathbbm{1}_{B_{r_1}} $. The proof is complete.
\end{proof}

\subsection{\texorpdfstring{Convergence rate of the outgoing resolvent}{Rate of Repsilon}}\label{subsecrateout}

To analyze the limiting behavior of $R_{\varepsilon}(z)$, we take $R=r$ in Proposition \ref{prp55} and get
\begin{equation}\label{528E}
    R_{\varepsilon}(z) \mathbbm{1}_{B_r} = Q_{\varepsilon}(z) \big(I + K_{\varepsilon}(z)  \big)^{-1} \mathbbm{1}_{B_r}.
\end{equation}
This reduces the problem to studying the limits of $(L^r_{\varepsilon})^{-1}$ and $(\varepsilon^2\Delta_{D_{\varepsilon}} + z^2)^{-1}$, which is considerably simpler, as both operators are defined on bounded domains. The corresponding results are established as follows.

Recall that $A_0$ is the homogenized coefficient defined by \eqref{homogenizedcoe}. We consider $L^r_0:= \nabla \cdot A_0 \nabla$ as an unbounded operator on $L^2(B_r)$ with the domain
\begin{equation}\label{domainD2}
    \mathcal{D}(L^r_0):= \{ u \in  H_0^1(B_r): \nabla \cdot A_0 \nabla u \in L^2(B_r)\} .
\end{equation}

\begin{theorem} \label{estimate_interior_problem}
Let $v_{\varepsilon} = (L^r_{\varepsilon})^{-1}f$ and $v_0 = (L^r_0)^{-1}f$ for any $f \in L^2(B_r)$. Then,
\begin{equation}
    \| v_\varepsilon - v_0 - (\varepsilon^2 \Delta_{D_{\varepsilon}})^{-1} f-\varepsilon \chi^{\varepsilon}  \cdot \eta_\varepsilon S_\varepsilon (  \nabla v_0 ) \|_{\mathcal{H}_{\varepsilon}(B_r)} \leq C\varepsilon^{1/2} \|f \|_{L^2(B_r)},
\end{equation}
where $C>0$ is a constant that depends only on $d,r,D$ and $\Omega$.
\end{theorem}

\begin{proof}
Let
\begin{equation}
    z_\varepsilon = v_\varepsilon - v_0 -\varepsilon \chi^{\varepsilon} \cdot \eta_\varepsilon S_\varepsilon (  \nabla v_0 ),
\end{equation}
and
\begin{equation}
    w_{\varepsilon} = v_\varepsilon - v_0 - (\varepsilon^2 \Delta_{D_{\varepsilon}})^{-1} f-\varepsilon \chi^{\varepsilon}  \cdot \eta_\varepsilon S_\varepsilon (  \nabla v_0 ).
\end{equation}
Note that $z_{\varepsilon}, w_{\varepsilon} \in H_0^1(B_r)$. Using integration by parts, we get
\begin{equation}
    -\int_{B_r} A_\varepsilon \nabla w_\varepsilon \cdot \overline{ \nabla w_\varepsilon} \, dx = \big\langle L_{\varepsilon} z_{\varepsilon}, w_\varepsilon \big\rangle_{H^{-1}(B_r),H^1(B_r)} + \int_{ D_\varepsilon } \nabla (\Delta_{D_{\varepsilon}})^{-1} f \cdot \overline{  \nabla w_\varepsilon} \,dx.
\end{equation}

Let $z=0$, $U = B_r$, $\psi = w_{\varepsilon}$, $\varphi = \mathcal{P}_{\varepsilon} \big(w_{\varepsilon}|_{B_r\setminus \overline{D_{\varepsilon}}} \big)$ in Lemma \ref{basiclemmasca}. Note that in this case $\Lambda \equiv 1$ in $Y$. Again, by applying a similar argument to the proof of Theorem \ref{L2velocity}, we have
\begin{equation}
    \begin{aligned}
        \left|\int_{B_r} A_\varepsilon \nabla w_\varepsilon \cdot \overline{ \nabla w_\varepsilon} \, dx \right| &\leq C\mathcal{E}(v_0) \| \nabla w_{\varepsilon}\|_{ L^2(\Omega \setminus \overline{D_\varepsilon} ) } 
        \\
        &\quad\, +C\left|  \int_{D_\varepsilon} f (\overline{\psi} - \overline{\varphi})\,dx + \int_{ D_\varepsilon } \nabla (\Delta_{D_{\varepsilon}})^{-1} f \cdot \overline{  \nabla \psi } \,dx\right|,
    \end{aligned}
\end{equation}
where, for the term in the second line, we have
\begin{equation}
    \int_{D_\varepsilon} f (\overline{\psi} - \overline{\varphi})\,dx + \int_{ D_\varepsilon } \nabla (\Delta_{D_{\varepsilon}})^{-1} f \cdot \overline{  \nabla \psi } \,dx= \int_{ D_\varepsilon } \nabla (\Delta_{D_{\varepsilon}})^{-1} f \cdot \overline{  \nabla \varphi } \,dx.
\end{equation}
Since $\| \nabla (\Delta_{D_{\varepsilon}})^{-1} f\|_{L^2(D_{\varepsilon})} \leq C \varepsilon \|f \|_{L^2(D_{\varepsilon})}$, we obtain
\begin{equation}
    \begin{aligned}
        \left|\int_{B_r} A_\varepsilon \nabla w_\varepsilon \cdot \overline{ \nabla w_\varepsilon} \, dx \right| \leq C \big( \mathcal{E}( v_0) + \varepsilon \|f\|_{ L^2(D_\varepsilon)} \big) \| \nabla w_{\varepsilon}\|_{ L^2(B_r \setminus \overline{D_\varepsilon } ) } .
    \end{aligned}
\end{equation}
It follows that
\begin{equation}
    \|\nabla w_\varepsilon\|_{L^2(B_r\setminus \overline{D_\varepsilon})} +\varepsilon\|\nabla w_\varepsilon\|_{L^2(D_\varepsilon)} \leq C \mathcal{E}( v_0) + C\varepsilon \|f\|_{ L^2(D_\varepsilon)}.
\end{equation}
Consequently, using Poincar\'{e} inequality, we have
\begin{equation}
    \begin{aligned}
        \|w_{\varepsilon} \|_{L^2(B_r)} &\leq \big\| \mathcal{P}_{\varepsilon} \big(w_{\varepsilon} |_{B_r\setminus \overline{D_{\varepsilon}}} \big) \big\|_{L^2(B_r)} + \big\| w_{\varepsilon} -  \mathcal{P}_{\varepsilon} \big(w_{\varepsilon} |_{B_r\setminus \overline{D_{\varepsilon}}} \big) \big\|_{L^2(D_{\varepsilon})} \\
        & \leq C\big\| \nabla\mathcal{P}_{\varepsilon} \big(w_{\varepsilon} |_{B_r\setminus \overline{D_{\varepsilon}}} \big) \big\|_{L^2(B_r)} + C\varepsilon\big\| \nabla w_{\varepsilon} -  \nabla \mathcal{P}_{\varepsilon} \big(w_{\varepsilon} |_{B_r\setminus \overline{D_{\varepsilon}}} \big) \big\|_{L^2(D_{\varepsilon})}  \\
        & \leq C \| \nabla w_{\varepsilon}  \|_{L^2(B_r\setminus \overline{D_{\varepsilon}})} + C\varepsilon \| \nabla w_{\varepsilon} \|_{L^2(D_{\varepsilon})} \\
        &\leq C \mathcal{E}( v_0) + C\varepsilon \|f\|_{ L^2(D_\varepsilon)}.
    \end{aligned}
\end{equation}
The proof is complete by applying Corollary \ref{col46L} to $v_0$.
\end{proof}

We now study the limit of $(\varepsilon^2 \Delta_{D_{\varepsilon}} + z^2)^{-1}$. Heuristically, $(\varepsilon^2\Delta_{D_{\varepsilon}} +z^2)^{-1}$ is an integral operator with a highly oscillatory Schwartz kernel. Therefore, one does not expect $(\varepsilon^2\Delta_{D_{\varepsilon}} +z^2)^{-1}$ to converge in operator norm on $L^2$ space as $\varepsilon \rightarrow 0$. However, by using the two-scale convergence method \cite{allaire1992homogenization}, it is not difficult to show that,
for any $g \in L^2(\Omega)$ and $z\in \mathbb{C}$ such that $z^2$ is not the eigenvalue of $-\Delta_D$,
\begin{equation}
    (\varepsilon^2\Delta_{D_{\varepsilon}} +z^2)^{-1} g \quad\textrm{two-scale converges to } (\Delta_D +z^2)^{-1}[1] (y) g(x).
\end{equation}
However, qualitative two-scale convergence is not sufficient for our purpose, since we wish to quantify the convergence of $(\varepsilon^2\Delta_{D_{\varepsilon}} +z^2)^{-1}$ and obtain convergence rates under some appropriate operator norms. For this, we employ the periodic unfolding method; the details are presented in Appendix \ref{secunfold}.

\begin{lemma}\label{estimateforRepsilon}
    For $z\in \mathbb{C} \setminus \Sigma_D$, let $s= \mathrm{dist}\,(z,\Sigma_D)$. We have
    \begin{equation}\label{lapDep}
        \|(\varepsilon^2 \Delta_{D_{\varepsilon}} + z^2)^{-1}\|_{L^2(D_\varepsilon)} \leq \frac{C }{ s } .  
    \end{equation}
Moreover, for a Lipschitz domain $U\subset \mathbb{R}^d$, we define the homogeneous $\dot{H}^1_0(U)$ Sobolev norm as
\begin{equation}
    \|  \phi \|_{\dot{H}^1_0(U)}:=\| \nabla \phi \|_{L^2(U)},
\end{equation}
and define the homogeneous $\dot{H}^{-1}(U)$ Sobolev norm as the dual norm of the homogeneous $\dot{H}^1_0(U)$ Sobolev norm, i.e.,
\begin{equation}
    \| g \|_{\dot{H}^{-1}(U)} : = \sup_{\|\phi \|_{\dot{H}^1_0(U)} =1 } \langle g,\phi \rangle_{H^{-1}(U),H_0^1(U)}. 
\end{equation}
Then, we have
\begin{equation}\label{normH-1}
    \|  (\varepsilon^2 \Delta_{D_{\varepsilon}} + z^2)^{-1} \|_{\dot{H}^{-1}(D_\varepsilon)} \leq \frac{C }{ s } .
\end{equation}
Here, $C>0$ is a constant that depends only on $d$ and $D$.
\end{lemma}
\begin{proof}
    By the definition \eqref{defdeltaDepsp} of $(\varepsilon^2 \Delta_{D_{\varepsilon}} + z^2)^{-1}$, for any $g \in L^2(D_{\varepsilon})$, we have
\begin{equation}
    \begin{aligned}
        \|(\varepsilon^2 \Delta_{D_{\varepsilon}} + z^2)^{-1} g\|^2_{L^2(D_\varepsilon)}
    &=  \varepsilon^d  \sum_{\mathbf{m}\in J_\varepsilon} \sum_{j=1}^\infty \frac{1}{|z^2-\lambda_j^2|^2} \Big| \big\langle g(\varepsilon(\cdot+\mathbf{m}) ),\varphi_j \big\rangle_{L^2(D)} \Big|^2\\
    & \leq \frac{C\varepsilon^d}{ s^2}  \sum_{\mathbf{m}\in J_\varepsilon}  \|g\left(\varepsilon(\cdot+\mathbf{m})\right)\|_{L^2(D)}^2 \\
    & = \frac{C }{ s^2 } \|g\|^2_{L^2(D_{\varepsilon})},
    \end{aligned}
\end{equation}
which proves \eqref{lapDep}. To prove \eqref{normH-1}, by the duality, it suffices to prove that
\begin{equation}
    \|  (\varepsilon^2 \Delta_{D_{\varepsilon}} + z^2)^{-1} g\|_{\dot{H}^1_0(D_\varepsilon)} \leq \frac{C }{ s } \|g\|_{\dot{H}^1_0(D_\varepsilon)} .
\end{equation}
Since $\langle \nabla \varphi_i , \nabla \varphi_j \rangle_{L^2(D)} = \lambda_j^2 \delta_{ij}$, we have
\begin{equation}
    \begin{aligned}
        \| \nabla (\varepsilon^2 \Delta_{D_{\varepsilon}} + z^2)^{-1} g\|^2_{L^2(D_\varepsilon)} &=\varepsilon^{d-2}  \sum_{\mathbf{m}\in J_\varepsilon} \sum_{j=1}^\infty \frac{\lambda_j^2}{|z^2-\lambda_j^2|^2} \Big| \big\langle g(\varepsilon(\cdot+\mathbf{m}) ),\varphi_j \big\rangle_{L^2(D)} \Big|^2\\
    & \leq \frac{C\varepsilon^d}{ s^2}  \sum_{\mathbf{m}\in J_\varepsilon}  \|\nabla g\left(\varepsilon(\cdot+\mathbf{m})\right)\|_{L^2(D)}^2 \\
    & = \frac{C }{ s^2 } \|\nabla g\|^2_{L^2(D_{\varepsilon})}. 
    \end{aligned}
\end{equation}
This yields \eqref{normH-1}. The proof is complete.
\end{proof}
 
Recall that $\beta(z)$ is the Zhikov's function defined in \eqref{defJikovfunc}.

\begin{proposition}\label{thm49}
   For $z\in \mathbb{C} \setminus \Sigma_D$, let  $s= \mathrm{dist}\,(z,\Sigma_D)$. Then we have
\begin{equation}
\| (\varepsilon^2 \Delta_{D_{\varepsilon}} +z^2)^{-1}  - \beta(z) \|_{L^2(\Omega) \rightarrow H^{-1}(\Omega)} \leq \frac{C }{ s } \varepsilon^{1/2},
\end{equation}
where $C>0$ is a constant that depends only on $d, D$ and $\Omega$.
\end{proposition}

\begin{proof}
   By the definition of $\mathcal{U}_{\varepsilon}$, it is clear that
\begin{equation}
        (\varepsilon^2 \Delta_{D_{\varepsilon}} +z^2)^{-1}  \mathcal{U}_{\varepsilon} = \mathcal{U}_{\varepsilon} (\Delta_{D,y}+z^2)^{-1}  .
\end{equation}
Therefore, 
    \begin{equation}
    \begin{aligned}
          &(\varepsilon^2 \Delta_{D_{\varepsilon}} +z^2)^{-1}   - \beta(z)  
            \\
            &=(\varepsilon^2 \Delta_{D_{\varepsilon}} +z^2)^{-1}  (\langle \cdot \rangle_Y - \mathcal{U}_{\varepsilon})  + (\mathcal{U}_{\varepsilon} - \langle \cdot \rangle_Y)(\Delta_{D,y}+z^2)^{-1}[1](y) .
    \end{aligned}
    \end{equation}
Since $H^{-1}(\Omega)$ norm and $\dot{H}^{-1}(\Omega)$ norm are comparable, we have
\begin{equation}
    \begin{aligned}
        &\|(\varepsilon^2 \Delta_{D_{\varepsilon}} +z^2)^{-1}  (\langle \cdot \rangle_Y - \mathcal{U}_{\varepsilon}) \|_{L^2(\Omega) \rightarrow H^{-1}(\Omega) } \\
        &\leq C\|(\varepsilon^2 \Delta_{D_{\varepsilon}} +z^2)^{-1}  \circ\mathbbm{1}_{D_{\varepsilon} }  \circ(\langle \cdot \rangle_Y - \mathcal{U}_{\varepsilon}) \|_{L^2(\Omega) \rightarrow \dot{H}^{-1}(D_{\varepsilon}) } \\
        &\leq  C\|(\varepsilon^2 \Delta_{D_{\varepsilon}} +z^2)^{-1} \|_{\dot{H}^{-1}(D_{\varepsilon})} \| \mathbbm{1}_{D_{\varepsilon} }\|_{\dot{H}^{-1}(\Omega) \rightarrow \dot{H}^{-1}(D_{\varepsilon})}\| \langle \cdot \rangle_Y - \mathcal{U}_{\varepsilon} \|_{L^2(\Omega) \rightarrow H^{-1}(\Omega)  } \\
        &\leq \frac{C}{s} \varepsilon^{1/2},
    \end{aligned}
\end{equation}
where we used \eqref{normH-1} and Lemma \ref{fold} in the last inequality. Again, by Lemma \ref{fold}, we have
\begin{equation}
\begin{aligned}
   & \| (\mathcal{U}_{\varepsilon} - \langle \cdot \rangle_Y)(\Delta_{D,y}+z^2)^{-1}[1](y) \|_{L^2(\Omega) \rightarrow H^{-1}(\Omega)}  \\
    & \leq \|\mathcal{U}_{\varepsilon} - \langle \cdot \rangle_Y \|_{L^2(\Omega) \rightarrow H^{-1}(\Omega)}\| (\Delta_{D,y}+z^2)^{-1}[1]\|_{L^2( D)} \\
    &\leq \frac{C}{ s } \varepsilon^{1/2} .
\end{aligned}
\end{equation}
The desired conclusion follows.
\end{proof}

We are in the position of studying the limit of $R_{\varepsilon}(z)$. For any $z\in \mathbb{C}$, we define
\begin{align}
    &T_0(z):= z^2 \beta(z) \mathbbm{1}_{\Omega} + z^2 \mathbbm{1}_{ \R^d \setminus \overline \Omega}  R_\Delta (z), \label{DefT0}\\
    &W_0(z) := (L_0^r)^{-1} \big( I - T_0(z) \big) ,\\
    &Q_0(z): = \chi W_0(z) + (1-\chi) R_\Delta(z) \label{DefQ0z} , \\
    & K_0(z):=\chi z^2 \big(W_{0}(z) -\mathbbm{1}_{\R^d\setminus \overline{\Omega }} R_\Delta(z) \big) + [\Delta,\chi] \big(W_{0}(z) - R_{\Delta}(z) \big). \label{DefK0}
\end{align}

For simplicity of notation, below we denote the operator norm $\|\cdot \|_{L^2(B_r) \rightarrow L^2(B_R)} $ by $\|\cdot\|_{r,R}$. When $r=R$, we abbreviate $\|\cdot\|_r := \|\cdot \|_{r,r}$.

\begin{theorem}\label{thmRBR}
For $z\in \mathbb{C}\setminus \Sigma_D$, let $s= \mathrm{dist}\,(z,\Sigma_D)$. Then for any $R>0$ such that $\overline{\Omega} \subset B_R$, we have
\begin{equation}\label{con571}
\begin{aligned}
    &\Big\| R_{\varepsilon}(z) - R^{(1)}_{\varepsilon}(z) \Big\|_{r,R} \\
    &\leq C\varepsilon^{1/2}  \left(1+  \frac{ 1 }{ s }  \right)^2 \frac{\big\| \big( I+K_0(z) \big)^{-1}  \big\|^2_r}{1- \varepsilon^{1/2}  (1+  s^{-1} ) \big\| \big( I+K_0(z) \big)^{-1}  \big\|_r } ,
\end{aligned}
\end{equation}
where $C>0$ is a constant that depends only on $d,r,R,|z|,D$, and 
\begin{equation}\label{defR1ep}
\begin{aligned}
    R^{(1)}_{\varepsilon}(z)&:=\big( Q_0(z) + (\varepsilon^2\Delta_{D_{\varepsilon}} +z^2)^{-1} \big)\big( I+K_0(z) \big)^{-1}  \\
    &\ = (L_0+z^2\mu_0^z)^{-1} \mu_0^z  + (\varepsilon^2\Delta_{D_{\varepsilon}} +z^2)^{-1} (I+K_0(z) )^{-1}.
\end{aligned}
\end{equation}
\end{theorem}
\begin{proof}
We divide the proof into three steps.

\emph{Step 1}. We first verify the last equality in \eqref{defR1ep}. We have
\begin{equation}
\begin{aligned}
    &(L_0+z^2\mu^z_0) Q_0(z) \\
    & =[\Delta,\chi]\big( W_0(z)-R_\Delta(z) \big) + (1-\chi) + \chi(I-T_0(z)) +z^2 \mu_0^z \chi W_0(z)\\
    & = I+K_0(z) -z^2\beta(z)\mathbbm{1}_{\Omega} + z^2(\mu_0^z-1)\chi W_0(z) ,
\end{aligned}
\end{equation}
it follows from $ \mu_0^z=1-z^2\beta(z)\mathbbm{1}_\Omega $ that
\begin{equation}
\begin{aligned}
    (L_0+z^2\mu^z_0) Q_0(z) & =(I+K_0(z)) + (\mu^z_0-1 ) \mathbbm{1}_\Omega ( I+z^2 W_0(z) ) \\
    & = \mu_0^z (I+K_0(z))  .
\end{aligned}
\end{equation}
This gives the last equality in \eqref{defR1ep}.

\emph{Step 2}. In this step, we show that
\begin{align}
    &\label{Qestimate0}
    \|Q_{\varepsilon}(z) - Q_0(z) -  (\varepsilon^2 \Delta_{D_{\varepsilon}} +z^2)^{-1} \|_{r,R} \leq C\left(1+ \frac{1}{ s}  \right)\varepsilon^{1/2}, \\
    &\label{Kesti}
    \|K_{\varepsilon}(z) - K_0(z) \|_{r,R}  \leq C\left(1+  \frac{ 1 }{s }  \right)\varepsilon^{1/2}.
\end{align}

 Since $T_{\varepsilon}(z)  = z^2(\varepsilon^2 \Delta_{D_{\varepsilon}} +z^2)^{-1} + z^2 \mathbbm{1}_{\mathbb{R}^d\setminus \overline{\Omega}} R_{\Delta}(z)$, and
 \begin{equation}
     T_{\varepsilon}(z) - T_0(z) = z^2 \big((\varepsilon^2 \Delta_{D_{\varepsilon}} +z^2)^{-1} - \beta(z) \mathbbm{1}_{\Omega} \big),
 \end{equation}
by Lemma \ref{estimateforRepsilon} and Theorem \ref{thm49}, we get 
\begin{equation}
    \label{normTep}
    \| T_{\varepsilon} (z)\|_{L^2(B_r)} \leq C\left(1+\frac{1  }{ s }\right), 
\end{equation}
and
\begin{equation}
\label{TepsT0}
    \| T_{\varepsilon}(z) - T_0(z) \|_{L^2(\mathbb{R}^d) \rightarrow H^{-1}(\mathbb{R}^d)} \leq \frac{C  }{ s } \varepsilon^{1/2}. 
\end{equation}
We can easily verify 
\begin{equation}
    (\varepsilon^2 \Delta_{D_{\varepsilon}})^{-1}  ( I - T_{\varepsilon}(z) )= (\varepsilon^2 \Delta_{D_{\varepsilon}} +z^2)^{-1},
\end{equation}
thus,
\begin{equation}
    \begin{aligned}
        &W_\varepsilon(z)-W_0(z)-  (\varepsilon^2 \Delta_{D_{\varepsilon}} +z^2)^{-1}
    \\
    &=  \big\{ (L^r_\varepsilon)^{-1}- ( L^r_0 )^{-1} - (\varepsilon^2 \Delta_{D_{\varepsilon}})^{-1}  \big\}  \big( I - T_\varepsilon (z) \big) + ( L^r_0 )^{-1}  \big( T_0(z)-T_{\varepsilon}(z) \big).
    \end{aligned}
\end{equation}
By \eqref{normTep}-\eqref{TepsT0}, Theorem \ref{estimate_interior_problem} and the elliptic regularity of $(L^r_0)^{-1}$, we obtain
\begin{equation}
   \label{Westima001}
    \| W_\varepsilon(z)-W_0(z)-  (\varepsilon^2 \Delta_{D_{\varepsilon}} +z^2)^{-1} \|_{L^2(B_r)\rightarrow L^2(\mathbb{R}^d) } \leq C\left(1+ \frac{ 1 }{ s }  \right) \varepsilon^{1/2}.
\end{equation}
The desired estimates \eqref{Qestimate0}-\eqref{Kesti} follow from \eqref{Westima001}, and
\begin{equation}
    Q_{\varepsilon}(z)-Q_0(z) -  (\varepsilon^2 \Delta_{D_{\varepsilon}} +z^2)^{-1}= \chi \big( W_\varepsilon(z)-W_0(z)-  (\varepsilon^2 \Delta_{D_{\varepsilon}} +z^2)^{-1} \big),
\end{equation}
and
\begin{equation}
       K_{\varepsilon}(z)  - K_0(z)  = \big(\chi z^2 + [\Delta,\chi] \big)  \big(W_{\varepsilon} (z)- W_0(z) - (\varepsilon^2 \Delta_{D_{\varepsilon}} +z^2)^{-1} \big). 
\end{equation}

\emph{Step 3}. By Step 2, we have
\begin{equation}\label{Q1rR}
    \begin{aligned}
    \|Q_{\varepsilon}(z) \|_{r,R} &\leq \|Q_0(z) \|_{r,R} + \|  (\varepsilon^2 \Delta_{D_{\varepsilon}} +z^2)^{-1} \|_r + C\left(1+ \frac{1 }{ s }  \right)\varepsilon^{1/2} \\
    &\leq C\left(1+ \frac{1 }{ s } \right)  ,
    \end{aligned}
\end{equation}
where we have used Lemma \ref{estimateforRepsilon} in the second inequality. By \eqref{528E}, we have
\begin{equation}\label{RrR}
\begin{aligned}
   & R_{\varepsilon}(z) \mathbbm{1}_{B_r} - R^{(1)}_{\varepsilon}(z) \mathbbm{1}_{B_r} \\
   & = \big\{ Q_{\varepsilon}(z) - Q_0(z) -  (\varepsilon^2 \Delta_{D_{\varepsilon}} +z^2)^{-1}  \big\}\big( I+K_0(z) \big)^{-1} \mathbbm{1}_{B_r} \\
   & \quad\, + Q_{\varepsilon}(z) \big\{ \big( I+K_\varepsilon(z) \big)^{-1}  - \big( I+K_0(z) \big)^{-1}  \big\} \mathbbm{1}_{B_r}.
\end{aligned}
\end{equation}
We also have
\begin{equation}\label{KrR}
\begin{aligned}
    &\big\| \big( I+K_\varepsilon(z) \big)^{-1}  - \big( I+K_0(z) \big)^{-1}  \big\|_r  \\
    & \leq \frac{ \big\| \big( I+K_0(z) \big)^{-1}  \big\|_r^2 }{ 1- \big\| \big( I+K_0(z) \big)^{-1}  \big\|_r \cdot \|K_{\varepsilon}(z) - K_0(z) \|_r  } \|K_{\varepsilon}(z) - K_0(z) \|_r .
\end{aligned}
\end{equation}
The desired conclusion \eqref{con571} follows by combining \eqref{Q1rR}-\eqref{KrR} and Step 2.
\end{proof}

For technical reasons, we give a slightly different decomposition of $R_{\varepsilon}(z)$. Define the Dirichlet Laplacian on the exterior domain $\mathbb{R}^d\setminus \overline{\Omega}$ by 
\begin{equation}
    \Delta_{\R^d\setminus \overline{\Omega} } : \mathcal{D}( \Delta_{\R^d\setminus \overline{\Omega} } )\subset L^2( \R^d\setminus \overline{\Omega} ) \rightarrow L^2( \R^d \setminus \overline{\Omega} ),
\end{equation}
where
\begin{equation}
    \mathcal{D}( \Delta_{\R^d\setminus \overline{\Omega}} ) := \big\{ u\in H^1_0( \R^d\setminus \overline{\Omega} ): \Delta u \in L^2( \R^d\setminus \overline{\Omega})  \big\}.
\end{equation}
For $\mathrm{Im}\,z >0$, let 
\begin{equation}
    \widetilde{R}_{\Delta}(z) = \big(\Delta_{\R^d\setminus \overline{\Omega} } +z^2 \big)^{-1}
\end{equation}
denote the outgoing resolvent related to the problem of scattering by the obstacle $\Omega$. It is well-known that $\widetilde{R}_{\Delta}(z) $ admits a meromorphic continuation to the entire complex plane as a family of operators from $L^2_{\mathrm{comp}}(\mathbb{R}^d) $ to $\mathcal{D}_{\mathrm{loc}}( \Delta_{\R^d\setminus \overline{\Omega}} )$, where
  \begin{equation}
        \mathcal{D}_{\mathrm{loc}}( \Delta_{\R^d\setminus \overline{\Omega}} ) : = \{ u \in  H_{\mathrm{loc}}^1(\mathbb{R}^d\setminus \overline{\Omega} ): u|_{\partial \Omega} =0,\,\Delta u \in L_{\mathrm{loc}}^2( \R^d\setminus \overline{\Omega})  \} .
    \end{equation}
Moreover, all poles of $\widetilde{R}_{\Delta}(z)$ lie in the lower half-plane $\{z\in \mathbb{C}:\mathrm{Im}\,z<0\}$. See, for instance, \cite[Chapter 9]{taylor2010partial}. We denote by $\Sigma_{\mathbb{R}^d\setminus \overline{\Omega}}$ the set of poles of $\widetilde{R}_{\Delta}(z)$.

We identify $\widetilde{R}_{\Delta}(z)g$ with its extension by zero to $\Omega$. Following the same argument as Lemma \ref{lemmaresolvent}, we define
\begin{align}
    & \widetilde{T}_\varepsilon(z) := z^2(\varepsilon^2 \Delta_{D_{\varepsilon}} +z^2)^{-1}  + z^2 \mathbbm{1}_{\mathbb{R}^d\setminus \overline{\Omega}} \widetilde{R}_{\Delta}(z) , \\
    &\widetilde{W}_\varepsilon(z) : =  (L_\varepsilon^r)^{-1}\big(I-\widetilde{ T}_\varepsilon(z)\big) ,\\
    &\widetilde{Q}_\varepsilon(z) : =  \chi \widetilde{W}_\varepsilon(z) +( 1-\chi ) \widetilde{R}_\Delta(z) ,\\
    &\widetilde K_\varepsilon(z): =  \chi \big( z^2\widetilde{W}_\varepsilon-\widetilde{T}_\varepsilon(z) \big) + [\Delta ,\chi ] \big(\widetilde{W}_\varepsilon (z)- \widetilde R_\Delta(z) \big) .
\end{align}
The following is an analogue of Lemma \ref{lemmaresolvent}, and its proof follows the same lines:
\begin{equation}\label{596til}
    (L_\varepsilon+z^2)\widetilde Q_\varepsilon(z) = I + \widetilde K_\varepsilon(z),\qquad \mathrm{Im}\,z>0.
\end{equation}
We also have the analogue of Lemma \ref{inverIK}:
\begin{lemma}
    $I + \widetilde{K}_\varepsilon(z)$ is invertible on $L^2(\mathbb{R}^d)$ for $\mathrm{Im}\,z>0$. 
\end{lemma}
\begin{proof}
It is clear that $\widetilde{K}_{\varepsilon}(z)$ is a compact operator on $L^2(\mathbb{R}^d)$. By the Fredholm theory, it suffices to show that $I+\widetilde{K}_{\varepsilon}(z)$ is injective. Assume that $(I+\widetilde{K}_{\varepsilon}(z))g =0$ for some $g \in L^2(\mathbb{R}^d)$.
     
Since the spectrum of $L_{\varepsilon}$ lies in $(-\infty,0]$, it follows from \eqref{596til} that $\widetilde{Q}_{\varepsilon}(z)g =0$. Thus, 
 \begin{equation}
    \widetilde{W}_{\varepsilon}(z) g =0 \quad \mathrm{in} \ \{ \chi =1 \} , \qquad \mathrm{and} \qquad \widetilde{R}_\Delta(z)g=0 \quad\mathrm{in} \ \{ \chi = 0 \} ,
\end{equation}
    which implies $ \widetilde W_\varepsilon (z)g,\,\widetilde R_\Delta(z)g   \in H^1_0( B_r\setminus \overline{\Omega} ) $. 
By the same argument in Lemma \ref{inverIK}, $g \equiv0$ in $\Omega$. Moreover,
    \begin{equation}
        \Delta \big(\widetilde W_\varepsilon(z) g -\widetilde R_\Delta(z)g\big)=\big(I-z^2 \widetilde R_\Delta(z) \big)g -\Delta \widetilde R_\Delta(z)g=0  \quad \mathrm{in}\  B_r\setminus\overline{ \Omega }.
    \end{equation}
    This implies that $ \widetilde W_\varepsilon(z) g -\widetilde R_\Delta(z)g = 0 $ in $B_r\setminus\overline{ \Omega }$. Thus,
    \begin{equation}
        \widetilde R_\Delta(z)g =\widetilde Q_\varepsilon(z)g -\chi\big( \widetilde W_\varepsilon(z) g -\widetilde R_\Delta(z)g \big) =0 \qquad \mathrm{in}\ \mathbb{R}^d\setminus \overline{\Omega}.
    \end{equation}
    It follows that $g\equiv 0$ in $\mathbb{R}^d\setminus \overline{\Omega}$. Combining this with $g \equiv0$ in $\Omega$, we conclude that $I+\widetilde{K}_{\varepsilon}(z)$ is injective. The proof is complete.
\end{proof}

Following the same argument of Theorem \ref{meromorphic_continuation}, we obtain the following decomposition of $R_{\varepsilon}(z)$, which is different from \eqref{Repdec2}:
\begin{equation}
    R_{\varepsilon}(z) = \widetilde{Q}_{\varepsilon}(z) \big\{I + \widetilde{K}_{\varepsilon}(z) \rho\big\}^{-1} \big\{ I-\widetilde{K}_{\varepsilon}(z) (1-\rho) \big\},
\end{equation}
where $\rho \in C_0^{\infty} (\mathbb{R}^d; [0,1])$ is a cut-off such that $B_r \subset \{\rho=1\}$. 

Following \eqref{DefT0}-\eqref{DefK0}, we can also define
\begin{align}
& \widetilde{T}_0(z):= z^2 \beta(z) \mathbbm{1}_{\Omega} + z^2 \mathbbm{1}_{ \R^d \setminus \overline \Omega}  \widetilde{R}_\Delta (z), \\
   &\widetilde{W}_0(z) :=   (L^r_0)^{-1} \big(I-\widetilde{ T}_0(z) \big),\\
    &\widetilde{Q}_0(z): =  \chi  \widetilde{ W}_0(z)  +( 1-\chi ) \widetilde{R}_\Delta(z),\\
    & \widetilde K_0(z):=  \chi z^2\big( \widetilde W_0(z) - \mathbbm{1}_{\R^d\setminus \overline{\Omega }} \widetilde{R}_\Delta(z)  \big) + [\Delta ,\chi ] \big( \widetilde W_0(z)- \widetilde R_\Delta (z)  \big).
\end{align}
By the same arguments, we obtain the following analogues of Theorem \ref{thmRBR}.

\begin{theorem}\label{thm511Rtilde}
For $z\in \mathbb{C}\setminus (\Sigma_D\cup \Sigma_{\mathbb{R}^d\setminus \overline{\Omega}})$, let $s= \mathrm{dist}\,(z,\Sigma_D\cup \Sigma_{\mathbb{R}^d\setminus \overline{\Omega}})$. Then for any $R>0$ such that $\overline{\Omega} \subset B_R$, we have
\begin{equation}\label{con57111}
\begin{aligned}
    &\Big\| R_{\varepsilon}(z) - R^{(2)}_{\varepsilon}(z) \Big\|_{r,R} \\
    &\leq C\varepsilon^{1/2}  \left(1+  \frac{ 1 }{ s }  \right)^2 \frac{\big\| \big( I+\widetilde{K}_0(z) \big)^{-1}  \big\|^2_r}{1- \varepsilon^{1/2}  (1+  s^{-1} ) \big\| \big( I+\widetilde{K}_0(z) \big)^{-1}  \big\|_r } ,
\end{aligned}
\end{equation}
where $C>0$ is a constant that depends only on $d,r,R,|z|,D$, and 
\begin{equation}\label{defR1ep11}
\begin{aligned}
    R^{(2)}_{\varepsilon}(z)&:=\big( \widetilde{Q}_0(z) + (\varepsilon^2\Delta_{D_{\varepsilon}} +z^2)^{-1} \big)\big( I+\widetilde{K}_0(z) \big)^{-1}  \\
    &\ =  (L_0+z^2\mu_0^z)^{-1}\mu_0^z + (\varepsilon^2\Delta_{D_{\varepsilon}} +z^2)^{-1} (I+\widetilde{K}_0(z) )^{-1}.
\end{aligned}
\end{equation}
\end{theorem}

\subsection{\texorpdfstring{Uniform resolvent estimates and proofs of Theorem \ref{mainresult1} and \ref{thm5121main}}{Uniform L2 resolvent estimate}}

In this section, we prove Theorem \ref{mainresult1} and \ref{thm5121main}. To this end, we need to introduce the outgoing resolvent and the scattering resonances associated with the homogenized scattering problem \eqref{equ0}. This requires additional effort, as the effective refraction $\mu^z$ is a nonlinear function of $z$. We present the details in the Appendix \ref{appendixA}.

We note that the set of the scattering resonances associated with the homogenized scattering problem \eqref{equ0} are denoted by $\Sigma_{\mathrm{hom}}$, see Definition \ref{DefA2z}.

\begin{lemma}
\label{pole_K01}
We have
\begin{equation}
    \big\{ \mathrm{ Poles\ of \ } (I+K_0(z) )^{-1} \mathrm{ \ in\ }\mathbb{C}^{\circ}\big\} \subset \Sigma_{\mathrm{hom}} \cup \Sigma_r.
\end{equation}
Moreover, $\Sigma_r \subset \mathbb{R}$.
\end{lemma}

\begin{proof}
We divide the proof into two steps. 

\emph{Step 1}. We first show that $0$ is not a pole of $(I+K_0(z))^{-1}$. Since $K_0(z)$ is compact, it remains to prove that $ I+K_0(0) $ is injective. Assume $g\in L^2_{\mathrm{comp}}(\R^d) $ satisfy $(I+K_0(0))g=0$. By \eqref{defR1ep}, we have $ L_0 Q_0(0)g=0 $ in $\mathbb{R}^d$. Since $Q_0(0)g  = R_\Delta(0)g  $ in $\R^d \setminus B_r$, we get, as $x\rightarrow \infty$,
\begin{equation}\label{asymp_Q00}
    Q_0(0)g(x) =\left\{
    \begin{aligned}
    &O(\log |x|), && d=2,\\
    &    O(|x|^{2-d}), && d\ge 3,  
    \end{aligned}
    \right. \quad \mathrm{and} \quad \nabla Q_0(0)g(x) = O( |x|^{1-d} )  .
\end{equation}
Therefore, $ \nabla Q_0(0)g\in L^2(\R^d) $. Multiplying by $  \overline{Q_0(0)g} $ on both sides of the equation $ L_0 Q_0(0)g=0 $ in $\mathbb{R}^d$ and integrating over $ B_R $ for any $ R>r $, we obtain
\begin{equation}
    0 = \int_{B_R} A_0 \nabla Q_0(0)g \cdot \nabla \overline{Q_0(0)g} \;dx - \int_{\partial B_R } \overline{Q_0(0)g}(x)\frac{\partial Q_0(0)g}{\partial \mathbf{n}} \Big|_-(x)\, dS(x).
\end{equation}
Let $ R\rightarrow \infty $, using \eqref{asymp_Q00}, the second term on the right-hand side vanishes. By the ellipticity of $ A_0 $, we get $ \nabla Q_0(0)g\equiv0 $ in $\mathbb{R}^d$. Since $ Q_0(0)g $ decays at infinity, we conclude that $ Q_0(0)g\equiv 0 $ in $\mathbb{R}^d$.

By the definition of $Q_0(0)$, we have
\begin{equation}\label{pole_K01_eq11}
     (L_0^r )^{-1} g =0 \quad \mathrm{in} \ \{\chi = 1\}\supset\supset\overline{\Omega} ,
\end{equation}
and
\begin{equation}
    R_\Delta(0)g = 0 \qquad \mathrm{in} \ \{\chi = 0\} \supset \mathbb{R}^d\setminus \overline{B_r} .
\end{equation}
This implies that $ u:=R_\Delta(0)g-(L_0^r)^{-1} g \in H^1_0(B_r)\cap H^2(B_r) $. It is clear that $\Delta u =0$ in $B_{r} \setminus \partial \Omega$. Hence, $u \in H^1_0(B_r)\cap H^2(B_r)$ allows us to obtain $\Delta u=0$ in $B_r$. Consequently, $u \equiv 0$ in $\mathbb{R}^d$. It then follows from
\begin{equation}
    R_\Delta(0)g = Q_0(0)g + \chi u=0 \qquad \mathrm{in}\ \mathbb{R}^d
\end{equation}
that $g\equiv 0$ in $\mathbb{R}^d$. Therefore, $ I+K_0(0) $ is injective, and $0$ is not a pole of $(I+K_0(z))^{-1}$.

\emph{Step 2}. We now assume that $ z\in \mathbb{C}^{\circ}\setminus \{0\}$ is a pole of $( I+K_0  (z) )^{-1}$. Since $K_0(z)$ is compact, there exists a nonzero $g \in L^2(B_r)$ such that $(I+K_0(z)) g =0$ in $B_r$. By \eqref{defR1ep}, we have
\begin{equation}
(L_0+z^2\mu^z_0) Q_0(z)g =\mu_0^z(I+K_0(z)) g =0.
\end{equation}
If $Q_0(z)g \neq 0$, by Theorem \ref{thmA3u}, $Q_0(z)g$ is a resonant state. It follows that $z \in \Sigma_{\mathrm{hom}}$. We now assume $Q_0(z)g=0$. Then, we have
\begin{equation}\label{pole_K01_eq1}
    W_0(z)g = (L_0^r )^{-1} ( I-T_0 (z)) g =0 \quad \mathrm{in} \ \{\chi = 1\}\supset\supset\overline{\Omega} ,
\end{equation}
and
\begin{equation}
    R_\Delta(z)g = 0 \qquad \mathrm{in} \ \{\chi = 0\} \supset \mathbb{R}^d\setminus \overline{B_r} .
\end{equation}
This implies that $ W_0(z) g,\,R_\Delta(z)g  \in H^1_0(B_r)\cap H^2(B_r) $, and $ g =\big(I+K_0(z) \big)g - z^2W_0(z)g = 0$ in $ \Omega$.

Let $u=R_\Delta(z)g - W_0 (z)g$, by the definition of $ W_0 (z)$, we get
     \begin{equation}\label{Deltau1}
         \Delta u = \Delta R_\Delta (z) g -( -z^2 R_\Delta (z)g +g ) =0 \quad \mathrm{in}\ B_{r} \setminus \overline{\Omega } .
     \end{equation}
The above identity can be verified directly for $\mathrm{Im}\,z>0$, and extended to the lower half-plane by holomorphic continuation. Similarly, by \eqref{pole_K01_eq1}, we have
\begin{equation}\label{Deltau2}
    (\Delta +z^2) u = g - 0 =0 \quad \mathrm{in}\ \Omega.
\end{equation}
Hence, $u \in H^1_0(B_r)\cap H^2(B_r)$ allows us to combine \eqref{Deltau1}-\eqref{Deltau2} to obtain
\begin{equation}
    (\Delta +z^2 \mathbbm{1}_{\Omega} ) u =0 \qquad \mathrm{in}\ B_r.
\end{equation}
To conclude that $z\in \Sigma_r$, it remains to verify that $u$ is nonzero. Otherwise, by the definition of $Q_0(z)$, we have
\begin{equation}
    R_\Delta(z)g = Q_0(z)g + \chi u=0 \qquad \mathrm{in}\ \mathbb{R}^d.
\end{equation}
This implies that $g\equiv 0$ in $\mathbb{R}^d$, which contradicts the assumption that $g$ is nonzero. $\Sigma_r \subset \mathbb{R}$ is guaranteed by Proposition \ref{uH01R}. The proof is complete.
\end{proof}

\begin{lemma}
\label{pole_K0116}
We have
\begin{equation}
    \big\{ \mathrm{Poles\ of \ } (I+\widetilde{K}_0(z) )^{-1} \mathrm{\ in\ }\mathbb{C}^{\circ}\big\} \subset \Sigma_{\mathrm{hom}} \cup \Sigma_{\mathbb{R}^d\setminus \overline{\Omega}}.
\end{equation}
\end{lemma}
    
\begin{proof}
By the same argument as Step 1 of the proof of Lemma \ref{pole_K01}, we can prove that $0$ is not a pole of 
$( I+\widetilde{K}_0  (z) )^{-1}$. We omit the details. Hereafter we assume that $z\in \mathbb{C}^{\circ}\setminus \big(\Sigma_{\mathbb{R}^d\setminus \overline{\Omega}}\cup\{0\} \big)$ is a pole of $( I+\widetilde{K}_0  (z) )^{-1}$, our goal is to show that $z \in \Sigma_{\mathrm{hom}}$. Since $\widetilde{K}_0(z)$ is compact, there exists a nonzero $g \in L^2(B_r)$ such that $(I+\widetilde{K}_0(z)) g =0$ in $B_r$. We divide the rest of the proof into two steps. 

\emph{Step 1}. We show $\widetilde{Q}_0(z)g\neq0$. By contradiction, suppose that $\widetilde{Q}_0(z)g =0$. Then, we have
\begin{equation} 
    \widetilde{W}_0(z)g = (L_0^r )^{-1} \big( I-\widetilde{T}_0(z)  \big) g =0 \quad \mathrm{in} \ \{\chi = 1\}\supset\supset\overline{\Omega} ,
\end{equation}
and
\begin{equation}
    \widetilde{R}_\Delta(z)g = 0 \qquad \mathrm{in} \ \{\chi = 0\} \supset \mathbb{R}^d\setminus \overline{B_r} .
\end{equation}
This implies that $ \widetilde{W}_0(z)g ,\,\widetilde{R}_\Delta(z)g  \in H^1_0(B_r\setminus \overline{\Omega})$, and $ g \equiv 0$ in $ \Omega$.

Let $u=\widetilde{R}_\Delta(z)g - \widetilde{W}_0(z) g$, by the definition of $ \widetilde{W}_0(z) $, we get
     \begin{equation} 
         \Delta u = \Delta \widetilde{R}_\Delta (z) g -( -z^2 \widetilde{R}_\Delta (z)g +g ) =0 \quad \mathrm{in}\ B_{r} \setminus \overline{\Omega } .
     \end{equation}
The above identity can be verified directly for $\mathrm{Im}\,z>0$, and extended to $\mathbb{C} \setminus \Sigma_{\mathbb{R}^d\setminus \overline{\Omega}} $ by holomorphic continuation. Therefore, $u=0$ in $B_{r} \setminus \overline{\Omega }$. Moreover, by definition of $\widetilde Q_0(z)$, we get
\begin{equation}
    \widetilde R_\Delta (z)g  = \widetilde Q_0 (z)g  + \chi u =0 \qquad  \mathrm{in} \ \mathbb{R}^d.
\end{equation}
Thus, $g\equiv 0$ in $\mathbb{R}^d$. This contradicts the assumption that $g$ is nonzero. Therefore, $\widetilde{Q}_0(z)g\neq0$. 

\emph{Step 2}. We show that $z\in \Sigma_{\mathrm{hom}}$. By \eqref{defR1ep11}, we have
\begin{equation}
    (L_0+z^2\mu^z_0) \widetilde{Q}_0(z)g= \mu_0^z(I+\widetilde{K}_0(z)) g  =0.
\end{equation}
To prove that $z\in \Sigma_{\mathrm{hom}}$, it remains to show that $\widetilde{Q}_0(z)g$ is a resonant state in the sense of Definition \ref{DefA2z}.
Take $ R>r $ and cutoff functions $ \rho_1 ,\rho\in C_c^\infty (\R^d;[0,1])$ such that $ \rho_1\equiv1 $ on $ B_R$ and $ \rho\equiv1 $ on $ \mathrm{supp}\, \rho_1$. Then, we have
\begin{equation}
    (\Delta+z^2)(1-\rho)\widetilde R_\Delta(z)\rho_1 = -[\Delta,\rho] \widetilde R_\Delta(z)\rho_1,
\end{equation}
which yields
\begin{equation}\label{5112R}
    -R_\Delta(z)[\Delta,\rho] \widetilde R_\Delta(z)\rho_1= (1-\rho)\widetilde R_\Delta(z)\rho_1 .
\end{equation}
The above identity can be verified directly for $\mathrm{Im}\,z>0$, and extended to $\mathbb{C} \setminus \Sigma_{\mathbb{R}^d\setminus \overline{\Omega}} $ by  meromorphic continuation. Define
\begin{equation}
     h := -[\Delta,\rho] \widetilde R_\Delta(z)\rho_1 g= -[\Delta,\rho] \widetilde R_\Delta(z) g,
\end{equation}
it follows from \eqref{5112R} that
\begin{equation}
    R_\Delta(z) h= (1-\rho)\widetilde R_\Delta(z)g = \widetilde  R_\Delta(z)g \qquad \mathrm{in} \ \mathbb{R}^d \setminus \mathrm{supp}\,\rho.
\end{equation}
Therefore,
\begin{equation}
    \widetilde{Q}_0(z) g =  \chi  \widetilde{ W}_0(z) g +( 1-\chi ) R_\Delta(z) h \qquad \mathrm{in} \ \mathbb{R}^d \setminus \mathrm{supp}\,\rho.
\end{equation}
By Theorem \ref{thmA3u}, $\widetilde{Q}_0(z) g$ is a resonant state of $L_0+z^2\mu^z_0$. Therefore, $z \in \Sigma_{\mathrm{hom}}$.
\end{proof}

We are in the position of giving the proof of Theorem \ref{mainresult1} and \ref{thm5121main}.

\hfill

\noindent \emph{Proof of Theorem \ref{thm5121main}}. Since $z\in \mathbb{C}\setminus (\Sigma_D \cup \Sigma_{\mathrm{hom}})$, and $\Sigma_r \cap \Sigma_{\mathbb{R}^d\setminus \overline{\Omega}} = \emptyset$, we have
\begin{equation}
    z \in \Big(\mathbb{C}\setminus (\Sigma_D \cup \Sigma_{\mathrm{hom}} \cup \Sigma_r)  \Big) \cup \Big(\mathbb{C}\setminus (\Sigma_D \cup \Sigma_{\mathrm{hom}} \cup \Sigma_{\mathbb{R}^d\setminus \overline{\Omega}} )  \Big).
\end{equation}

For $z \in  \mathbb{C}\setminus (\Sigma_D \cup \Sigma_{\mathrm{hom}} \cup \Sigma_r)  $, by Lemma \ref{pole_K01}, $(I+K_0(z))^{-1}$ exists. The uniform boundedness of $\|R_{\varepsilon}(z)\|_{L^2(B_r)\rightarrow L^2(B_r)}$ follows immediately from Theorem \ref{thmRBR} and \eqref{lapDep}. 

For $z\in \mathbb{C}\setminus (\Sigma_D \cup \Sigma_{\mathrm{hom}} \cup \Sigma_{\mathbb{R}^d\setminus \overline{\Omega}} )$, by Lemma \ref{pole_K0116}, $(I+\widetilde{K}_0(z))^{-1}$ exists. The uniform boundedness of $\|R_{\varepsilon}(z)\|_{L^2(B_r)\rightarrow L^2(B_r)}$ follows immediately from Theorem \ref{thm511Rtilde} and \eqref{lapDep}. The proof is complete.  \hfill $\square$

\hfill

\noindent \emph{Proof of Theorem \ref{mainresult1}}. 
By Proposition \ref{propSigmahomnega}, $k \notin \Sigma_{\mathrm{hom}} \cup \Sigma_D $. It follows from Theorem \ref{thm5121main} that $\|R_{\varepsilon}(z) \|_{L^2(B_r)\rightarrow L^2(B_r)} \leq C$ for all $\varepsilon\in (0,1)$.
By Theorem \ref{blethm}, $\mathcal{E}(\widehat{u}_0) \leq C\varepsilon^{1/2} \| \uinc \|_{L^2(B_{r+1})}$. 
Therefore, by Theorem \ref{L2velocity}, we have
\begin{equation}
\begin{aligned}
    \| u_\varepsilon - \Lambda_\varepsilon \widehat{u}_0 \|_{L^2(B_r)} & \leq \| \varepsilon \chi^{\varepsilon}   \cdot\eta_{\varepsilon} S_\varepsilon ( \nabla \widehat{u}_0 )   \|_{L^2(B_r)} +C\mathcal{E}(\widehat{u}_0)  (1+\| R_{\varepsilon}(k) \|_{L^2(B_r)\rightarrow L^2(B_r)} ) \\
    & \leq C\varepsilon^{1/2} \| \uinc \|_{L^2(B_{r+1})}.
\end{aligned}
\end{equation}
We get \eqref{L2formulaLam}. For \eqref{L2formulaLam11}, using the same argument as above, except using Theorem \ref{H1_velocity} instead of Theorem \ref{L2velocity}.
\hfill $\square$

\subsection{\texorpdfstring{Proof of Theorem \ref{maintheoremresonance}: convergence rate of the scattering resonances}{Poles of R0+S}}

\hfill

\noindent \emph{Proof of Theorem \ref{maintheoremresonance}}. 
We divide the proof into three steps.

\emph{Step 1}. In this step, we assume that $z_0\in \Sigma_{D,0}$. Then, $\mu_0^z$ is holomorphic at $z_0$. By Lemma \ref{pole_K0116}, $(I+\widetilde{K}_0(z) )^{-1}$ is holomorphic at $z_0$. Since $\Sigma_{D,0}\subset \mathbb{R}$, $(L_0+z^2\mu_0^z)^{-1} $ is holomorphic at $z_0$. Choose $r$ sufficiently small such that the circle $\Gamma:=\{z\in \mathbb{C}:|z-z_0|=r\}$ contains no other point in $\Sigma_{\mathrm{hom}} \cup \Sigma_D$. By the definition \eqref{defR1ep11} of $R^{(2)}_{\varepsilon}(z)$, for any $0<\delta<1$, we have
\begin{equation}
    \oint_{\delta\Gamma} R^{(2)}_{\varepsilon}(z)\,dz = \mathbbm{1}_{\{z_0\}} \big( (\varepsilon^2\Delta_{D_{\varepsilon}} +z^2)^{-1}  \big)(I+\widetilde{K}_0(z_0) )^{-1},
\end{equation}
where $\delta \Gamma:= \{z\in \mathbb{C}:|z-z_0|=\delta r\}$, and $\mathbbm{1}_A(T)$ denotes the spectral projection of the operator $T$ onto the Borel measurable set $A$. In particular, 
\begin{equation}
    \left\| \oint_{\delta \Gamma} R^{(2)}_{\varepsilon}(z)\,dz \right\|_{r,R} \geq  \frac{1}{\| I+\widetilde{K}_0(z_0)  \|_{r,R} }>0, \qquad \forall \varepsilon>0.
\end{equation}
We apply Theorem \ref{thm511Rtilde} to obtain
\begin{equation}\label{con5711123}
\begin{aligned}
    \left\| \oint_{\delta\Gamma}R_{\varepsilon}(z) \,dz \right\|_{r,R} &\geq \left\| \oint_{\delta\Gamma}R^{(2)}_{\varepsilon}(z) \,dz \right\|_{r,R}  -C\varepsilon^{1/2}  \left(1+  \frac{ 1 }{ \delta }  \right)^2 \frac{\delta}{1- \varepsilon^{1/2}  (1+  \delta^{-1} )   } \\
    &\geq C-C\varepsilon^{1/2}  \left(1+  \frac{ 1 }{ \delta }  \right)^2 \frac{\delta}{1- \varepsilon^{1/2}  (1+  \delta^{-1} )   } , \qquad \forall \varepsilon>0.
\end{aligned}
\end{equation}
From this, we conclude that for any $\delta = \varepsilon^{\alpha}$, $0\leq \alpha<1/2$,
\begin{equation}
    \left\| \oint_{\delta\Gamma}R_{\varepsilon}(z) \,dz \right\|_{r,R} \geq C>0,
\end{equation}
which means there are poles of $R_{\varepsilon}(z) $ contained in $\delta\Gamma$. Therefore, there exists $\varepsilon_i\rightarrow 0$ and resonance $z_i $ of $R_{\varepsilon_i}(z)$ such that $|z_i -z_0|<C\varepsilon^{\frac{1}{2}-} $.

\emph{Step 2}. Define
\begin{equation}
    R^{(1)}_0(z):= \big( Q_0(z) + \beta(z) \mathbbm{1}_{\Omega}\big)\big( I+K_0(z) \big)^{-1}.
\end{equation}
By Proposition \ref{thm49}, we have
\begin{equation}\label{5136}
    \|R^{(1)}_\varepsilon(z) -R^{(1)}_0(z) \|_{L^2(B_r) \rightarrow H^{-1}(B_R)} \leq \frac{C}{s} \varepsilon^{1/2} \big\| \big( I+K_0(z) \big)^{-1}  \big\|_r .
\end{equation}

We show that
\begin{equation}
    \Sigma_{\mathrm{hom}} \subset \big\{\textrm{Poles of } R^{(1)}_0(z) \textrm{ in }\mathbb{C}^{\circ}\big\} .
\end{equation}
By contradiction, assume that $z_0 \in \Sigma_{\mathrm{hom}} $ and not a pole of $R^{(1)}_0(z) $. By the identity $Q_0(z)(I+K_0(z))^{-1} = (L_0+z^2\mu_0^z)^{-1}\mu_0^z $, $z_0$ is a pole of $(I+K_0(z))^{-1}$. Since $Q_0(z)$ and $\beta(z)$ are holomorphic, by the Laurent expansion, there exists an integer $m>0$ and a nonzero finite-rank operator $A:L^2_{\mathrm{comp}}(\mathbb{R}^d) \rightarrow L^2_{\mathrm{loc}}(\mathbb{R}^d)$ such that for any $g\in \mathrm{ran}\,A$, we have
\begin{align}
    & (Q_0(z_0)+\beta(z_0)\mathbbm{1}_{\Omega})g=0,  \label{ranA1}\\
    & (I+K_0(z_0))g =0  \label{ranA2}.
\end{align}
We have $\mathrm{Im}\, z_0 <0$ by Proposition \ref{propSigmahomnega}. Then, by a direct computation, we get that $\beta(z_0)\neq 0$. Therefore, \eqref{ranA1} implies that $g\in H_0^1(\Omega) $. By the definitions \eqref{DefQ0z}-\eqref{DefK0} of $Q_0(z)$ and $K_0(z)$, we have
\begin{align}
    & W_0(z_0) g + \beta(z_0) g =0 && \mathrm{in}\ \Omega, \\
    & g +z_0^2W_0(z_0) g =0 && \mathrm{in}\ \Omega.
\end{align}
These imply that $\mu_0^{z_0} g = g - z_0^2 \beta(z_0) g =0 $ in $\Omega$. Since $\mathrm{Im}\,z_0<0$, $\mu_0^{z_0}\neq 0$. Therefore, $g \equiv 0$ in $\mathbb{R}^d$, we get a contradiction.

\emph{Step 3}. We now assume that $z_0\in \Sigma_{\mathrm{hom}}$. By Step 2, $z_0$ is a pole of $R^{(1)}_0(z)$ of order $n$, where $n\geq 1$ is an integer. Then, $z_0$ is a pole of $(I+K_0(z))^{-1}$ of order $m$, where $m\geq n$ is an integer.
Choose $r$ sufficiently small such that the circle $\Gamma:=\{z\in \mathbb{C}:|z-z_0|=r\}$ contains no other point in $R^{(1)}_0(z)$. Then, for any $0<\delta<1$, we have
\begin{equation}\label{5147}
    \left\| \oint_{\delta \Gamma} (z-z_0)^{n-1} R^{(1)}_0(z)\,dz \right\|_{L^2(B_r) \rightarrow H^{-1}(B_R)}= C>0,
\end{equation}
and
\begin{equation}\label{5148}
    \begin{aligned}
        &\left\| \oint_{\delta \Gamma} (z-z_0)^{n-1} R_{\varepsilon} (z)\,dz \right\|_{L^2(B_r) \rightarrow H^{-1}(B_R)} \\
        & \geq  \left\| \oint_{\delta \Gamma} (z-z_0)^{n-1} R^{(1)}_0(z)\,dz \right\|_{L^2(B_r) \rightarrow H^{-1}(B_R)} \\
        &\quad\, - \left\| \oint_{\delta \Gamma} (z-z_0)^{n-1} ( R^{(1)}_{\varepsilon} - R^{(1)}_0)(z)\,dz \right\|_{L^2(B_r) \rightarrow H^{-1}(B_R)} \\
        &\quad\, - \left\| \oint_{\delta \Gamma} (z-z_0)^{n-1} ( R_{\varepsilon} - R^{(1)}_{\varepsilon})(z)\,dz \right\|_{r,R} .
    \end{aligned}
\end{equation}
By \eqref{5136}, we have
\begin{equation}
    \left\| \oint_{\delta \Gamma} (z-z_0)^{n-1} ( R^{(1)}_{\varepsilon} - R^{(1)}_0)(z)\,dz \right\|_{L^2(B_r) \rightarrow H^{-1}(B_R)} \leq C\delta^{n-m-1} \varepsilon^{1/2}.
\end{equation}
By Theorem \ref{thmRBR}, we have
\begin{equation}\label{5150}
    \left\| \oint_{\delta \Gamma} (z-z_0)^{n-1} ( R_{\varepsilon} - R^{(1)}_{\varepsilon})(z)\,dz \right\|_{r,R} \leq C\frac{\varepsilon^{1/2}\delta^{n-2m-2}}{1-\varepsilon^{1/2} \delta^{-m-1}}.
\end{equation}
From \eqref{5147}-\eqref{5150}, we conclude that for any $\delta = \varepsilon^{\alpha}$, $0\leq \alpha<\frac{1}{2(2m+2-n)}$,
\begin{equation}
    \left\| \oint_{\delta \Gamma} (z-z_0)^{n-1} R_{\varepsilon} (z)\,dz \right\|_{L^2(B_r) \rightarrow H^{-1}(B_R)}\geq C>0,
\end{equation}
which means there are poles, of order at least $n$, of $R_{\varepsilon}(z) $ contained in $\delta\Gamma$. Therefore, there exists $\varepsilon_i\rightarrow 0$ and resonance $z_i $ of $R_{\varepsilon_i}(z)$ such that $|z_i -z_0|<C\varepsilon^{\nu-} $.  \hfill $\square$

\section{\texorpdfstring{Optimal $L^2$ convergence rate and the far-field pattern}{Optimal L2 rate for scattering problem}}\label{secoptmail}

\subsection{\texorpdfstring{Proof of Theorem \ref{mainresult2}: optimal $L^2$ convergence rate}{optimal rate}}

In this section, we use the dual method (\cite[Chapter 7]{shen2018periodic}) to establish the optimal $L^2$ convergence rate for the homogenization of the scattering problem \eqref{maineq}, i.e., Theorem \ref{mainresult2}. We need the following lemma.

\begin{lemma}\label{auxeq_convergevelocity}
Fix $k\in (0,\infty) \setminus \Sigma_D$. Let $ v_\varepsilon = R_\varepsilon(k)f$ and $\widehat{v}_0 = ( L_0 +k^2 \mu^k_0 )^{-1} \mu_0^k f$ for any $ f\in L^2( B_r )$. Then, for any $r>0$ such that $\overline{\Omega}\subset B_r$, we have
\begin{equation}
    \big\|v_\varepsilon- \Lambda_\varepsilon \widehat{v}_0-\varepsilon\chi^\varepsilon \cdot \eta_\varepsilon S_\varepsilon ( \nabla \widehat{v}_0 ) - ( \varepsilon^2 \Delta_{D_\varepsilon }+k^2 )^{-1} f \big\|_{\mathcal{H}_\varepsilon(B_r)} \leq C \varepsilon^{\frac{1}{2}}\|f\|_{L^2(B_r)},
\end{equation}
where $C>0$ is a constant that depends only on $d,k,r,D$ and $\Omega$.
\end{lemma}

\begin{proof}
Let 
\begin{equation}
    w_\varepsilon = v_\varepsilon - \Lambda_\varepsilon \widehat{v}_0 - \varepsilon \chi^{\varepsilon} \cdot \eta_{\varepsilon} S_{\varepsilon}(\nabla \widehat{v}_0) - (\varepsilon^2 \Delta_{D_{\varepsilon}} +k^2)^{-1} f .
\end{equation}
For simplicity, we denote $\mathcal{P}_{\varepsilon} (h|_{B_r\setminus \overline{D_{\varepsilon}}} )$ by $\mathcal{P}_{\varepsilon} h$ for a function $h:B_r\rightarrow \mathbb{C}$ below. We divide the rest of the proof into two steps.

\emph{Step 1}. we use the duality argument to estimate $\| w_{\varepsilon} \|_{L^2(B_r)}$. For any $ g\in L^2(B_r)$, take $ \psi_\varepsilon=R_\varepsilon(k)g $, then
\begin{equation}
    \begin{aligned}
        \int_{B_r}  w_\varepsilon g\,dx
        = \big\langle(L_\varepsilon+k^2)w_\varepsilon,\psi_\varepsilon \big\rangle_{H^{-1}(B_r),H^1(B_r)}  + \int_{\partial B_r} \frac{\partial \psi_\varepsilon}{\partial \mathbf{n}} w_\varepsilon - \frac{\partial w_\varepsilon}{\partial \mathbf{n}} \psi_\varepsilon\,dS(x),
    \end{aligned}
\end{equation}
and by \eqref{smiden} the last term on the right-hand side vanishes. By Lemma \ref{basiclemmasca} and a similar argument in the proof of Theorem \ref{L2velocity}, we obtain
\begin{equation}\label{aux_velocity_eq1}
    \begin{aligned}
        \left|\int_{B_r} w_\varepsilon g\, dx \right| & \leq C \mathcal{E}(\widehat{v}_0)\|g\|_{L^2(B_r)}  (1+\| R_{\varepsilon}(k) \|_{L^2(B_r) } )\\
&\quad\,  + \left| \int_{D_\varepsilon} \varepsilon^2 \nabla(\varepsilon^2 \Delta_{D_{\varepsilon}} +k^2)^{-1} f \cdot \nabla \mathcal{P} \psi_\varepsilon\,dx \right| \\
&\quad\,  + k^2 \left|\int_{ \Omega}\big\{ (\varepsilon^2 \Delta_{D_{\varepsilon}} +k^2)^{-1}    - \beta(k)  \big\} f\, \mathcal{P} \psi_\varepsilon \,dx \right|.
    \end{aligned}
\end{equation}
Using Theorem \ref{thm5121main}, Corollary \ref{col46L}, Lemma \ref{estimateforRepsilon} and Proposition \ref{thm49}, we get
\begin{equation}
    \begin{aligned}
        \left|\int_{B_r} w_\varepsilon g\, dx \right| & \leq C \varepsilon^{1/2} \|f \|_{L^2(B_r)} \|g \|_{L^2(B_r)} + \varepsilon^2 \| \nabla(\varepsilon^2 \Delta_{D_{\varepsilon}} +k^2)^{-1}  f\|_{L^2(D_{\varepsilon})} \| \nabla \mathcal{P} \psi_{\varepsilon} \|_{L^2(D_{\varepsilon})}\\
&\quad\,  + k^2\| (\varepsilon^2 \Delta_{D_{\varepsilon}} +k^2)^{-1}  - \beta(k) \|_{L^2(\Omega) \rightarrow H^{-1}(\Omega)} \| \mathcal{P }\psi_\varepsilon \|_{H^1(\Omega)} \\
& \leq C \varepsilon^{1/2} \|f \|_{L^2(B_r)} \|g \|_{L^2(B_r)} .
    \end{aligned}
\end{equation}
Since $ g $ is chosen arbitrarily, we have
\begin{equation}\label{step1wepsL2}
    \|w_\varepsilon \|_{L^2(B_r) } \leq C\varepsilon^{ \frac{1}{2}} \| f\|_{L^2(B_r)} .
\end{equation}

\emph{Step 2}. Next, we estimate the $ \mathcal{H}_\varepsilon(B_r) $ norm of $ w_\varepsilon $. Using \eqref{Renega}, we have
\begin{equation}\label{aux_velocity_eq3}
    \begin{aligned}
        &\|\nabla w_\varepsilon\|^2_{L^2(B_r\setminus \overline{D_\varepsilon} )} \\
        & \leq \mathrm{Re} \left\{ -\big\langle(L_\varepsilon+k^2)w_\varepsilon,\mathcal{P} w_\varepsilon \big\rangle_{H^{-1}(B_r),H^1(B_r)}  + k^2\int_{B_r} w_\varepsilon \overline{\mathcal{P}w_\varepsilon}\,dx
        - \varepsilon^2\int_{D_\varepsilon} \nabla w_\varepsilon \cdot \overline{\nabla\mathcal{P} w_\varepsilon}\,dx \right\} \\
        & \leq \left| \big\langle(L_\varepsilon+k^2)w_\varepsilon,\mathcal{P} w_\varepsilon \big\rangle_{H^{-1}(B_r),H^1(B_r)}  \right| + k^2 \|w_\varepsilon\|_{L^2(B_r)} \| \mathcal{P}_{\varepsilon} w_{\varepsilon} \|_{L^2(B_r)} \\
        &\quad\,+ \varepsilon^2  \|\nabla w_\varepsilon\|_{L^2(D_\varepsilon) }   \|\nabla \mathcal{P} w_\varepsilon\|_{L^2( D_\varepsilon)}.
    \end{aligned}
\end{equation}
For the first term on the right-hand side of the last line, we take $ \psi=\varphi =\overline{\mathcal{P}\psi_\varepsilon} $ in Lemma \ref{basiclemmasca}, using an estimate similar to that in Theorem \ref{H1_velocity}, to obtain
\begin{equation}
    \begin{aligned}
         \left| \big\langle(L_\varepsilon+k^2)w_\varepsilon,\mathcal{P} w_\varepsilon \big\rangle_{H^{-1}(B_r),H^1(B_r)} \right|& \leq
         C\varepsilon^{\frac{1}{2}} \|f \|_{L^2(B_r)} \|w_\varepsilon\|_{\mathcal{H}_\varepsilon(\Omega)} .
    \end{aligned}
\end{equation}
Substituting this into \eqref{aux_velocity_eq3}, we get
\begin{equation}
    \begin{aligned}
        \|\nabla w_\varepsilon\|^2_{L^2(B_r\setminus D_\varepsilon)} 
         &\leq C \varepsilon^{\frac{1}{2}} \|f \|_{L^2(B_r)} \|w_\varepsilon\|_{\mathcal{H}_\varepsilon(\Omega)} + C\|w_\varepsilon\|_{L^2(B_r)}^2 \\
         &\quad\, +C\big( \|w_\varepsilon \|_{L^2(B_r)} + \varepsilon^2 \|\nabla w_\varepsilon\|_{L^2(D_\varepsilon) }  \big) \|\nabla w_\varepsilon\|_{L^2(B_r\setminus D_\varepsilon)} .
    \end{aligned}
\end{equation}
By the Cauchy-Schwarz inequality and the estimate of $\|w_{\varepsilon}\|_{L^2(B_r)}$ in Step 1, we get
\begin{equation}\label{aux_converge1}
    \begin{aligned}
        \|\nabla w_\varepsilon\|_{L^2(B_r\setminus \overline{D_\varepsilon})} \leq C\varepsilon^{\frac{1}{2}} \|f\|_{L^2(B_r)}+C\varepsilon^2\|\nabla w_\varepsilon \|_{L^2(D_\varepsilon) }.
    \end{aligned}
\end{equation}
Let $ \overline\xi= w_\varepsilon-\mathcal{P}_\varepsilon w_\varepsilon\in H^1_0(D_\varepsilon) $, by a direct computation we get
\begin{equation}
    \begin{aligned}
       & \varepsilon^2\int_{D_\varepsilon} \nabla w_\varepsilon \cdot \nabla \xi\, dx\\
     & =\int_{D_\varepsilon } k^2 \big( w_\varepsilon  + \varepsilon v_{1,\varepsilon} ) \xi  -\varepsilon^2 (\Lambda_\varepsilon \nabla \widehat{v}_0+\varepsilon \nabla v_{1,\varepsilon})\cdot\nabla \xi  + \varepsilon (\nabla \Lambda)^\varepsilon\cdot \nabla \widehat{v}_0\,\xi \,dx,
    \end{aligned}
\end{equation}
where $ v_{1,\varepsilon} =\chi^{\varepsilon} \cdot \eta_{\varepsilon} S_{\varepsilon}(\nabla \widehat{v}_0)$. Using the Poincar\'{e} inequality and a similar argument to the proof of Theorem \ref{H1_velocity}, we then obtain
\begin{equation}
\begin{aligned}
     \varepsilon\|\nabla \xi \|_{L^2(D_\varepsilon)} &\leq C \varepsilon \|\nabla \mathcal{ P}_\varepsilon w_\varepsilon\|_{L^2(D_\varepsilon)}+ C\|w_\varepsilon \|_{L^2(D_\varepsilon ) } 
        + C\varepsilon^{\frac{1}{2}} \mathcal{E}(\widehat{v}_0) \\
        & \leq C\varepsilon \|\nabla w_\varepsilon \|_{L^2(\Omega\setminus D_\varepsilon)} + C\varepsilon^{\frac{1}{2}} \|f \|_{L^2(B_r)} ,
\end{aligned}
\end{equation}
Combining this with \eqref{aux_converge1} yields
\begin{equation}\label{613}
    \varepsilon\|\nabla w_\varepsilon\|_{L^2(D_\varepsilon)} \leq C \varepsilon \|\nabla w_\varepsilon\|_{L^2(\Omega\setminus \overline{D}_\varepsilon)}+ C\varepsilon^{ \frac{1}{2}} \|f\|_{L^2(\Omega)}
    \leq C \varepsilon^{ \frac{1}{2}} \|f\|_{L^2(\Omega)} +C\varepsilon^3\|\nabla w_\varepsilon \|_{L^2(D_\varepsilon) }.
\end{equation}
The desired conclusion follows by combining \eqref{step1wepsL2}, \eqref{aux_converge1} and \eqref{613}.
\end{proof}

\noindent \emph{Proof of Theorem \ref{mainresult2}}. Denote 
\begin{equation}
     z_\varepsilon=u_\varepsilon-\Lambda_\varepsilon \widehat{u}_0 - \varepsilon \chi^{\varepsilon}  \cdot\eta_{\varepsilon} S_\varepsilon ( \nabla \widehat{u}_0 ) ,
\end{equation}
$ v_\varepsilon= R_\varepsilon(k)f$, $\widehat{v}_0=(L_0+k^2\mu^k_0)^{-1} \mu^k_0 f $, $ v_{1,\varepsilon}= \chi^\varepsilon \cdot \eta_\varepsilon S_\varepsilon(\nabla \widehat{v}_0) $, and
\begin{equation}
    w_\varepsilon= v_\varepsilon-\Lambda_\varepsilon\widehat{v}_0 - \varepsilon v_{1,\varepsilon}-( \varepsilon^2\Delta +k^2 )^{-1} f .
\end{equation}
By \eqref{smiden}, we have
\begin{equation}\label{6171}
    \begin{aligned}   
         \int_{B_r} z_\varepsilon f\,dx  &=\langle (L_\varepsilon+k^2)z_\varepsilon , v_\varepsilon\rangle_{ H^{-1}(B_r) , H^1(B_r) }\\ 
     & = \big\langle (L_\varepsilon+k^2)z_\varepsilon , w_\varepsilon \big\rangle_{ H^{-1}(B_r) , H^1(B_r) } \\
    &\quad\, + \big\langle (L_\varepsilon+k^2)z_\varepsilon , \Lambda_\varepsilon \widehat{v}_0 \big\rangle_{ H^{-1}(B_r) , H^1(B_r) } \\
    & \quad\,  + \big\langle (L_\varepsilon+k^2)z_\varepsilon ,  \varepsilon v_{1,\varepsilon} \big\rangle_{ H^{-1}(B_r) , H^1(B_r) } \\
    & \quad\,  + \big\langle (L_\varepsilon+k^2)z_\varepsilon ,  (\varepsilon^2\Delta_{D_{\varepsilon}} +k^2)^{-1}f \big\rangle_{ H^{-1}(B_r) , H^1(B_r) }.
    \end{aligned}
\end{equation}
We divide the rest of the proof into four steps. 

\emph{Step 1}. For the first term on the right-hand side of \eqref{6171}, since 
\begin{equation}
     \left\|\frac{\partial z_\varepsilon}{\partial \mathbf{n}} \right\|_{L^2(\partial B_r)} \|w_\varepsilon\|_{L^2(\partial B_r)} \leq C \|z_\varepsilon\|_{H^2(B_r\setminus \overline{B_{r-1}})} \|w_\varepsilon\|_{ H^1( B_r\setminus \overline{ B_{r-1}})} ,
\end{equation}
by standard elliptic estimates and Theorem \ref{H1_velocity}, Corollary \ref{auxeq_convergevelocity}, Theorem \ref{blethm}, we get
\begin{equation}
    \begin{aligned}
        \left| \big\langle (L_\varepsilon+k^2)z_\varepsilon , w_\varepsilon \big\rangle_{ H^{-1}(B_r) , H^1(B_r) } \right|\leq C \|z_\varepsilon\|_{\mathcal{H}_\varepsilon(B_r)}\|w_\varepsilon\|_{\mathcal{H}_\varepsilon(B_r)}\leq C\varepsilon \|\uinc\|_{L^2(B_{r+1}) } \|f\|_{L^2(B_r)} .
    \end{aligned}
\end{equation}

\emph{Step 2}. We estimate the second term on the right-hand side of \eqref{6171}. Given $ \phi\in C_0^{\infty}(\R^d) $, we have
\begin{equation}\label{velocityimproved1}
\begin{aligned}
    &\varepsilon \big\|\nabla \big(\eta_{\varepsilon}  S_{\varepsilon} (\nabla\widehat{u}_0) \big) \nabla\phi\big\|_{L^1(\Omega)} \\
    & \leq C\|\nabla\widehat{u}_0\|_{ L^2(\mathcal{O }_{(d+2)\varepsilon})} \|\nabla\phi\|_{L^2(\mathcal{O}_{ (d+1)\varepsilon })} + C\varepsilon \| \nabla^2\widehat{u}_0 \|_{L^2(\Omega\setminus \mathcal{O}_{ (d-1)\varepsilon } )} \|\nabla\phi\|_{L^2(\Omega)}\\
    & \leq C \varepsilon^{ \frac{1}{2} }\|\widehat{u}_0\|_{H^2(\Omega) } \big( \| \nabla  \phi \|_{L^2(\mathcal{O}_{ (d+1)\varepsilon })} + \varepsilon^{\frac{1}{2} }  \|\nabla\phi\|_{L^2(\Omega)} \big),
\end{aligned}
\end{equation}
and the following estimate which is similar to \eqref{Stack}:
\begin{equation}\label{velocityimproved2}
\begin{aligned}
     &\left| \int_{\Omega}  (\mu_0^k -\Lambda_{\varepsilon} )\widehat{u}_0  \phi\,dx\right| \\
     &\leq C\left| \sum_{\mathbf{m}\in J_\varepsilon} \int_{\varepsilon (Y+\mathbf{m})}  (\mu_0^k -\Lambda_{\varepsilon} )\widehat{u}_0  \phi \,dx\right| + C\| \Lambda \|_{L^\infty(Y ) } \|\widehat{u}_0\|_{L^2(\mathcal{O}_{d\varepsilon})}  \| \phi \|_{L^2( \mathcal{ O}_{d\varepsilon}) } \\
     & \leq C \varepsilon\| \widehat{ u }_0 \|_{H^1( \Omega )}\|\phi \|_{H^1(\Omega)},
\end{aligned}
\end{equation}
where in the above two estimates, we used Lemma \ref{blelem}.

Taking $\overline\psi=\Lambda_\varepsilon \widehat{v}_0$ and $\overline\varphi=\widehat{v}_0$ in Lemma~\ref{basiclemmasca}, and following the proof of Theorem~\ref{L2velocity} with the estimates therein replaced, where appropriate, by \eqref{velocityimproved1}–\eqref{velocityimproved2}, and invoking Lemma~\ref{propertySepsilondifference}, we obtain
\begin{equation}\label{OptmL2velocity_eq1}
    \begin{aligned}
        &\left|\big\langle ( L_{\varepsilon} +k^2) z_{\varepsilon}, \Lambda_\varepsilon \widehat{v}_0 \big\rangle_{H^{-1}(B_r),H^1(B_r)}\right| \leq C \varepsilon\|\widehat{u}_0\|_{H^2 (\Omega) }\| \widehat{v}_0\|_{H^2(\Omega)} \\
        & + \left| \int_{D_\varepsilon} - k^2 \Lambda_{\varepsilon} \widehat{u}_0 ( \Lambda_\varepsilon-1 ) \widehat{v}_0\,dx  
        + \int_{D_\varepsilon} \varepsilon (\nabla\Lambda)^{\varepsilon} \widehat{u}_0\cdot  \nabla \big( ( \Lambda_\varepsilon-1 ) \widehat{v}_0 \big) \,dx \right| .
    \end{aligned}
\end{equation}
Since $ \Lambda_\varepsilon-1\in H^1_0(D_\varepsilon) $, we have
\begin{equation}
\begin{aligned}
     &\int_{D_\varepsilon} \varepsilon (\nabla\Lambda)^{\varepsilon} \widehat{u}_0\cdot  \nabla \big( ( \Lambda_\varepsilon-1 ) \widehat{v}_0 \big) \,dx \\
     &= \int_{D_\varepsilon} k^2 \Lambda_{\varepsilon} \widehat{u}_0 ( \Lambda_\varepsilon-1 ) \widehat{v}_0 - \varepsilon (\nabla \Lambda)^\varepsilon \cdot \nabla\widehat{u}_0 (\Lambda_\varepsilon-1) \widehat{v}_0 \,dx,
\end{aligned}
\end{equation}
Thus, by \eqref{OptmL2velocity_eq1}, we get
\begin{equation}
    \begin{aligned}
        \left|\big\langle ( L_{\varepsilon} +k^2) z_{\varepsilon}, \Lambda_\varepsilon \widehat{v}_0 \big\rangle_{H^{-1}(B_r),H^1(B_r)}\right| \leq C\varepsilon\|\widehat{u}_0\|_{H^2 (\Omega) }\| \widehat{v}_0\|_{H^2(\Omega)}.
    \end{aligned}
\end{equation}

\emph{Step 3}. Taking $ \overline\psi=\overline{\varphi}=\varepsilon v_{1,\varepsilon} $ in Lemma \ref{basiclemmasca}, by \eqref{velocityimproved1}, we have
\begin{equation}
\begin{aligned}
    &\varepsilon \big\|\nabla \big(\eta_{\varepsilon}  S_{\varepsilon} (\nabla\widehat{u}_0) \big)\cdot \nabla \varepsilon v_{1,\varepsilon} \big\|_{L^1(\Omega)} \\
    & \leq C \varepsilon^{ \frac{3}{2} }\|\widehat{u}_0\|_{H^2(\Omega) } \big( \| \nabla v_{1,\varepsilon} \|_{L^2(\mathcal{O}_{ (d+1)\varepsilon })} + \varepsilon^{\frac{1}{2} } \|\nabla v_{1,\varepsilon} \|_{L^2(\Omega)} \big).
\end{aligned}
\end{equation}
Then, similar to Step 2, we get
\begin{equation}
    \begin{aligned}
        &\left|\big\langle ( L_{\varepsilon} +k^2) z_{\varepsilon},\varepsilon v_{1,\varepsilon} \big\rangle_{H^{-1}(B_r),H^1(B_r)} \right| \\
        & \leq C \varepsilon^{ \frac{3}{2} }\|\widehat{u}_0\|_{H^2(\Omega) } \big(\|\nabla v_{1,\varepsilon} \|_{ L^2(\mathcal{O }_{(d+1) \varepsilon})}+ \varepsilon^{\frac{1}{2}}\| v _{1,\varepsilon}\|_{H^1(\Omega)} \big).
    \end{aligned}
\end{equation}
With
\begin{equation}
    \varepsilon^{\frac{1}{2} }\|\nabla v_{1,\varepsilon}\|_{L^2( \mathcal{O}_{ (d+1)\varepsilon} )} \leq C\varepsilon^{\frac{1}{2}} \| \nabla^2 \widehat{v}_0 \|_{L^2(\Omega)} +C \varepsilon^{ -\frac{1}{2} }\| \nabla\widehat{v}_0 \|_{L^2( \mathcal{O}_{(d+1)\varepsilon} )} \leq C \| \widehat{v}_0 \|_{H^2(\Omega)},
\end{equation}
and
\begin{equation}
    \begin{aligned}
        \varepsilon \|v_{1,\varepsilon}\|_{H^1(\Omega)} \leq C \| \nabla \widehat{v}_0 \|_{ L^2( \Omega ) } + C\varepsilon \| \nabla^2 \widehat{v}_0 \|_{ L^2(\Omega) }\leq C \|\widehat{v}_0 \|_{H^2(\Omega)} ,
    \end{aligned}
\end{equation}
we get
\begin{equation}
        \left|\big\langle ( L_{\varepsilon} +k^2) z_{\varepsilon},\varepsilon v_{1,\varepsilon} \big\rangle_{H^{-1}(B_r),H^1(B_r)} \right|  \leq C\varepsilon\|\widehat{u}_0\|_{H^2 (\Omega) }\| \widehat{v}_0\|_{H^2(\Omega)}.
\end{equation}

\emph{Step 4}. Taking $ \overline\psi= (\varepsilon^2\Delta_{D_{\varepsilon}} +k^2)^{-1} f$ and $\varphi=0 $ in Lemma \ref{basiclemmasca}, we get
\begin{equation}
    \begin{aligned}
        &\big\langle ( L_{\varepsilon} +k^2) z_{\varepsilon}, \psi \big\rangle_{H^{-1}(B_r),H^1(B_r)} \\
        &= -k^2\int_{D_\varepsilon}  \Lambda_{\varepsilon} \widehat{u}_0 (\varepsilon^2 \Delta_{D_{\varepsilon}} +k^2)^{-1} f \,dx  
        + \int_{D_\varepsilon } \varepsilon (\nabla\Lambda)^\varepsilon \widehat{u}_0  \cdot  \nabla (\varepsilon^2 \Delta_{D_{\varepsilon}} +k^2)^{-1} f \,dx\\
        & \quad\,+ \int_{D_\varepsilon } \varepsilon^2 \big\{ \Lambda_{\varepsilon}\nabla\widehat{u}_0 +\varepsilon \nabla( \chi^{\varepsilon} \eta_{\varepsilon}  \cdot S_{\varepsilon} (\nabla\widehat{w}_0) ) \big\} \cdot  \nabla (\varepsilon^2 \Delta_{D_{\varepsilon}} +k^2)^{-1}f \,dx.
    \end{aligned}
\end{equation}
Since $ (\varepsilon^2 \Delta_{D_{\varepsilon}} +k^2)^{-1} f \in H^1_0(D_\varepsilon) $, we have
\begin{equation}
\begin{aligned}
    &\int_{D_\varepsilon} \varepsilon (\nabla\Lambda)^\varepsilon \widehat{u}_0  \cdot \nabla (\varepsilon^2 \Delta_{D_{\varepsilon}} +k^2)^{-1}f- k^2 \Lambda_{\varepsilon} \widehat{u}_0(\varepsilon^2 \Delta_{D_{\varepsilon}} +k^2)^{-1} f \,dx  \\
    &= -\varepsilon\int_{D_\varepsilon}  ( \nabla \Lambda)^\varepsilon \cdot \nabla \widehat{ u}_0 (\varepsilon^2 \Delta_{D_{\varepsilon}} +k^2)^{-1} f \, dx.
\end{aligned}
\end{equation}
Hence,
\begin{equation}
    \begin{aligned}
        &\left|\big\langle ( L_{\varepsilon} +k^2) z_{\varepsilon}, \psi \big\rangle_{H^{-1}(B_r),H^1(B_r)} \right| \\
        &\leq C 
        \varepsilon \big( \|\nabla \widehat{u}_0 \|_{L^2(D_\varepsilon)} +\varepsilon \| \nabla (\chi^{\varepsilon} \eta_{\varepsilon}  \cdot S_{\varepsilon} (\nabla\widehat{w}_0))\|_{ L^2(D_\varepsilon) }\big) \|f\|_{ L^2 (D_\varepsilon)}.
    \end{aligned}
\end{equation}

Combine all of the above, along with the standard elliptic estimate for $C^{1,1}$  domain $\Omega $:
\begin{equation}
     \|\widehat{v}_0\|_{H^2(\Omega)}\leq C \|f\|_{L^2(B_r)},
\end{equation}
we obtain the desired conclusion \eqref{251optimal}. \hfill$\square$

\subsection{Proof of Theorem \ref{thmamplitd}: convergence rate of the far-field pattern}

Let $u_{\varepsilon}$ be the solution of the scattering problem \eqref{maineq}. Let  $u^s_{\varepsilon} := u_{\varepsilon} - \uinc$ be the scattered wave. Then, $u^s_{\varepsilon}$ admits the asymptotics \eqref{633} at infinity, with far-field pattern given by
\begin{equation}\label{6341}
u^{\infty}_\varepsilon(\theta) = \left(\frac{k}{2}\right)^{\frac{d-3}{2}} \frac{e^{-\mathrm{i}\frac{\pi}{4}(d+1)} }{4\pi^{\frac{d-1}{2}}} \int_{ \partial B_r } \left( \mathrm{e}^{-\mathrm{i}k\theta \cdot y} \frac{\partial u^s_\varepsilon}{ \partial \mathbf{n}_y } ( y) - u^s_\varepsilon( y ) \frac{\partial (\mathrm{e}^{ -\mathrm{i}k \theta \cdot y })}{ \partial \mathbf{n}_y }  \right) dS(y).
\end{equation}
Here $R>0$ satisfies $\overline{\Omega} \subset B_{R-1}$. In fact, according to the Green's formula \eqref{Greenrep} and the representation \eqref{320Gk} of the fundamental solution $G^k$, \eqref{633} follows from the following asymptotics as $|x|,\,r\rightarrow \infty$:
\begin{align}
    & |x-y|= |x|\left(1- 2\theta\cdot\frac{y}{|x|}+\frac{|y|^2}{|x|^2} \right)^{\frac{1}{2}}= |x|- \theta\cdot y+O(|x|^{-1}), \\
    &\frac{x-y}{|x-y|}= \frac{\theta - |x|^{-1}y}{ 1+O(|x|^{-1}) }= \theta + O(|x|^{-1}), \\
    & H^1_{\frac{d-2}{2}}(kr) = \left(\frac{2}{\pi kr}\right)^{\frac{1}{2}} e^{\mathrm{i}kr -\mathrm{i}\frac{\pi}{4}(d-1) }+ O(r^{-3/2}) \\
    & ( H^1_{\frac{d-2}{2}})'(kr) = \mathrm{i}\left(\frac{2}{\pi kr}\right)^{\frac{1}{2}} e^{\mathrm{i}kr -\mathrm{i}\frac{\pi}{4}(d-1) } + O(r^{-3/2}).
\end{align}

The far-field pattern for the homogenized scattering problem \eqref{equ0} can be obtained similarly:
\begin{equation}\label{63410}
u^{\infty}_0(\theta) =\left(\frac{k}{2}\right)^{\frac{d-3}{2}} \frac{e^{-\mathrm{i}\frac{\pi}{4}(d+1)} }{4\pi^{\frac{d-1}{2}}} \int_{ \partial B_r } \left( \mathrm{e}^{-\mathrm{i}k\theta \cdot y} \frac{\partial u^s_0}{ \partial \mathbf{n}_y } ( y) - u^s_0( y ) \frac{\partial (\mathrm{e}^{ -\mathrm{i}k \theta \cdot y })}{ \partial \mathbf{n}_y }  \right) dS(y).
\end{equation}
where $u^s_0:=\widehat{u}_0 - \uinc $.

\noindent\emph{Proof of Theorem \ref{thmamplitd}}.
 By the formulas \eqref{6341} and \eqref{63410} for the far-field pattern, we have
\begin{equation}
\begin{aligned}
    &u^\infty_\varepsilon (\theta)- u^\infty_0 ( \theta ) \\
    &= \left(\frac{k}{2}\right)^{\frac{d-3}{2}} \frac{e^{-\mathrm{i}\frac{\pi}{4}(d+1)} }{4\pi^{\frac{d-1}{2}}} \int_{ \partial B_r } \left( \mathrm{e}^{-\mathrm{i}k\theta \cdot y} \frac{\partial (u_\varepsilon-\widehat{u}_0) }{ \partial \mathbf{n}_y } ( y) - (u_\varepsilon-\widehat{u}_0 )( y ) \frac{\partial (\mathrm{e}^{ -\mathrm{i}k \theta \cdot y })}{ \partial \mathbf{n}_y }  \right) dS(y),
\end{aligned}
\end{equation}
Since $( \Delta + k^2 )( u_{\varepsilon} - \widehat{u}_0)= 0 $ in $\mathbb{R}^d\setminus \overline{\Omega} $, by applying the trace theorem and the standard elliptic estimates, for any integer $n\geq 0$, we obtain
\begin{equation}
    \left|\partial_{\theta}^n\big( u^\infty_\varepsilon (\theta)- u^\infty_0 ( \theta ) \big)\right| 
 \leq C   \| u_\varepsilon - \widehat{u}_0  \|_{H^2( B_{r+1}\setminus \overline{B_r} )} .
\end{equation}
The proof is complete by applying Theorem \ref{mainresult1} and \ref{mainresult2}, and by noting Remark \ref{rem1rk}.
\hfill $\square$

\section{Concluding remark}\label{secrema}

In this paper, we have developed a comprehensive mathematical theory for the homogenization of scattering by a periodic array of high-contrast subwavelength resonators. We established quantitative convergence of both the scattered field and the scattering resonances to those of an effective model. A key ingredient is a new framework for the meromorphic continuation of the outgoing resolvent. The techniques introduced here would allow us to extend the results in the paper to a random setting and to more complex systems in fluid mechanics, electromagnetics, elasticity, and so on. We also expect our approach to be useful to obtain effective medium theories for other classes of high-contrast resonators, such as air-bubble arrays. These extensions will be the subject of forthcoming work.


\section*{Acknowledgments}
This work was partially supported by the NSFC Grant No.\,12571220 and by the New Cornerstone Investigator Program 100001127. 

\appendix

\section{Periodic unfolding method}\label{secunfold}

In 1990, Arbogast, Douglas and Hornung \cite{arbogast1990derivation} defined a ``dilation'' operation to study homogenization for a periodic medium with double porosity. It turns out that the dilation technique reduces two-scale convergence to weak convergence in an appropriate space. Combining this approach with ideas from Finite Element approximations, Cioranescu, Damlamian and Griso \cite{MR1921004} proposed the periodic unfolding method to study homogenization of multiscale periodic problems. For more details, we refer to \cite{allaire_bloch_1998, cioranescu_periodic_2018}. The most striking advantage of the periodic unfolding method is that it lifts two-scale convergence in $L^2(\Omega)$ to convergence in the unfolded space $L^2(\Omega \times Y)$. This refinement enables us to discuss the convergence rates.

For our purposes, we do not need to introduce the full machinery of the periodic unfolding method. We only need the so-called local averaging operator.

	\begin{definition}
		For $x \in \mathbb{R}^d$, there exists a unique element in $\mathbb{Z}^d$ denoted $[x]_Y$ such that
        \begin{equation}
            x = [x]_Y + \{x\}_Y , \qquad \{x\}_Y \in Y.
        \end{equation}
	The local averaging operator $\mathcal{U}_{\varepsilon} :L^2(\Omega \times Y) \rightarrow L^2(\Omega)$ is defined by
			\begin{equation}
				(\mathcal{U}_{\varepsilon} \phi ) (x) :=\sum_{\mathbf{m}\in J_{\varepsilon}} \mathbbm{1}_{\varepsilon(Y+\mathbf{m})} (x)\int_Y \phi \left( \varepsilon  \left[ \frac{x}{\varepsilon}\right]_Y+\varepsilon z,  \left\{ \frac{x}{\varepsilon}\right\} _Y \right)\,dz,\qquad x \in \Omega,
			\end{equation}
where $J_{\varepsilon}$ is defined in \eqref{jepsi}.  
\end{definition}

\begin{lemma}\label{fold}
    Let $\langle \cdot \rangle_Y : L^2(\Omega \times Y)\rightarrow L^2(\Omega)$ be the averaging operator:
    \begin{equation}
        \langle \cdot \rangle_Y : \phi \mapsto \int_Y \phi (\cdot,y)\,dy.
    \end{equation}
    We have
    \begin{equation}
        	\| \mathcal{U}_{\varepsilon} - \langle \cdot \rangle_Y \|_{L^2(\Omega \times Y) \rightarrow H^{-1}(\Omega)}\leq C \varepsilon^{\frac{1}{2}},
    \end{equation}
    where $C>0$ is a constant that depends only on $\Omega$.
\end{lemma}
\begin{proof}
    Given $\phi \in L^2(\Omega \times Y)$ and $g \in H_0^1(\Omega)$, we have
    \begin{equation}\label{(a)formula}
    \begin{aligned}
       & \left| \int_{\Omega} \left( \mathcal{U}_{\varepsilon} \phi (x)-  \int_Y \phi(x,s)\,ds \right)g(x)\,dx \right|\\
       &\leq C \sum_{\mathbf{m}\in J_{\varepsilon}} \left|\int_{\varepsilon(Y+\mathbf{m})} \left( \int_Y \phi \left( \varepsilon  \left[ \frac{x}{\varepsilon}\right]_Y+\varepsilon z,  \left\{ \frac{x}{\varepsilon}\right\} _Y \right)\,dz - \int_Y \phi (x,s)\,ds\right) g(x)\,dx \right| \\
       & \quad\, + C\int_{\mathcal{O}_{d\varepsilon}}  \int_Y |\phi (x,s)g(x)| \,ds dx.
    \end{aligned}
\end{equation}
Using Lemma \ref{blelem}, the last line above is bounded as follows:
\begin{equation}\label{451}
\int_{\mathcal{O}_{d\varepsilon}}  \int_Y |\phi (x,s)g(x)| \,ds dx \leq C \varepsilon^{1/2}\| \phi  \|_{L^2(\Omega \times Y)} \| g \|_{H^1(\Omega)}.
\end{equation}
We denote the remaining term on the right-hand side of \eqref{(a)formula} by $I$. From the change of variable $x=\varepsilon \mathbf{m}+\varepsilon y$, we have
\begin{equation}
\begin{aligned}
    I = \varepsilon^d    \sum_{\mathbf{m}\in J_{\varepsilon}}  \left| \int_Y \Phi^{\mathbf{m}}_{\varepsilon}(y) g  (\varepsilon \mathbf{m}+\varepsilon y) \,dy \right|,
\end{aligned}
\end{equation}
where
	\begin{equation}
				\Phi^{\mathbf{m}}_{\varepsilon}(y):=\int_Y \phi (\varepsilon {\mathbf{m}}+\varepsilon z,y)\,dz-\int_Y \phi (\varepsilon {\mathbf{m}}+\varepsilon y ,s)\,ds .
			\end{equation}
Note that $\int_Y \Phi^{\mathbf{m}}_{\varepsilon}(y)\,dy= 0$, so by the Poincar\'{e} inequality, we get
\begin{equation}\label{454}
\begin{aligned}
I& = \varepsilon^d   \sum_{\mathbf{m}\in J_{\varepsilon}}  \left| \int_Y \Phi^{\mathbf{m} }_{\varepsilon }(y)\left( g  (\varepsilon \mathbf{m} +\varepsilon y)  - \int_Y g  (\varepsilon \mathbf{m} +\varepsilon w) \,dw \right)\,dy \right|  \\
& \leq  \varepsilon^d   \sum_{\mathbf{m}\in J_{\varepsilon}}  \| \Phi^{\mathbf{m}}_{\varepsilon }\|_{L^2(Y)} \left\| g (\varepsilon \mathbf{m} +\varepsilon \cdot)  - \int_Y g  (\varepsilon \mathbf{m} +\varepsilon w) \,dw \right\|_{L^2(Y)} \\
& \leq C\varepsilon^{\frac{d}{2}+1}  \left( \sum_{\mathbf{m}\in J_{\varepsilon}}   \int_Y |\Phi^{\mathbf{m}}_{\varepsilon}(y)|^2 \,dy \right)^{1/2}  \| g \|_{H^1(\Omega)} .
\end{aligned}
\end{equation}
By the definition of $\Phi^{\mathbf{m}}_{\varepsilon}$, we have
			\begin{equation}
				\begin{aligned}
					|\Phi^{\mathbf{m}}_{\varepsilon}(y)|^2 \leq 2\int_Y |\phi (\varepsilon \mathbf{m}+\varepsilon z,y)|^2\,dz + 2\int_Y |\phi (\varepsilon \mathbf{m} +\varepsilon y,s)|^2\,ds,
				\end{aligned}
			\end{equation}
thus,
			\begin{equation}\label{456}
				\begin{aligned}
					 \sum_{\mathbf{m}\in J_{\varepsilon}}   \int_Y |\Phi^{\mathbf{m}}_{\varepsilon}(y)|^2 \,dy &\leq C  \sum_{\mathbf{m}\in J_{\varepsilon}}   \int_Y \int_Y   |\phi (\varepsilon \mathbf{m}+\varepsilon z,y)|^2\,dzdy \\
                     &\leq C  \varepsilon^{-d} \int_{\Omega \times Y} |\phi (x,y)|^2\,dxdy. 
				\end{aligned}
			\end{equation}
The proof is complete by combining \eqref{(a)formula}, \eqref{451}, \eqref{454} and \eqref{456} together.
\end{proof}

\section{Meromorphic continuation of the homogenized operator}\label{appendixA}

We first construct the meromorphic continuation of $L_0+z^2\mu_0^z $ to $\mathbb{C}^{\circ}$. Consider $L_0= \nabla\cdot A_0\nabla$ as an unbounded operator on $L^2(\mathbb{R}^d)$ with the domain
\begin{equation}\label{domainD0}
    \mathcal{D}(L_{0}) := \{ u \in  H^1(\mathbb{R}^d): \nabla \cdot A_{0} \nabla u \in L^2(\mathbb{R}^d) \} ,
\end{equation}
where $\nabla \cdot A_{0} \nabla u$ is understood in the sense of distribution. 

The following lemma is an analogue of Lemma \ref{lemmaresolvent}, we omit the proof.

\begin{lemma}\label{lemmaresolventx}
Let $r>0$ such that $\overline{\Omega} \subset B_r$. Let $\chi \in C_0^{\infty} (\mathbb{R}^d; [0,1])$ such that 
\begin{equation}
    \mathrm{supp}\,\chi \subset B_r \qquad \mathrm{and} \qquad \overline{\Omega} \subset \{\chi =1\}^{\circ}.
\end{equation} 
Define
\begin{equation}\label{defWx}
    Z_0(z):=(L_0^r)^{-1}\big(I-z^2R_\Delta (z)\big).
\end{equation}
Then, the following identity holds:
\begin{equation}\label{resolvent_equlity}
    \begin{aligned}
          (L_{0} + z^2\mu_0^{z}) D_{0}(z) = I + F_{0}(z) , \quad  \mathrm{Im} \, z >0,
    \end{aligned}
\end{equation}
where
\begin{equation}\label{defQx}
    D_{0}(z) := \chi Z_0(z) +  (1-\chi)R_{\Delta}(z),
\end{equation}
and 
\begin{equation}\label{defF}
 F_{0}(z) := \chi z^2 \big\{ \mu_0^{z}Z_0(z)-R_{\Delta } (z)\big\} + [\Delta,\chi] \big(Z_0(z) - R_{\Delta}(z) \big)  .
\end{equation}
\end{lemma}

\begin{lemma}\label{InverIF}
    $I + F_0(z)$ is invertible on $L^2(\mathbb{R}^d)$ for $\mathrm{Im}\,z>0$. 
\end{lemma}
\begin{proof}

It is clear that $F_0(z)$ is a compact operator on $L^2(\mathbb{R}^d)$. By the Fredholm theory, it suffices to show that $I+F_0(z)$ is injective. Assume that $(I+F_0(z))g =0$ for some $g \in L^2(\mathbb{R}^d)$. The rest of the proof is divided into two steps.

\emph{Step 1}. By a direct computation, we have
\begin{equation}\label{comz2}
    z^2 \mu_0^z = z^2 \left( 1-\theta - \sum_{j=1}^{\infty} \frac{\lambda_j^2}{z^2- \lambda_j^2} \left| \int_D \varphi_j(y)\,dy \right|^2  \right),
\end{equation}
where $\theta = |D|$ denotes the Lebesgue measure of $D$, and $\{(\varphi_j, \lambda_j^2)\}_{j=1}^{\infty}$ are the eigenpairs of $-\Delta_D$, with $\{\varphi_j\}_{j=1}^\infty$ forming an orthonormal basis of $L^2(D)$. If $z^2<0$, we get $z^2\mu_0^z<0$ by \eqref{comz2}. Taking the imaginary part of \eqref{comz2}, we get
\begin{equation}
    \mathrm{Im}\, (z^2\mu_0^z)= \mathrm{Im}\,(z^2) \left( 1-\theta+\sum_{j=1}^\infty \frac{\lambda_j^4  }{\big|z^2-\lambda^2_j\big|^2}  \left| \int_D \varphi_j(y)\,dy \right|^2 \right),
\end{equation}
thus, $ \mathrm{Im}\, (z^2\mu_0^z) $ has the same sign as $\mathrm{Im}\,(z^2)$. We conclude that $z^2\mu_0^z \in \mathbb{C} \setminus [0,\infty)$ for $\mathrm{Im}\,z>0$.

\emph{Step 2}. Since the spectrum of $L_0$ contains in $(-\infty,0]$, it follows from Step 1 and \eqref{resolvent_equlity} that $D_0(z)g =0$. Thus,
\begin{equation}
    Z_0 (z) g =0 \quad \mathrm{in} \ \{ \chi =1 \} , \qquad \mathrm{and} \qquad R_\Delta(z)g=0 \quad\mathrm{in} \ \{ \chi = 0 \} ,
\end{equation}
which implies that $  Z_0(z) g$ and $R_\Delta(z)g$ belong to $H^2(B_r)\cap H^1_0(B_r) $.
It follows from the definition \eqref{defF} of $F_0(z)$ that
\begin{equation}
    g = -F_0(z)g = z^2 \big\{ R_{\Delta } (z) -\mu_0^{z}Z_0(z)\big\}g = z^2R_{\Delta}(z)g  \qquad \mathrm{in}\  \Omega.
\end{equation}
Let $u=Z_0 (z)g - R_{\Delta}(z)g$, we get
\begin{equation}\label{1uap}
    \Delta  u =-\Delta R_{\Delta}(z)g = -g+z^2 R_{\Delta} g = 0  \qquad \mathrm{in}\ \Omega,
\end{equation}
and
\begin{equation}\label{2uap}
    \begin{aligned}
        \Delta u &= \Delta ( L^r_0)^{-1} \big(I- z^2 R_{\Delta}(z) \big) g - \Delta R_{\Delta}(z) g \\
        & =\big(I- z^2R_{\Delta}(z) \big)g -\Delta R_{\Delta}(z)g  \\
        &=0
    \end{aligned} \qquad \mathrm{in} \ B_r\setminus \overline{\Omega}.
\end{equation}
From these two equations \eqref{1uap}-\eqref{2uap} and the regularity $u\in H^2(B_r)$, we find that $u \in H_0^1(B_r)$ satisfies 
\begin{equation}
    \Delta  u=0 \qquad \mathrm{in}\ B_r.
\end{equation}
We get $u\equiv 0$ in $B_r$. Thus,
\begin{equation}
    R_{\Delta}(z)g  = D_0(z)g -\chi u = 0 \quad \mathrm{in} \ \mathbb{R}^d.
\end{equation}
Therefore, $g\equiv0$ in $\mathbb{R}^d$. This completes the proof.
\end{proof}

The following lemma is an analogue of Theorem \ref{meromorphic_continuation}, we omit the proof.

\begin{theorem}
    $(L_0+z^2\mu_0^z)^{-1}:L^2(\mathbb{R}^d) \rightarrow \mathcal{D}(L_0)$ is holomorphic for $\mathrm{Im}\,z>0$, and admits a meromorphic continuation to $\mathbb{C}^{\circ}$, as a family of operators from $L^2_{\mathrm{comp}}(\mathbb{R}^d)$ to $\mathcal{D}_{\mathrm{loc}}(L_0)$, where
    \begin{equation}
        \mathcal{D}_{\mathrm{loc}}(L_0) : = \big\{ u \in  H_{\mathrm{loc}}^1(\mathbb{R}^d): \rho \in C_0^{\infty}(\mathbb{R}^d), \rho|_{B_r} \equiv 1 \Rightarrow\nabla \cdot A_0 \nabla (\rho u) \in L^2(\mathbb{R}^d) \big\} .
    \end{equation}
This meromorphic continuation is called the outgoing resolvent of $L_0+z^2\mu_0^z$, and is denoted by $R_0(z)$. Moreover, we have the identity
\begin{equation}\label{Repdec2ap}
    R_0(z) = D_0(z) \big\{I + F_0(z) \rho\big\}^{-1} \big\{ I-F_0(z) (1-\rho) \big\},
\end{equation}
where $\rho \in C_0^{\infty} (\mathbb{R}^d; [0,1])$ is a cut-off such that $B_r \subset \{\rho=1\}$.
\end{theorem}

\begin{definition}\label{DefA2z}
The poles of $R_0(z)$ are called the scattering resonances associated with $L_0+z^2\mu_0^z$, denoted by $\Sigma_{\mathrm{hom}}$. A resonant state of $ L_0+z^2\mu_0^z $ is defined as a function $ u\in H^1_{\mathrm{loc}}(\R^d) $ such that
\begin{equation}
    u\in \Pi_{z}(L^2_{\mathrm{comp}}) \quad\ \mathrm{and} \quad\ (L_0+z^2\mu_0^z)u=0.
\end{equation}
Here 
\begin{equation}
    \Pi_z : = \frac{1}{2\pi \mathrm{i}}\oint_{z_0}  R_0(z)\frac{z^2\mu_0^z - z_0^2\mu_0^{z_0}}{z^2-z_0^2} 2z\,dz,
\end{equation}
and the integral $\oint_{z_0}$ denotes the contour integral that enclosing $z_0$ and no other poles of $R_0(z)$.
\end{definition}
 
\begin{proposition}\label{propSigmahomnega}
    $\Sigma_{\mathrm{hom}} \subset \{\mathrm{Im}\,z<0\}$.
\end{proposition}
\begin{proof}
    By Lemma \ref{InverIF}, $\Sigma_{\mathrm{hom}} \subset \{\mathrm{Im}\,z\leq 0\}$. For $ k=0 $, since $ D_0(0)=Q_0(0) $ and $ F_0(0)=K_0(0) $, following the proof of Theorem \ref{pole_K01}, we obtain that $0$ is not a pole of $R_0(z)$. It remains to show that, for any $k\in \mathbb{R} \setminus \{0\}$, $ k $ is not a pole of $(I+F_0(k))^{-1} $.
    
    Assume that there exists $ g\in L^2_{\mathrm{comp}}(\R^d) $ such that $ (I+F_0(k))g=0 $, then $ (L_0+k^2\mu^k_0)D_0(k)g=0 $. Denote $u=D_0(k)g$, multiply both sides by $ \overline{u} $ and integrate over $ B_R$, we get
    \begin{equation}
        \int_{B_R} A_0 \nabla u \cdot \overline{ \nabla u} - k^2 \mu_0^k |u|^2 \,dx = \int_{\partial B_R } \overline{u} \frac{\partial u}{\partial \mathbf{n} } \, dS.
    \end{equation}
    Since \( k\in\R \), the right-hand side is real. Thus, taking the imaginary parts of both sides, we get
\begin{equation}
    \mathrm{Im}\int_{\partial B_R } \overline{u} \frac{\partial u}{\partial \mathbf{n} } \, dS=0,
\end{equation}
which implies that \( u=0 \) in $ \mathbb{R}^d\setminus \overline{B_R}$ by Rellich's Uniqueness Theorem (Lemma \ref{rellichunique}). By unique continuation principle (\cite[Chapter 8]{colton1998inverse}), $ u=0 $ in $ \mathbb{R}^d\setminus \overline{\Omega}$. Hence, $u$ satisfies
\begin{equation}
    \left\{
    \begin{aligned} 
        &\nabla \cdot A_0 \nabla u + k^2 \mu_0^k u = 0 && \mathrm{in}\ \Omega,\\
        &u|_- = u|_+ =0 && \mathrm{on}\ \partial \Omega,\\
        &\frac{\partial u }{\partial \mathbf{n}} \Big|_-= \frac{\partial u }{\partial \mathbf{n}}\Big|_+ = 0&& \mathrm{on}\ \partial \Omega.
    \end{aligned}
    \right.
\end{equation}
This implies that $ u = D_0(k)g\equiv 0 $ in $\mathbb{R}^d$. Then
\begin{equation}
    Z_0(k)g=  D_0(k)g = 0 \quad \mathrm{in} \ \{\chi=1\}, \quad \mathrm{and}\qquad R_\Delta(k)g=0 \quad \mathrm{in}\ \{\chi=0\}.
\end{equation}
Thus, $ Z_0(k)g,\,R_\Delta(k)g \in H^1_0(B_r)\cap H^2(B_r)$. By definition of $ Z_0(k) $, $ Z_0(k)g-R_\Delta (k)g \in H^1_0(B_r) $ solves
\begin{equation}
    \Delta\big(Z_0(k)g-R_\Delta (k)g\big)=g-R_\Delta(k) g - \Delta R_\Delta(k) g =0 \quad \mathrm{in} \ B_r.
\end{equation}
Thus, $ Z_0(k)g-R_\Delta (k)g=0 $ in $ B_r $. We conclude that
\begin{equation}
    R_\Delta(k)g = D_0(k)g-\chi\big(Z_0(k)g-R_\Delta (k)g \big)=0 \qquad \mathrm{in}\ \mathbb{R}^d,
\end{equation}
and $ g\equiv0 $ in $\mathbb{R}^d$. The proof is complete.
\end{proof}

\begin{theorem}\label{thm_meroexpansion}
Let $z_0\in \mathbb{C}^{\circ} $ be a pole of $ R_0(z) $. Then $ R_0(z) $ admits the Laurent expansion near $ z_0 $:
    \begin{equation}\label{meroexpansion}
        R_0(z)= \sum_{k=1}^{K_0} \frac{A_k}{(z^2-z_0^2)^k} +A_0(z),
    \end{equation}
where $ A_0(z) $ is holomorphic near $z_0 $, $K_0$ is a positive integer that depends only on $z_0$, and $ \{ A_k \}_{k=1}^{K_0} $ is a family of finite-rank operators satisfying:  for any $ \rho\in C_0^\infty(\R^d) $ such that $ \rho|_{\Omega} \equiv 0$,
\begin{equation}
 A_k\rho= \Pi_{z_0}(-\Delta-z_0^2)^{k-1} \rho ,
\end{equation}
for any $1\le k\le K_0$, where we used the notation $ A_j=0 $ for $ j\ge K_0+1 $. Moreover,
\begin{equation}\label{meroexpansion_3}
        R_0(z)\rho= \Pi_{z_0} \sum_{k=1}^{K_0} \frac{ (-\Delta-z_0^2)^{k-1}}{(z^2-z_0^2)^k} \rho +A_0(z)\rho.
    \end{equation}
\end{theorem}
\begin{proof}
    Since $ R_0(z) $ is a meromorphic operator-valued family on $ \mathbb{C} $, there exist some finite-rank operators $ A_k , 1\le k\le K_0$ such that
    \begin{equation}\label{meroexpansion_4}
        R_0(z)= \sum_{k=1}^{K_0} \frac{A_k}{(z^2-z_0^2)^k} +A_0(z).
    \end{equation}
For any $ \rho\in C_0^\infty(\R^d) $, by analytic continuation, we have
    \begin{equation}
        R_0(z)(L_0+z^2\mu_0^z)\rho=\rho, \qquad z\in \mathbb{C}^\circ.
    \end{equation}
Inserting \eqref{meroexpansion_4} into the equality above, modulo terms that are holomorphic near $z$, and denoting $ f(z)=z^2\mu_0^z $, we have
\begin{equation}\label{meroexpansion_eq1}
    \begin{aligned}
        0 &\simeq \rho  \\
        & \simeq \sum_{k=1}^{K_0} \frac{A_k(L_0+z_0^2\mu^{z_0}_0 +z^2\mu_0^z -z_0^2 \mu_{0}^{z_0} )}{( z^2-z_0^2)^k} \rho\\
        & \simeq \sum_{k=1}^{K_0} \frac{1 }{(z^2-z_0^2)^k} \left(A_k( L_0 +z_0^2\mu_0^{z_0}) + A_{k+1} \mathbbm{1}_{\mathbb{R}^d\setminus \overline{\Omega}} + A_{k+1} \frac{f(z)-f(z_0)}{z^2-z_0^2}\mathbbm{1}_{ \Omega} \right) \rho,
    \end{aligned}
\end{equation}
where $\simeq$ denotes equality modulo terms that are holomorphic near $z_0$. It follows that for any $ \rho\in C_0^\infty(\R^d) $ such that $ \rho|_{\Omega} \equiv0$,
\begin{equation}
     A_{k+1}\rho= - A_k(\Delta+z^2)\rho,\quad1\le k<K_0, \quad \mathrm{and} \quad  A_{K_0}(\Delta+z^2)\rho = 0.
\end{equation}
In particular, noting that $ (\Delta+z^2)\rho\in C_0^\infty(\R^d)$ and $ (\Delta+z^2)\rho=0 $ on $\Omega$, we have
\begin{equation}
    A_k \rho = A_1(-\Delta-z^2)^{k-1}\rho, \quad \mathrm{and}\quad 0= A_{K_0}(-\Delta-z^2)\rho= A_1(-\Delta-z^2)^{K_0} \rho .
\end{equation}
Moreover, by a direct computation,
\begin{equation}
    \Pi_{z_0}\rho = \frac{1}{2 \pi \mathrm{i}}\oint_{z_0} R(z) \rho \,2z\,dz= A_1\rho,
\end{equation}
which proves \eqref{meroexpansion_3}. The proof is complete.
\end{proof}

\begin{theorem}\label{thmA3u}
Assume $ z_0\in \mathbb{C}^{\circ} $ . If $ u\in H^1_{\mathrm{loc}}(\R^d) $ solves $ (L_0+z_0^2\mu_0^{z_0})u=0 $ in $\R^d $ and there exist $ g\in L^2_{\mathrm{comp}}(\R^d) $ and $ R>0 $ such that
\begin{equation}
    u|_{\R^d\setminus B_R}= R_\Delta(z_0)g|_{\R^d\setminus B_R},
\end{equation}
then $ u$ is a resonant state.
\end{theorem}
\begin{proof}
    Assume $ u\in H^1_{\mathrm{loc}}(\R^d) $ solves $ (L_0+z^2\mu_0^z)u=0 $ in $\R^d $, and there exist $ g\in L^2_{\mathrm{comp}}(\R^d) $ and $ R>0 $ such that 
    \begin{equation}\label{resonanceeqiuv_1}
        u|_{\R^d\setminus B_R}=R_{\Delta}(z_0)g|_{\R^d\setminus B_R}.
    \end{equation}
    Applying $ L_0+z^2\mu_0^{z_0} $ on both sides, we see that $ \mathrm{supp}\,g\subset B_R $.

    By meromorphic continuation, we have, for any $ \rho \in C_0^\infty(\R^d) $, 
    \begin{equation}
        R_0(z)(L_0+z^2\mu_0^{z})\rho = \rho, \quad z\in \mathbb{C}^\circ.
    \end{equation}
Hence, using $ (L_0+z_0^2\mu_0^{z_0} )u=0 $, we get
\begin{equation}\label{resonanceequiv_2}
\begin{aligned}
    \rho u & = R_0(z )(L_0+z_0^2\mu_0^{z_0} + z^2\mu_0^{z} - z_0^2\mu_0^{z_0})\rho u\\
    &= R_0(z)[\Delta,\rho]  u + R_0(z)(z^2\mu_0^{z} - z_0^2\mu_0^{z_0} )\rho u.
\end{aligned}
\end{equation}

Taking cut-off function $ \rho_1\in C_0^\infty(\R^d;[0,1]) $, $ \rho_1\equiv 1 $ on $B_R$, let $ \rho $ be a cut-off function such that $ \rho\in C_0^\infty(\R^d;[0,1]) $, $ \rho\equiv 1 $ on $ \mathrm{supp}\, \rho_1 $, then
\begin{equation}
    (L_0+z^2\mu_0^{z})(1-\rho) R_{\Delta}(z)\rho_1 = -[\Delta,\rho]R_{\Delta}(z)\rho_1, \quad  \mathrm{Im}\,z>0,
\end{equation}
applying $ R_0(z) $ on both sides, we have, for $\mathrm{Im}\,z>0,$
\begin{equation}\label{resonanceequiv_4}
    (1-\rho) R_{\Delta}(z)\rho_1 =-R_0(z)[\Delta,\rho]R_{\Delta}(z)\rho_1 , 
\end{equation}
and, by meromorphic continuation, this holds for $ z\in\mathbb{C}^\circ $. Apply \eqref{resonanceequiv_4} on $ g$, since $\mathrm{supp}\,g \subset B_R$,
\begin{equation}\label{resonanceequiv_3}
\begin{aligned}
    (1-\rho) R_{\Delta}(z)g & =- R_0(z)[\Delta,\rho]\big( R_{\Delta}(z_0)g+ (R_{\Delta}(z)-R_{\Delta}(z_0))g\big)\\
    & = - R_0(z)[\Delta,\rho]u- R_0(z)[\Delta,\rho] (R_{\Delta}(z)-R_{\Delta}(z_0))g,
\end{aligned}
\end{equation}
where the second equality comes from \eqref{resonanceeqiuv_1} and the fact that $ \mathrm{supp}\, [\Delta,\rho] \subset \mathbb{R}^d\setminus B_R $.

Let
\begin{equation}
    u(z):= (1-\rho) R_{\Delta}(z)g +\rho u,
\end{equation}
    then $u(z_0)=u$, using \eqref{resonanceequiv_2} and \eqref{resonanceequiv_3}, we get
    \begin{equation}
        u(z)= R_0(z)(z^2\mu_0^{z} - z_0^2\mu_0^{z_0})\rho u- R_0(z)[\Delta,\rho] (R_{\Delta}(z)-R_{\Delta}(z_0))g.
    \end{equation}
Then,
\begin{equation}\label{B43}
    \begin{aligned}
        u&= u(z)\\
        &=\frac{1}{2\pi \mathrm{i}}\oint_{z_0} \frac{u(z)}{z^2-z_0^2} 2z \,dz\\
             & = \frac{1}{2\pi \mathrm{i}}\oint_{z_0} R_0(z) \frac{z^2\mu_0^{z} - z_0^2\mu_0^{z_0}}{z^2-z_0^2} 2z\,dz\, \rho u 
             - \frac{1}{2\pi \mathrm{i}}\oint_{z_0} R_0(z)[\Delta,\rho] \frac{R_{\Delta}(z)-R_{\Delta}(z_0) }{z^2-z_0^2} 2z \,dz\, g\\
             & =\Pi_{z_0}\rho u - \frac{1}{2\pi \mathrm{i}}\oint_{z_0} R_0(z)[\Delta,\rho] \frac{R_{\Delta}(z)-R_{\Delta}(z_0) }{z^2-z_0^2} 2z \,dz\, g.
    \end{aligned}
\end{equation} 

    Since $ \mathrm{supp}\, g \subset B_R $ and $ R_\Delta(z)g \in C^\infty(\mathbb{R}^d\setminus \overline{B_R}) $ for any $ z $, we get that
    \begin{equation}
        [\Delta,\rho] \frac{R_{\Delta}(z)-R_{\Delta}(z_0) }{z^2-z_0^2}  g \in C_0^\infty(\R^d) 
    \end{equation}
    vanishes on $ \Omega $. Then, By Theorem \ref{thm_meroexpansion}, the second term of the last line in \eqref{B43} becomes
        \begin{equation}
            \begin{aligned}
                &- \frac{1}{2\pi \mathrm{i}}\oint_{z_0} R_0(z)[\Delta,\rho] \frac{R_{\Delta}(z)-R_{\Delta}(z_0) }{z^2-z_0^2} 2z \,dz\, g \\
                &=\Pi_{z_0} \bigg\{-\sum_{k=1}^{K_0} (-\Delta-z_0^2)^{k-1}[\Delta,\rho] \frac{1}{2\pi \mathrm{i}}\oint_{z_0}\frac{R_{\Delta}(z)-R_{\Delta}(z_0) }{(z^2-z_0^2)^{k+1}} 2z\,dz \, g \bigg\}.
            \end{aligned}
        \end{equation}
Thus, $ u=\Pi_{z_0}v $ with
\begin{equation}
    v=\rho u - \sum_{k=1}^{K_0} (-\Delta-z_0^2)^{k-1}[\Delta,\rho] \frac{1}{2\pi \mathrm{i}}\oint_{z_0}\frac{R_{\Delta}(z)-R_{\Delta}(z_0) }{(z^2-z_0^2)^{k+1}} 2z \,dz\, g  ,
\end{equation}
which belongs to $L^2_{\mathrm{comp}}(\R^d)$. The proof is complete.
\end{proof}

\bibliographystyle{amsalpha}
\bibliography{reference}

\end{document}